\documentclass[12pt]{book}
\begin{document}
\sloppy
\pagestyle{myheadings}
\markboth{Meta Math!}{}
\pagenumbering{roman}

\section*{}

\textbf{\Huge
META MATH!\vspace{3mm}
\\
\Large
The Quest for Omega
\vspace{5mm}
\\
by Gregory Chaitin}

\section*{}

\emph{\textbf{Gregory Chaitin} has devoted his life to the attempt to understand what mathematics
can and cannot achieve, and is a member of the digital philosophy/digital physics movement.
Its members believe that
the world is built out of digital information, out of 0 and 1 bits,
and they view 
the universe as a giant information-processing machine, a giant digital computer.
In this book on the history of ideas, Chaitin
traces digital philosophy back to the 
nearly-forgotten
17th century genius Leibniz. 
He also tells us how he discovered the celebrated Omega number, which marks the current boundary
of what mathematics can achieve.
This book is an opportunity to get inside the head of a creative mathematician
and see what makes him tick,
and opens a window for its readers 
onto a glittering world of high-altitude thought that few intellectual mountain climbers
can ever glimpse.}

\chapter*{Cover}

William Blake: \emph{The Ancient of Days,} 1794.
\\
Relief etching with watercolor, $9\,\frac{1}{8} \times 6\,\frac{7}{8}$ inches.
\\
British Museum, London.

\chapter*{Preface}
\markright{Preface}

\section*{}
\section*{}
Science is an open road:
each question you answer raises \textbf{ten} new questions,
and much more difficult ones!  Yes, I'll tell you some things I discovered, but
the journey is endless, and mostly I'll share with you, the reader,
my doubts and preoccupations and what I think are promising and challenging
new things to think about.

It would be easy to spend many lifetimes working on any of a number
of the questions that I'll discuss.  That
is how the good questions are.  You can't answer them in five minutes,
and it would be no good if you could.

Science is an adventure.  I don't believe in spending years studying the work
of others, years learning a complicated field before I can contribute a tiny
little bit.  I prefer to stride off in totally new directions, where imagination is,
at least initially, much more important than technique, because the techniques
have yet to be developed.
It takes all kinds of people to advance knowledge, the pioneers, and those
who come afterwards and patiently work a farm.
This book is for pioneers!

Most books emphasize what the author knows. I'll try to emphasize what I would
like to know, what I'm hoping someone will discover, and just how much there is
that is fundamental and that we \textbf{don't know!}

And yes, I'm a mathematician, but I'm really interested in everything: what is
life, what's intelligence, what is consciousness, does the universe contain randomness, are space and
time continuous or discrete.  To me math is just the fundamental tool of
philosophy, it's a way to work out ideas, to flesh them out, to build models, \textbf{to understand}!
As Leibniz said, without math you cannot really understand philosophy, without philosophy
you cannot really understand mathematics, and with neither of them, you can't really
understand a thing!
Or at least that's my credo, that's how I operate.

On the other hand, as someone said a long time ago, ``a mathematician who is not something
of a poet will never be a good mathematician.'' And ``there is no permanent place in the world
for ugly mathematics'' (G. H. Hardy).  To survive, mathematical ideas must be beautiful, they must be
seductive, and they must be illuminating, they must help us to understand, they must inspire us.
So I hope this little book
will also convey something of this more personal aspect of mathematical creation,
of mathematics as a way of celebrating the universe, as a kind of love-making!
I want you to fall in love with mathematical ideas, to begin to feel seduced by them,
to see how easy it is to be entranced and to want to spend years in their company, 
years working on mathematical projects.

And it is a mistake to think that a mathematical idea can survive merely because it is \textbf{useful},
because it has practical applications.  On the contrary, what is useful varies as a function
of time, while ``a thing of beauty is a joy forever'' (Keats).
Deep theory is what is really useful, not the ephemeral usefulness of practical applications!

Part of that beauty, an essential part, is the clarity and sharpness that the mathematical
way of thinking about things promotes and achieves.  Yes, there are also mystic and poetic
ways of relating to the world, and to create a new math theory, or to discover new mathematics,
you have to feel comfortable with vague, unformed, embryonic ideas, even as you try to
sharpen them.  But one of the things about math that seduced me as a child was the black/whiteness,
the clarity and sharpness of the world of mathematical ideas, that is so different from the
messy (but wonderful!)\ world of human emotions and interpersonal complications!
No wonder that scientists express their understanding in mathematical terms, when they can!

As has been often said, to understand something is to make it mathematical, and I hope that
this may eventually even happen to the fields of psychology and sociology, someday.  That is my bias,
that the math point of view can contribute to everything, that it can help to clarify anything.
Mathematics is a way of characterizing or expressing \textbf{structure}.  And the universe
seems to be built, at some fundamental level, out of mathematical structure.  To speak metaphorically,
it appears that God is a mathematician, and that the structure of the world---God's thoughts!---are
mathematical, that this is the cloth out of which the world is woven, the wood out of which the world is built\ldots

When I was a child the excitement of relativity theory (Einstein!)\ 
and quantum mechanics was trailing off,
and the excitement of DNA and molecular biology had not begun.  What was the new big 
thing then? The Computer!  Which some people referred to then as ``Giant Electronic Brains.''
I was fascinated by computers as a child.  First of all, because they were a great toy, an infinitely
malleable artistic medium of creation.  
I loved programming!  But most of all, because the computer was (and still is!)\ a
wonderful new \textbf{philosophical} and mathematical concept.  
The computer is even more revolutionary as an \textbf{idea},
than it is as a practical device that alters society---and we all know 
how much it has changed our lives.  Why do I say this?  Well, the computer changes epistemology,
it changes the meaning of ``to understand.''  To me, you understand something only if you can
program it. (You, not someone else!)  Otherwise you don't really understand it, you only 
\textbf{think} you understand it.

And, as we shall see, the computer changes the way you do mathematics, it changes the kind
of mathematical models of the world that you build.  In a nutshell, now God seems to be a programmer,
not a mathematician! The computer has provoked a paradigm shift: it suggests a digital philosophy,
it suggests a new way of looking at the world,
in which everything is discrete and nothing is continuous,
in which everything is digital information, 0's and 1's.
So I was naturally attracted to this revolutionary new idea.

And what about the so-called real world of money and taxes and disease and death and war?
What about this ``best of all possible worlds, in which everything is a necessary evil!?''
Well, I prefer to ignore that insignificant world and concentrate instead on the world of
ideas, on the quest for understanding.  Instead of looking down into the mud, how about looking
up at the stars!  Why don't you try it? Maybe you'll like that kind of a life too!

Anyway, you don't have to read these musings from cover to cover.
Just leap right in wherever you like, and on a first reading 
please skip anything that seems too difficult.
Maybe it'll turn out afterwards that it's actually not so difficult\ldots\
I think that basic ideas are simple. And I'm not really interested in complicated ideas,
I'm only interested in fundamental ideas.
If the answer is extremely complicated, 
I think that probably means that we've asked the wrong question!

No man is an island, and practically every page of this book has benefited from
discussions with Fran\c{c}oise Chaitin-Chatelin during the past decade;
she has wanted me to write this book for that long.  
Gradiva, a Portuguese publishing house, arranged a stimulating visit
for me to that country in January 2004.
During that trip my lecture at the University of Lisbon on Chapter V
of this book was captured on digital video and is now available on the web.

I am also grateful to Jorge Aguirre for inviting me to present this book as a course
at his summer school at the University of R\'{\i}o Cuarto in C\'ordoba, Argentina, in
February 2004; the students' comments there were extremely helpful.  And Cristian Calude has
provided me with a delightful environment in which to finish this book at
his Center for Discrete Math and Theoretical Computer Science at the University
of Auckland. In particular, I thank Simona Dragomir for all her help.  

Finally, I am greatly indebted to Nabil Amer and to the Physics Department
at the IBM Thomas J. Watson Research Center in Yorktown Heights, New York 
(my home base between trips) for their support of my work.
\vspace{\baselineskip}
\\
\emph{---Gregory Chaitin}
\vspace{.5\baselineskip}
\\
Auckland, New Zealand, March 2004

\markright{}
\newpage
\section*{}
\section*{}
\section*{}
\section*{}
\begin{center}
\emph{\`A la femme de mes r\^eves, the eternal muse,\\ and the search for inexpressible beauty\ldots}
\end{center}

\tableofcontents

\chapter*{Quotes by Leibniz/Galileo}
\addcontentsline{toc}{chapter}{Quotes by Leibniz/Galileo}
\markright{Quotes by Leibniz/Galileo}
\pagenumbering{arabic}

\textbf{\large Leibniz ---}
\vspace{\baselineskip}
\\
Sans les math\'ematiques on ne p\'en\`etre point au fond de la philosophie.
\\
Sans la philosophie on ne p\'en\`etre point au fond des math\'ematiques.
\\
Sans les deux on ne p\'en\`etre au fond de rien. 
\vspace{.5\baselineskip}
\\
Without mathematics we cannot penetrate deeply into philosophy.
\\
Without philosophy we cannot penetrate deeply into mathematics.
\\
Without both we cannot penetrate deeply into anything.
\vspace{\baselineskip}
\\
\textbf{\large Galileo ---}
\vspace{\baselineskip}
\\
La filosofia \`e scritta in questo grandissimo libro 
\\
che continuamente ci sta aperto innanzi a gli occhi (io dico l'universo) 
\\
ma non si pu\`o intender se prima non s'impara a intender 
\\
la lingua e conoscere i caratteri ne' quali \`e scritto. 
\\
Egli \`e scritto in lingua matematica e 
\\
i caratteri sono triangoli, cerchi, ed altre figure geometriche 
\\
senza i quali mezi \`e impossibile a intenderne umanamente parola; 
\\
senza questi \`e un aggirarsi vanamente per un oscuro laberinto.
\vspace{.5\baselineskip}
\\
Philosophy is written in this very great book  
\\
which always lies open before our eyes (I mean the universe), 
\\
but one cannot understand it unless one first learns to understand 
\\
the language and recognize the characters in which it is written. 
\\
It is written in mathematical language and 
\\
the characters are triangles, circles and other geometrical figures;
\\
without these means it is humanly impossible to understand a word of it; 
\\
without these there is only clueless scrabbling around in a dark labyrinth.

\chapter*{Franz Kafka: Before the Law}
\addcontentsline{toc}{chapter}{Franz Kafka: Before the Law}
\markright{Franz Kafka: Before the Law}

Before the law sits a gatekeeper. To this
gatekeeper comes a man from the country who asks to gain entry into the law.
But the gatekeeper says that he cannot grant him entry at the moment.
The man thinks about it and then asks if he will be allowed to come in
later on. ``It is possible,'' says the gatekeeper, ``but not now.''
At the moment the gate to the law stands open, as always, and the
gatekeeper walks to the side, so the man bends over in order to see through the
gate into the inside. When the
gatekeeper notices that, he laughs and says: ``If it tempts you so much, try it
in spite of my prohibition. But
take note: I am powerful. And I am only the most lowly gatekeeper.
But from room to room stand gatekeepers, each more powerful than the
other. I can't even endure even one glimpse of the third.''
The man from the country has not expected such difficulties: the law
should always be accessible for everyone, he thinks, but as he now looks more
closely at the gatekeeper in his fur coat, at his large pointed nose and his
long, thin, black Tartar's beard, he decides that it would be better to wait
until he gets permission to go inside. The
gatekeeper gives him a stool and allows him to sit down at the side in front of
the gate. There he sits for days and years. He makes many attempts to be let in, and he wears the
gatekeeper out with his requests. The
gatekeeper often interrogates him briefly, questioning him about his homeland
and many other things, but they are indifferent questions, the kind great men
put, and at the end he always tells him once more that he cannot let him inside
yet. The man, who has equipped
himself with many things for his journey, spends everything, no matter how
valuable, to win over the gatekeeper. The
latter takes it all but, as he does so, says, ``I am taking this only so that
you do not think you have failed to do anything.''
During the many years the man observes the gatekeeper almost
continuously. He forgets the other
gatekeepers, and this one seems to him the only obstacle for entry into the law.
He curses the unlucky circumstance, in the first years thoughtlessly and
out loud, later, as he grows old, he still mumbles to himself.
He becomes childish and, since in the long years studying the gatekeeper
he has come to know the fleas in his fur collar, he even asks the fleas to help
him persuade the gatekeeper. Finally
his eyesight grows weak, and he does not know whether things are really darker
around him or whether his eyes are merely deceiving him.
But he recognizes now in the darkness an illumination which
breaks inextinguishably out of the gateway to the law.
Now he no longer has much time to live.
Before his death he gathers in his head all his experiences of the entire
time up into one question which he has not yet put to the gatekeeper.
He waves to him, since he can no longer lift up his stiffening body.
The gatekeeper has to bend way down to him, for the great
difference has changed things to the disadvantage of the man. ``What do you
still want to know, then?'' asks the gatekeeper. ``You are insatiable.''
``Everyone strives after the law,'' says the man, ``so how is that in
these many years no one except me has requested entry?''
The gatekeeper sees that the man is already dying and, in order to reach
his diminishing sense of hearing, he shouts at him, ``Here no one else can gain
entry, since this entrance was assigned only to you.
I'm going now to close it.''
\vspace{\baselineskip}
\\
{\footnotesize
[This translation is by Ian Johnston of Malaspina College, Canada. 
Note that in Hebrew ``Law'' is ``Torah'', which also means ``Truth''.
Orson Wells delivers a beautiful reading of this parable 
at the very beginning of his film
version of Kafka's \emph{The Trial}.]}

\chapter*{Chapter I---Introduction}
\addcontentsline{toc}{chapter}{Chapter I---Introduction}
\markright{Chapter I---Introduction}

In his book \emph{Everything and More: A Compact History of Infinity}, 
David Foster Wallace refers to G\"odel as
``modern math's absolute Prince of Darkness'' (p.\ 275) and states that 
because of him
``pure math's been in mid-air for the last 70 years'' (p.\ 284). 
In other words, according to Wallace, 
since G\"odel published his famous paper in 1931, mathematics has been
suspended hanging in mid-air without anything like a proper foundation.

It is high time 
these dark thoughts were permanently laid to rest.
Hilbert's century-old vision of a static completely mechanical 
absolutely rigorous formal mathematics, was a misguided attempt intended to
demonstrate the absolute certainty of mathematical reasoning.
It is time for us to recover from this disease!  

G\"odel's 1931 work on incompleteness, Turing's 1936 work on uncomputability, and
my own work on the role of information, randomness and complexity have shown
increasingly emphatically that the role that Hilbert envisioned for formalism
in mathematics is best served by computer programming languages, which \textbf{are} in fact formalisms
that can be mechanically interpreted---but they are formalisms for computing and calculating,
not for reasoning, not for proving theorems, and most emphatically not for inventing
new mathematical concepts nor for making new mathematical discoveries. 

In my opinion, the view that math provides absolute certainty and is static and perfect
while physics is tentative and constantly evolving is a false dichotomy.  Math is actually
not that different from physics.  Both are attempts of the human mind to organize, to make
sense, of human experience; in the case of physics, experience in the laboratory, in the
physical world, and in the case of math, experience in the computer, in the mental mindscape
of pure mathematics.

And mathematics is far from static and perfect; it is constantly evolving, constantly changing,
constantly morphing itself into new forms.  New concepts are constantly transforming math and
creating new fields, new viewpoints, new emphasis, and new questions to answer. 
And mathematicians do in fact utilize unproved new principles suggested by
computational experience, just as a physicist would.

And in discovering and creating new mathematics, mathematicians do base themselves on
intuition and inspiration, on unconscious motivations and impulses, and on their aesthetic sense, just
like any creative artist would.  And mathematicians do not lead logical mechanical ``rational''
lives.  Like any creative artist, they are passionate emotional people who deeply care
about their art, they are unconventional eccentrics motivated by mysterious forces, not by money
nor by a concern for the ``practical applications'' of their work.

I know, because I'm one of these crazy people myself!  I've been obsessed by these questions
for my whole life, starting at an early age.  And I'll give you an insider's view of all of this,
a first-hand report from the front, where there is still a lot of fighting, a lot of pushing
and shoving, between different viewpoints.
In fact basic questions like this are never settled, never definitively put aside, they have
a way of resurfacing, of popping up again in transformed form, every few generations\ldots

So that's what this book is about: 
It's about reasoning questioning itself, and its limits and the role of
creativity and intuition, and the sources of new ideas and of new knowledge.
That's a big subject, and I only understand a little bit of it, the areas that I've worked in
or experienced myself.
Some of this \textbf{nobody} understands very well, it's a task for the future. 
How about \textbf{you}?!  Maybe you can do some work in this area. 
Maybe you can push the darkness back a millimeter or two!
Maybe you can come up with an important new idea, maybe you can imagine a new kind of question
to ask, maybe you can transform the landscape by seeing it from a different point of view!
That's all it takes, just one little new idea, and lots and lots of hard work to develop it and to
convince other people!
Maybe you can put a scratch on the rock of eternity!

Remember that math is a free creation of the human mind, and as Cantor---the inventor
of the modern theory of infinity described by Wallace---said,
the essence of math resides in its freedom, in the freedom to create.
But history judges these creations by their enduring beauty and by the extent
to which they illuminate other mathematical ideas or the physical universe, in
a word, by their ``fertility''.
Just as the beauty of a woman's breasts or the delicious curve of her hips 
is actually concerned with child-bearing, and isn't merely for the delight of painters
and photographers, so a math idea's beauty also has something to do with its ``fertility'',
with the extent to which it enlightens us, illuminates us, and inspires us with other ideas and
suggests unsuspected connections and new viewpoints.

Anyone can define a new mathematical concept---many mathematical papers do---but only the beautiful
and the fertile ones survive. It's sort of similar to what Darwin termed ``sexual selection'',
which is the way in which animals (including us) choose their mates for their beauty. 
This is a part of Darwin's original theory of evolution
that nowadays one usually doesn't hear much about, but in my opinion it is 
much to be preferred to the ``survival of the
fittest'' and the ``nature red in tooth and claw'' view of biological evolution. As an example of this
unfortunate neglect, in a beautiful edition
of Darwin's \emph{The Descent of Man} that I happen to possess, several chapters on sexual selection 
have been completely eliminated!

So this gives some idea of the themes that I'll be exploring with you.
Now let me outline the contents of the book.

\section*{Overview of the Book}

Here's our road to $\Omega$:
\begin{itemize}
\item
In the second chapter I'll tell you how the idea of the computer entered mathematics
and quickly established its usefulness.
\item
In the third chapter, we'll add the idea of algorithmic information, of measuring
the size of computer programs.
\item
The intermezzo briefly discusses physical arguments against infinite-precision real numbers.
\item
The fifth chapter analyzes such numbers from a mathematical point of view.
\item
Finally the sixth chapter presents my information-based analysis of what mathematical
reasoning can or cannot achieve.
Here's $\Omega$ in all its glory.
\item
A brief concluding chapter discusses creativity\ldots
\item
And there's a short list of suggested books, plays, and even musicals!
\end{itemize}
Now let's get to work\ldots

\chapter*{Chapter II---Three Strange Loves: Primes/G\"odel/LISP}
\addcontentsline{toc}{chapter}{Chapter II---Three Strange Loves: Primes/G\"odel/LISP}
\markright{Chapter II---Three Strange Loves: Primes/G\"odel/LISP}

I am a mathematician, and this is a book about mathematics. So I'd like to start
by sharing with you my vision of mathematics: why is it beautiful, how it advances,
what fascinates me about it. And in order to do this I will give some case histories.
I think that it is useless to make general remarks about mathematics without exhibiting
some specific examples.  So you and I are going to do some real mathematics together:
important mathematics, significant mathematics.  I will try to make it as easy as
possible, but I'm going to show you the real thing!  If you can't understand something,
my advice is to just skim it to get the flavor of what is going on.  Then, if you are
interested, come back later and try to work your way through it slowly with a pencil and paper, 
looking at examples and special cases.
Or you can just skip all the
math and read the general observations about the nature of mathematics and the
mathematical enterprise that my examples are intended to illustrate!

In particular, this chapter will try to make the case the computer isn't just a billion
(or is it trillion?)\ dollar industry, it is also---which is more important to me---an
extremely significant and fundamental new concept that changes the way you think about
mathematical problems.  Of course, you can use a computer to check examples or to do
massive calculations, but I'm not talking about that.  I'm interested in the computer as a new
\textbf{idea}, a new and fundamental philosophical concept that changes mathematics, that solves
old problems better and suggests new problems, that changes our way of thinking and helps us to
understand things better, that gives us radically new insights\ldots

And I should say right away that I completely disagree with those who say that the field
of mathematics embodies static eternal perfection, and that mathematical ideas are inhuman
and unchanging.  On the contrary, these case studies, these intellectual histories, illustrate
the fact that mathematics is constantly evolving and changing, and that our perspective,
even on basic and deep mathematical questions, often shifts in amazing and unexpected fashion.\footnote
{See also the history of proofs that there are transcendental numbers in Chapter V.}
All it takes is a new idea!  You just have to be inspired, and then work like mad to develop
your new viewpoint. People will fight you at first, but if you are right, then everyone
will eventually say that it was \textbf{obviously} a better way to think about the problem, and
that you contributed little or nothing! In a way, that is the greatest compliment. And that's
exactly what happened to Galileo: he is a good example of this phenomenon in the history of ideas.  
The paradigm shift that he fought so hard to achieve is now taken so absolutely, totally and
thoroughly for granted, that we can no longer even understand how much he actually contributed!
We can no longer conceive of any other way of thinking about the problem.

And although mathematical ideas and thought are constantly evolving, you will also see
that the most basic fundamental problems never go away. Many of these problems go back to
the ancient Greeks, and maybe even to ancient Sumer, although we may never know for sure.
The fundamental philosophical questions like the continuous versus the discrete or the limits of knowledge
are \textbf{never} definitively solved. Each generation formulates its own answer,
strong personalities briefly impose their views, but the feeling of satisfaction is always temporary, 
and then the process continues, it continues forever.  Because in order to be able to fool
yourself into thinking that you have solved a really fundamental problem, you have to shut your
eyes and focus on only one tiny little aspect of the problem.
Okay, for a while you can do that, you can and should make progress that way.  
But after the brief elation of ``victory'', you, or other people who come after you, 
begin to realize that the
problem that you solved was only a toy version of the real problem, one that leaves out
significant aspects of the problem, 
aspects of the problem that in fact you \textbf{had} to ignore in order to be able to get anywhere. 
And those forgotten aspects of the problem never go away entirely: 
Instead they just wait patiently outside your cozy little mental construct, biding their time,
knowing that at some point someone else is going to have to take them seriously, even if it takes
hundreds of years for that to happen!

And finally, let me also say that I think that the history of ideas is the best way
to learn mathematics.  I always hated textbooks.  I always hated books full of formulas,
dry books with no colorful opinions, with no personality! The books I loved were books where
the author's personality shows through, books with lots of words, explanations and ideas, not
just formulas and equations!  I still think that the best way to learn a new idea is to see
its history, to see why someone was forced to go through the painful and wonderful process of
giving birth to a new idea!  To the person who discovered it, a new idea seems inevitable,
unavoidable.  The first paper may be clumsy, the first proof may not be polished, but that is raw
creation for you, just as messy as making love, just as messy as giving birth!
But you \textbf{will} be able to see where the new idea comes from.
If a proof is ``elegant'', if it's the result of two-hundred years of finicky polishing, 
it will be as inscrutable as a direct divine revelation, and it's impossible to guess how
anyone could have discovered or invented it. It will give you no insight, no, probably none at all.  

Enough talk! Let's begin! I'll have much more to say after we see a few examples.

\section*{An Example of the Beauty of Mathematics: The Study of the Prime Numbers}

The primes 
\[
   2, 3, 5, 7, 11, 13, 17, 19, 23, 29, 31, 37 \ldots 
\]
are the unsigned whole numbers with
no exact divisors except themselves and 1.  It is usually better not to consider 1 a prime,
for technical reasons (so that you get unique factorization into primes, see below).
So 2 is the only even prime and 
\[
   9 = 3 \times 3, \; 15 = 3 \times 5, \; 21 = 3 \times 7, \; 25 = 5 \times 5,  
\]
\[
   27 = 3 \times 9, \; 33 = 3 \times 11, \; 35 = 5 \times 7 \ldots
\]
are \textbf{not} primes.
If you keep factorizing a number, eventually you have to reach primes, which
can't be broken down any more. For example:
\[
   100 = 10 \times 10  = (2 \times 5) \times (2 \times 5) = 2^2 \times 5^2. 
\]
Alternatively, 
\[
   100 = 4 \times 25 = (2 \times 2) \times (5 \times 5) = 2^2 \times 5^2.
\]
Note that the final results are the same.  That this is always the case was already
demonstrated by Euclid two thousand years ago, 
but amazingly enough a much simpler proof was discovered recently (in the last century).
To make things easier in this chapter, though, 
I'm not going to use the fact that this is always the case, 
that factorization into primes is unique. So
we can go ahead without delay; we don't need to prove that.

The ancient Greeks came up with all these ideas two millennia ago, and they have fascinated
mathematicians ever since.
What is fascinating is that as simple as the whole numbers and the primes are, it is easy
to state clear, straight-forward questions about them that \textbf{nobody} knows how to answer, 
not even in two-thousand
years, not even the best mathematicians on earth!

Now I'd like to mention two ideas that we'll discuss a lot later: \textbf{irreducibility} and 
\textbf{randomness}.

Primes are \emph{irreducible numbers}, irreducible via multiplication, via factoring\ldots

And
what's mysterious about the primes is that they seem to be scattered about in a haphazard manner.
In other words,
the primes exhibit some kind of \emph{randomness}, since the local details of the distribution
of primes shows no apparent order,
even though we can calculate them one by one.

By the way,
Chapters 3 and 4 of Stephen Wolfram's \emph{A New Kind of Science} give many other examples 
of simple rules that yield extremely complicated behavior.
But the primes were the first example of this phenomenon that people noticed.

For example,
it is easy to find an arbitrarily long gap in the primes.
To find $N - 1$ non-prime numbers in a row---these are called composite numbers---just multiply
all the whole numbers from 1 to $N$ together (that's called $N!$, $N$ factorial), and add in turn 2, 3,
4 up through $N$ to this product $N!$:
\[
   N! + 2, \; N! + 3, \; \ldots \; N! + (N-1), \; N! + N.
\]
None of these numbers are prime.
In fact, the first one is divisible by 2, the second one is divisible by 3, and the last one
is divisible by $N$.
But you don't really need to go that far out. 

In fact the size of the gaps in the primes seem to jump around in a fairly haphazard manner too.
For example, it looks like there are infinitely many twin primes, consecutive odd primes
separated by just one even number.  The computational evidence is very persuasive. But no
one has managed to prove this yet!

And it is easy to show that there are infinitely many primes---we'll do this three different ways
below---but in whatever direction you go you quickly get to results which are conjectured, but
which no one knows how to prove.
So the frontiers of knowledge are nearby, in fact, extremely close.

For example, consider ``perfect'' numbers like 6, which is equal to the sum of all its proper
divisors---divisors less than the number itself---since $6 = 3 + 2 + 1$. 
Nobody knows if there are infinitely many perfect numbers.
It is known that each Mersenne prime, each prime of the form $2^n - 1$, gives you
a perfect number, namely $2^{n - 1} \times (2^n - 1)$, 
and that there are many Mersenne primes
and that every even perfect number is generated from a Mersenne prime in this way. But no one
knows if there are infinitely many Mersenne primes.
And nobody knows if there are \textbf{any} odd perfect numbers. Nobody has ever seen one, and
they would have to be very large, but no one is sure what's going on here\ldots\footnote
{For more on this subject, see Tobias Dantzig's beautiful history of ideas 
\emph{Number, The Language of Science}.}

So right away you get to the frontier, to questions that nobody knows how to answer.
Nevertheless many young math students, children and teenagers who have been bitten by the math 
bug, work away on these problems, hoping that they may succeed where everyone else has failed.
Which indeed they may! 
Or at least they may find something interesting along the way, 
even if they don't get all the way to their goal.
Sometimes a fresh look is better, sometimes it's better not to know what everyone else has done,
especially if they were all going in the wrong direction anyway!

For example, there is the recent discovery of a fast algorithm to check if a number
is prime by two students and their professor in India, a simple algorithm 
that had been missed by all the experts.
My friend Professor Jacob Schwartz at the Courant Institute of New York University 
even had the idea of including a few famous unsolved
math problems on his final exams, 
in the hope that a brilliant student who was not aware that they were
famous problems might in fact manage to solve one of them!

Yes, miracles can happen. But they don't happen that often.
And the question of whether we are burdened by our current knowledge is in fact a serious one.

Are the primes the right concept? How about perfect numbers?
A concept is only as good as the theorems that it leads to!  Perhaps we have been
following the wrong clues.  
How much of our current mathematics is habit, and how much is essential?
For example, instead of primes, perhaps we should 
be concerned with the opposite, with ``maximally divisible numbers''! 
In fact the brilliantly intuitive mathematician Ramanujan, who was self-taught,
and Doug Lenat's artificial intelligence program AM (for Automated Mathematician)
both came up with just such a concept.
Would mathematics done on another planet by 
intelligent
aliens be similar or very different from ours?

As the great French mathematician Henri Poincar\'e said, 
``Il y a les probl\`emes que l'on se pose, et les probl\`emes qui se posent.''
There are
questions that one asks, and questions that ask themselves!  So just how inevitable are our 
current concepts?
If evolution were rerun, would humans reappear?  If the history of math were rerun, would
the primes reappear?  Not sure!
Wolfram has examined this question in Chapter 12 of his book, and he comes up with 
interesting examples that suggest that our
current mathematics is much more arbitrary than most people think.

In fact the aliens are right here on this very planet!  Those amazing Australian animals, and the
mathematicians of previous centuries, they are the aliens, and they are very different from us.
Mathematical style and fashion varies substantially as a function of time, even 
in the comparatively recent past\ldots

Anyway, here comes our first case study of the history of ideas in mathematics.
Why are there infinitely many primes?

\section*{Euclid's Proof That There Are Infinitely Many Primes}

We'll prove that there are infinitely many primes by assuming that there are only
finitely many and deriving a contradiction. This is a common strategy in mathematical
proofs, and it's called a proof by \emph{reductio ad absurdum}, which is Latin for
``reduction to an absurdity''.

So let's suppose the opposite of what we wish to prove, namely that there are only finitely
many primes, and in fact that $K$ is the very last prime number. Now consider
\[
 1 + K! = 1 + (1 \times 2 \times 3 \times \cdots \times K).
\]
This is the product of all the positive integers up to 
what we've assumed is the very last prime, plus one.
But when this number is divided by any prime, it leaves the remainder 1!  
So it must itself be a prime!
Contradiction!  
So our initial assumption that $K$ was the largest prime has got to be false.

Here's another way to put it. Suppose that all the primes that we know are less than or equal
to $N$. How can we show that there has to be a bigger prime? Well, consider $N! + 1$.
Factorize it completely, until you get only primes.  Each of these primes has to be greater than $N$,
because no number $\leq N$ divides $N! + 1$ exactly. 
So the next prime greater than $N$ has to be $\leq N! + 1$.

This masterpiece of a proof has never been equalled, even though it is 2000 years old!
The balance between the ends and the means is quite remarkable.
And it shows that people were just as intelligent 2000 years ago as they are now.
Other, longer proofs, however, illuminate other aspects of the problem, 
and lead in other directions\ldots\
In fact, there are many interesting proofs that there are infinitely many primes.  Let me tell
you two more that I like.
You can skim over these two proofs if you wish, or skip them altogether. 
The important thing is to understand Euclid's original proof!

\section*{Euler's Proof That There Are Infinitely Many Primes}

\emph{Warning: this is the most difficult proof in this book.
Don't get stuck here.
In the future I won't give any long proofs,
I'll just explain the general idea.  But this is a really beautiful piece of mathematics by
a wonderful mathematician. It may not be your favorite way to prove that there are
infinitely many primes---it may be overkill---but it shows how very far you can get
in a few steps with what is essentially just high-school mathematics.
But I really want to encourage you to make the effort to understand
Euler's proof.  I am not a believer in instant gratification.  It is well worth
it to make a sustained effort to understand this one piece of mathematics; 
it's like long, slow foreplay rewarded by a final orgasm of understanding!}

First let's sum what's called an infinite geometric series:
\[
   1 + r + r^2 + r^3 + \cdots = 1/(1-r)
\]
This works as long as the absolute value (the unsigned value) of $r$ is less than 1.
\emph{Proof:}

\begin{center}
\begin{tabular}{|c|}
\hline
\\
Let $S_n$ stand for $1 + r + r^2 + r^3 + \cdots + r^n.$
\\
$S_n - (r \times S_n) = 1 - r^{n + 1}$
\\
\\
Therefore $(1 - r) \times S_n = 1 - r^{n + 1}$
\\
and $S_n = (1 - r^{n + 1})/(1 - r)$. 
\\
\\
So $S_n$ tends to $1/(1 - r)$ as $n$ goes to infinity,
\\
because if $-1 < r < 1$, then $r^{n + 1}$ goes to zero 
\\
as $n$ goes to infinity (gets bigger and bigger).
\\
\\
\hline
\end{tabular}
\end{center}

So we've just summed an infinite series! Let's check if the result makes any sense.
Well, let's look at the special case $r = 1/2$.  Then our result says that
\[
   1 + 1/2 + 1/4 + 1/8 + 1/16 + 1/32 + 1/64 + \ldots 1/2^n + \ldots 
\]
\[
   = 1/(1-r) = 1/(1/2) = 2,
\]
which is correct. What happens if $r$ is exactly zero? Then the infinite series
becomes $1 + 0 + 0 + 0 \ldots = 1$, which is $1/(1-r)$ with $r = 0$. And if $r$ is just a tiny bit greater
than 0, then the sum will be just a tiny bit greater than 1. And what if $r$ is just a tiny
little bit less than 1? Then the infinite series starts off looking like $1 + 1 + 1 + 1 \ldots$
and the sum will be very big, which also agrees with our formula: $1/(1-r)$ is also very large
if $r$ is just a tiny bit less than 1.  And for $r = 1$ everything falls apart, and coincidentally
both the infinite series and our expression $1/(1-1) = 1/0$ for the sum give infinity.
So that should give us some confidence in the result.

Now we'll sum the so-called ``harmonic series'' of all reciprocals of the positive integers.

\begin{center}
\begin{tabular}{|c|}
\hline
\\
$1 + 1/2 + 1/3 + 1/4 + 1/5 \ldots = \infty$
\\
\\
In other words, it diverges to infinity, 
it becomes arbitrarily 
\\
large
if you sum enough terms of the harmonic series.
\\
\\
\hline
\end{tabular}
\end{center}

\emph{Proof:} compare the harmonic series with $1 + 1/2 + 1/4 + 1/4 + 1/8 + 1/8 + 1/8 + 1/8 + \ldots$
Each term of the harmonic series is greater than or equal to the corresponding term in this new
series, which is obviously equal to $1 + 1/2 + 1/2 + 1/2 + \ldots$ and therefore diverges to infinity.

But the harmonic series is less than or equal to the product of 
\[
   1/(1-1/p) = 1 + 1/p + 1/p^2 + 1/p^3 + \ldots
\]
taken over all the primes $p$.
(This is the sum of a geometric series with ratio $r = 1/p$.)
Why?  Because the reciprocal of every number can be expressed as a product 
\[
 1/(p^{\alpha} \times q^{\beta} \times r^{\gamma} \times \ldots)
\]
of the reciprocal of prime powers $p^{\alpha}$, $q^{\beta}$, $r^{\gamma}$ \ldots

In other words,
\[
    1 + 1/2 + 1/3 + 1/4 + 1/5 + 1/6 + 1/7 + 1/8 + 1/9 + 1/10 + \ldots 
\]
\[
   \leq (1 + 1/2 + 1/4 + 1/8 + 1/16 + \ldots) \times (1 + 1/3 + 1/9 + 1/27 + 1/81 + \ldots)
\]
\[
   \times (1 + 1/5 + 1/25 + 1/125 + 1/625 + \ldots)
   \times \ldots
\]
\textbf{Since the left-hand side of this inequality diverges to infinity, the right-hand side must too,
so there must be infinitely many primes!}

This leads to Euler's product formula for the Riemann zeta function, as I'll now explain.

Actually, prime factorization is unique. (Exercise for budding mathematicians: 
Can you prove this by the
method of infinite descent?  Assume that $N$ is the smallest positive integer that
has two different prime factorizations, and show that a smaller positive
integer must also have this property.  But in fact 1, 2, 3, 4, 5 all have unique 
factorization, and if we keep going down, we must eventually get down to there!  Contradiction!)\footnote
{For a solution see Courant and Robbins, \emph{What is Mathematics?}}

So, a special case of what is known as Euler's product formula,
$1 + 1/2 + 1/3 + 1/4 + 1/5 \ldots =$ (not $\leq$)
product of $1/(1-1/p)$ over all primes $p$.
This is a special case of the following more general formula:
\[
   \zeta(s) = 1 + 1/2^s + 1/3^s + 1/4^s + 1/5^s \ldots = 
\]
\[
   \mbox{product over all primes $p$ of $1/(1-1/p^s)$},
\]
which gives us two different expressions for
$\zeta(s)$, which is Riemann's famous zeta function. Above we have been considering
$\zeta(1)$, which diverges to infinity.
The modern study of the statistical distribution of the primes depends on
delicate properties of Riemann's zeta function $\zeta(s)$ for complex arguments $s = a + b \sqrt{-1}$,
which is a complicated business that I'm not going to discuss, and is where the famous Riemann
hypothesis arises.

Let me just tell you how far I got playing with this myself 
when I was a teenager.
You can do a fair amount using comparatively elementary methods and the fact that
the sum of the reciprocals of the primes diverges:
\[
   1/2 + 1/3 + 1/5 + 1/7 + 1/11 + 1/13 + 1/17 + 1/23 + \ldots = \infty
\]
This was established by Euler, and it
shows that the primes cannot be too sparse, 
or this infinite series would converge to a finite sum instead of diverging to infinity.

\section*{My Complexity-Based Proof That There Are Infinitely Many Primes}

Because if there were only finitely many different primes, 
expressing a number $N$ via a prime factorization 
\[
   N = 2^e \times 3^f \times 5^g \times \ldots 
\]
would be too concise!
This is too compressed a form to give each number $N$.  There are too many $N$, and not enough
expressions that concise to name them all!

In other words, most $N$ cannot be defined that simply, they are too complex for that.
Of course, \textbf{some} numbers can be expressed extremely concisely: For example,
$2^{99999}$ is a very small expression for a very large number.  And
$2^{2^{99999}}$ is an even more dramatic example.  But these are \textbf{exceptions},
these are atypical.

In general, $N$ requires order of $\log N$ characters, but a prime factorization 
\[
   N = 2^e \times 3^f \times 5^g \times \ldots 
\]
with
a fixed number of primes would only require order of $\log \log N$ characters, and this isn't enough
characters to give that many, the required number, of different $N$!

If you think of
$2^e \times 3^f \times 5^g \times \ldots$ 
as a computer program for generating $N$, and if there were only finitely many primes,
these programs would be too small, they would enable you to compress all $N$ enormously, which
is impossible, because in general the best way to specify $N$ via a computer program is to just
give it explicitly as a constant in a program with \textbf{no} calculation at all!

For those of you who do not know what the ``log'' function is, you can think of it 
as ``the number of digits needed to write the number $N$''; that makes it sound less technical.
It grows by 1 each time that $N$ is multiplied by ten.
$\log \log N$ grows even more slowly.
It increases by 1 each time that $\log N$ is multiplied by ten.

Needless to say, I'll explain these ideas much better later. The size of
computer programs is one of
the major themes of this book. This is just our first taste of this new spice!

\section*{Discussion of These Three Proofs}

Okay, I've just sketched three different---very different---proofs that there are infinitely
many prime numbers.  One that's 2000 years old, one that's about 200 years old, and one that's
about 20 years old! Notice how very, very different these proofs are!

So right away, I believe that this completely explodes the myth, dear to both Bourbaki
and Paul Erd\"os, that there is only \textbf{one} perfect proof for each mathematical fact, just
one, the most elegant one.
Erd\"os used to refer to ``the book'', God's book with the perfect proof of each theorem.
His highest praise was, ``that's a proof from the book!'' As for Bourbaki, that enterprising
group of French mathematicians that liked to attribute the output of their collective efforts
to the fictitious ``Nicolas Bourbaki'', they would have endless fights and revisions of their monographs
until everything was absolutely perfect.  Only perfection was acceptable, nothing less than that!

In my opinion this is a totalitarian doctrine. Mathematical truth is not totally objective.
If a mathematical statement is false, there will be no proofs, but if it is true, there will
be an endless variety of proofs, not just one!  Proofs are not impersonal, they express
the personality of their creator/discoverer just as much as literary efforts do.  If something
important is true, there will be \textbf{many} reasons that it is true, many proofs of that fact.
Math is the music of reason, and some proofs sound like jazz, others sound like a fugue.
Which is better, the jazz or the fugues? Neither: it's all a matter of taste; some people prefer
jazz, some prefer fugues, and there's no arguing over individual tastes.
In fact this diversity is a good thing: if we all loved the same woman it would be a disaster!

And each proof will emphasize different aspects of the problem, each proof will lead in different
directions.  
Each one will have different corollaries, different generalizations\ldots\
Mathematical facts are not isolated, they are woven into a vast spider's web of
interconnections.  

As I said,
each proof will illuminate a different aspect of the problem. Nothing is ever 
absolutely black or white; things are always very complicated.  
Trivial questions may have a simple answer: $2 + 2$ is definitely not 5.
But if you are asking a \textbf{real} question, 
the answers are more likely to be: ``on the one hand this and that, on 
the other hand so and so'', even in the world of pure mathematics, 
not to mention the real world, which is 
much, much messier than the imaginary mental mindscape of pure mathematics.

Another thing about the primes and elementary number theory, is how close you always are
to the frontiers of knowledge. Yes, sometimes you are able to prove something nice, like our
three proofs that the list of primes is endless.  But those are the good questions, and
they are in a minority!  Most questions that you ask are extremely difficult or impossible to answer,
and even if you can answer them, the answers are extremely complicated and lead you nowhere.
In a way, math isn't the art of answering mathematical questions, it is the art of asking
the right questions, the questions that give you insight, the ones that lead you 
in interesting directions, the ones that connect with lots of other interesting questions---the
ones with beautiful answers!

And the map of our mathematical knowledge resembles a highway running through the desert or a
dangerous jungle;
if you stray off the road, you'll be hopelessly lost and die!  In other words, the current
map of mathematics reflects what our tools are currently able to handle, not what is really out there.
Mathematicians don't like to talk about what they don't know, they like to talk about the questions
that current technique, current mathematical technology, is capable of handling.  
Ph.D. students who are too
ambitious never finish their thesis and disappear from the profession, unfortunately.
And you may think that mathematical reality is objective, that it's not a matter of opinion. 
Supposedly
it is clear whether a proof is correct or not. Up to a point!  But whether a piece of math
is correct isn't enough, the real
question is whether it's ``interesting'', and that
is absolutely and totally a matter
of opinion, and one that depends on the current mathematical fashions.  So fields become popular,
and then they become unpopular, and then they disappear and are forgotten! Not always, but sometimes.
Only really important mathematical ideas survive.

Changing direction a bit, let me say why it's good to have many different proofs of an important
mathematical result.  Each, as I said before, illuminates a different aspect of the problem and
reveals different connections and leads you in different directions.  
But it is also a fact, as was said by the mathematician George
P\'olya in his lovely little book \emph{How to Solve It} (which I read as a child), 
that it's better
to stand on two legs than on one.  If a result is important, you badly want to find different
ways to see that fact; that's much safer.  If you have only one proof, and it contains an error,
then you're left with nothing.
Having several proofs is not only safer, it also gives you more insight, 
and it gives you more understanding.
After all, the real goal of mathematics is to obtain insight, not just proofs.
A long, 
complicated proof that gives you no insight is not only psychologically unsatisfying, it's fragile,
it may easily turn out to be flawed.
And I prefer proofs with ideas, not proofs with lots of computation.

\section*{My Love/Hate Relationship with G\"odel's Proof}

\begin{quote}
``The theory of numbers, more than any other branch of mathematics, began by being
an experimental science. Its most famous theorems have all been conjectured,
sometimes a hundred years or more before they were proved; and they have been
suggested by the evidence of a mass of computations.''\emph{---G. H. Hardy,} quoted in Dantzig op.\ cit.\
\end{quote}

So number theory is an experimental science as well as a mathematical theory.
But theory lags far, far behind experiment!
And one wonders, will it ever catch up?
Will the primes refuse to be tamed?
Will the mystery remain?
I was fascinated reading about all of this as a child.

And then one day I discovered that a little book had just been published.
It was by Nagel and Newman, and it was called \emph{G\"odel's Proof}.
This was in 1958, and the book was an expanded version of an article that
I'd also seen, and that was published by the two of them in \emph{Scientific American} in 1956.
It was love at first sight!  Mad love, crazy love, obsessive love, what the French
call ``amour \`a la folie.''
Here in fact was a possible explanation for the difficulties that mathematicians were
experiencing with the primes: G\"odel's incompleteness theorem, which asserts
that any finite system of mathematical axioms, any mathematical theory, is \textbf{incomplete}. 
More precisely, he showed that there will always be arithmetic assertions,
assertions about the positive integers and addition and multiplication, what are called
number-theoretic assertions, that are true but unprovable!

I carried this book around with me constantly, absolutely and totally fascinated, mesmerized
by the whole idea.  There was only one small, tiny little problem, fortunately, which was
that for the life of me I couldn't understand G\"odel's proof of this wonderful
meta-mathematical result. It's called that because it's not a mathematical result, it's
a theorem \textbf{about} mathematics itself, about the limitations of mathematical methods.
It's not a result within any field of mathematics, it stands outside looking down at mathematics,
which is itself a field called metamathematics!

I wasn't an idiot, so why couldn't I understand G\"odel's proof?  Well, I could follow
it step by step, but it was like trying to mix oil and water. My mind kept resisting.
In other words, I didn't lack the necessary intelligence, I just plain didn't like G\"odel's proof
of his fabulous result.  
His original proof seemed too complicated, too fragile!  It didn't seem to get to the heart
of the matter, because it was far from clear how prevalent incompleteness might in fact be.

And this is where my own career as a mathematician takes off.
I loved reading about number theory, I loved computer programming 
(for example, programs for calculating primes), I disliked G\"odel's
proof, but I loved Turing's alternative approach to incompleteness using the idea of the computer.
I felt very comfortable with computers. I thought they were a great toy, and I loved
writing, debugging and running computer programs---FORTRAN programs, 
machine language programs---I thought that this was great mental entertainment!
And at 15 I had the idea---anticipated by Leibniz in 1686---of looking at the
size of computer programs and of defining
a random string of bits to be one for which there is no program for calculating it that
is substantially smaller than it is.

The idea was define a kind of logical, mathematical or structural randomness, as opposed to the kind
of physical randomness that
Einstein and Infeld's delightful book \emph{The Evolution of Physics} emphasized to be
an essential characteristic of quantum physics, the physics of the microcosm.
That's another book suitable for children and teenagers that I can highly recommend,
and that I also worshiped at that stage of my life.

Let me explain what happened better; 
I'll now reveal to you one of the secrets of mathematical creation!
I loved incompleteness, but not G\"odel's proof.  Why?  Because of the lack of balance
between the ends and the means, between the theorem and its proof.
Such a deep and important---philosophically important---mathematical result deserved
a deep proof that would give deep insight into the ``why'' of incompleteness, instead of a clever
proof that only permitted you to have a superficial understanding of what was going on.
That was my feeling, totally on intuitive grounds, pure instinct, pure intuition,
my subconscious, gut-level, emotional reaction to G\"odel's proof.

And so I set to work to make it happen! This was a totally subjective act of creation, because
I \textbf{forced} it to happen. How?  Well, by changing the rules of the game, 
by reformulating the problem,
by redefining the context in which incompleteness was discussed in such a way that there would
\textbf{be} a deep reason for incompleteness, in such a way that a deeper reason for incompleteness
could emerge! You see, within the context that G\"odel
worked, he had done the best that was possible.  If you were to keep the setup exactly the
same as the one he had dealt with, there \textbf{was no} deeper reason for incompleteness.  And so
I proceeded to change the question until I could get out a deep reason for incompleteness.
My instinct was that the original context in which the problem of incompleteness was formulated
had to be changed to one that permitted such deeper understanding---that it was the wrong context
if this wasn't possible!

Now you see why I say that the mathematician is a creator as much as a discoverer, and why I
say that mathematical creation is a totally personal act.  

On the other hand, I couldn't have
succeeded if my intuitions that there was a deeper reason were incorrect. 
Mathematical truth is malleable, but only up to a point!

Another way to put it, is that I wanted to eliminate the superficial details that I thought
obscured these deeper truths, and I proceeded to change the formulation of the problem in order
to make this happen.  So you might say that this is an act of pure invention, that I created
a deeper reason for incompleteness because I so badly wanted there to be one.  That's true up
to a point. Another way to put it is that my intuition whispered to me because the idea \textbf{wanted}
to be found, because it was a more natural, a less forced way to perceive the question of 
incompleteness.  So, from that point of view, I wasn't \textbf{making} anything happen, on the contrary!
I simply was acutely sensitive to the vague half-formulated but more natural view of incompleteness
that was partially obscured by G\"odel's original formulation.

I think that both views of this 
particular
act of creation are correct: On the one hand, there was a \emph{masculine}
component, in making something happen by ignoring the community consensus of how to think about
the problem. On the other hand, there was a \emph{feminine} component, in allowing my hypersensitive
intuition to sense a delicate new truth that no one else was receptive to, that no one else was
listening for.

The purpose of this book is to explain what I created/discovered to you. 
It took many years of work, culminating with 
the halting probability $\Omega$---sometimes called Chaitin's number---that's
the discovery that I'm most proud of.
It will take several chapters for me to explain all of this well, 
because I'll have 
to build the appropriate intellectual
framework for thinking about incompleteness and my $\Omega$ number.

The first step is for me to explain to you the starting point for my own work, which most
definitely was not G\"odel's proof of 1931, but was instead Turing's alternative 1936 approach to
incompleteness, in which the idea of computation plays a fundamental role.
And it's thanks to Turing that the idea of computation, that the idea of the computer, became
a new force in mathematical thinking.  So let me tell you how this happened, and then I'll illustrate
the power of this new idea by showing you the way in which it solved an outstanding open 
question called
Hilbert's 10th problem, that was formulated by Hilbert as one of a list of 23 challenge 
problems in 1900.
That'll keep us busy for the rest of this chapter.

In the next chapter I'll get back to my fundamental new idea, to what I added to Turing,
which is my definition of randomness and complexity that Leibniz---who happens to be the inventor of
the infinitesimal calculus---clearly and definitely anticipated in 1686.

\section*{Hilbert, Turing \& Post on Formal Axiomatic Systems \& Incompleteness}

The first step in the direction of being able to use math to study the power of math
was taken by David Hilbert about a century ago.  
So I view him as the creator of metamathematics.
It was his idea that in order to be able to study
what mathematics can achieve, we first have to specify completely the rules of the game.
It was his idea to create a formal axiomatic system or FAS for all of mathematics, one that
would eliminate all the vagueness in mathematical argumentation,
one that would eliminate any doubt whether a mathematical proof is correct or not.

How is this done?  What's in a formal axiomatic system?
Well, the general idea is that it's like Euclid's \emph{Elements}, except that you have
to be much, much more fussy about all the details!

The first step is to create a completely formal artificial language for doing mathematics.
You specify the alphabet of symbols that you're using, the grammar, the axioms, the rules of
inference, and a proof-checking algorithm:

\begin{center}
\begin{tabular}{|c|}
\hline
\\
\textbf{\emph{\large Hilbert Formal Axiomatic System}}
\\
\\
Alphabet
\\
Grammar
\\
Axioms
\\
Rules of Inference
\\
Proof-Checking Algorithm
\\
\\
\hline
\end{tabular}
\end{center}

Mathematical proofs must be formulated in this language with \textbf{nothing} missing, with
every tiny step of the reasoning in place.
You start with the axioms, and then you apply the rules of inference one by one, and
you deduce all the theorems!
This was supposed to be like a computer programming language: So precise that a machine
can interpret it, so precise that a machine can understand it, so precise that a machine
can check it. 
No ambiguity, none at all!
No pronouns.  No spelling mistakes!  Perfect grammar!

And part of this package is that a finite set of mathematical axioms or postulates 
are explicitly given, plus you use symbolic logic to deduce all the possible consequences of the
axioms.  The axioms are the starting point for any mathematical theory; they are taken
as self-evident, without need for proof.  The consequences
of these axioms, and the consequences of the consequences, and the consequences of that, 
and so forth and so on \emph{ad infinitum,} are called
the ``theorems'' of the FAS.\footnote
{These ingredients are already in Euclid, except that some of his axioms are tacit,
some steps in his proofs are skipped, and he uses Greek instead of symbolic logic, in other words,
a human language, not a machine language.}

And a key element of a FAS is that there is a proof-checking algorithm, a mechanical
procedure for checking if a proof is correct or not.  In other words, there 
is a computer program that can decide whether or not a proof follows all the rules.
So if you are using a FAS to do mathematics, then you do not need to have
human referees check whether a mathematics paper is correct before you publish it.
You just run the computer program, and it tells you whether or not there's a mistake!

Up to a point this doesn't seem too much to demand: It's just the idea that math can
achieve perfect rigor, that mathematical truth is black or white, that math provides
absolute certainty.
We will see!  

\begin{center}
\begin{tabular}{|c|}
\hline
\\
\textbf{\emph{\large Math = Absolute Certainty???}}
\\
\\
\hline
\end{tabular}
\end{center}

The first step in this drama, in the decline and fall of mathematical certainty,
was an idea that I learned from Alan Turing's famous 1936 paper introducing the computer---as
a mathematical idea, not with actual hardware!  
And the most interesting thing about this paper, is Turing's
famous halting problem, which shows that there are things that no computer can ever
calculate, no matter how cleverly you program it, no matter how patient you are in waiting
for an answer.  In fact Turing finds two such things: the halting problem; 
and uncomputable real numbers, which we'll discuss in Chapter V.

\begin{center}
\begin{tabular}{|c|}
\hline
\\
\textbf{\emph{\large Turing, 1936:}} Halting Problem
\\
\\
\hline
\end{tabular}
\end{center}

Here I'll just talk about the halting problem.  What's the halting problem? It's the question
of whether or not a computer program, one that's entirely self-contained, one that does no
input/output, will ever stop.
If the program needs any numbers, they have to be given in the program itself, not read in from
the outside world.

So the program grinds away step by step, and either goes on forever, or else it eventually stops,
and the problem is to decide which is the case in a finite amount of time, without waiting
forever for it to stop.

And Turing was able to show the extremely fundamental result that there is no way to
decide in advance whether or not a computer program will ever halt, not in a finite amount of time.
If it does halt, you can eventually discover that. The problem is to decide when to give up
and decide that it will never halt. But there is no way to do that.

Just like I'm not polluting your mind with G\"odel's proof, I won't say a word about
how Turing showed that the halting problem cannot be solved, that there is no algorithm
for deciding if a computer program will never halt.

I'll give you my own proof later, in the Chapter VI. 
I'll show that you can't prove that a program is ``elegant'', by which I mean that it's the smallest 
possible program
for producing the output that it does.
From the fact that you can't establish elegance,
I'll deduce the immediate consequence that the halting problem must be unsolvable.
Not to worry, it'll all be explained in due course.

But what I will explain right now is how
Turing derived incompleteness from the halting problem.  

\begin{center}
\begin{tabular}{|c|}
\hline
\\
\textbf{\emph{\large Turing:}} Halting Problem implies Incompleteness!
\\
\\
\hline
\end{tabular}
\end{center}

The idea is very simple.
Let's assume that we have a FAS that can always \textbf{prove} whether or not individual programs
halt.  Then you just run through all possible proofs in size order until you 
find a proof that the particular program that you're interested in never halts, or you
find a proof that it does in fact halt.  You work your way systematically through
the tree of all possible proofs, starting from the axioms. Or you can write down one by one
all the possible strings of characters in your alphabet, in size order, and apply the
proof-checking algorithm to them to filter out the invalid proofs and determine all
the valid theorems.  In either case it's slow work, very slow, 
but what do I care, I'm a theoretician!
This isn't supposed to be a practical approach; 
it's more like what physicists call a thought experiment,
``gedanken'' experiment in the original German.
I'm trying to prove a theorem, 
I'm trying to show that a FAS must have certain limitations, I'm not trying to do anything practical.

So, Turing points out, if we have a FAS that can always prove whether or not individual
programs halt, and if the FAS is ``sound'' which means that all the theorems are true, then
we can generate the theorems of the FAS one by one and use this to decide if any
particular program that we are interested in will ever halt!  But that's impossible, it cannot be,
as Turing proved in his 1936 paper and I'll prove in a totally different manner in Chapter VI.

So the FAS must be incomplete. In other words, there must be programs for which there is
no proof that the program halts and there is also no proof that it will never halt.
There is no way to put all the truth and only the truth about the halting problem into a FAS!

Turing's marvelous idea was to introduce the notion of computability, of distinguishing things
that cannot be calculated from those that can, and then deduce incompleteness from
uncomputability.

\begin{center}
\begin{tabular}{|c|}
\hline
\\
\textbf{\emph{\large Turing:}} Uncomputability implies Incompleteness!
\\
\\
\hline
\end{tabular}
\end{center}

Uncomputability is a deeper reason for incompleteness. This immediately makes incompleteness
seem much more natural, because, as we shall see in Chapter V, a lot of things are uncomputable,
they're everywhere, they're very easy to find. 

(Actually, what Turing has shown is more general than that, it's that
\begin{center}
\begin{tabular}{|c|}
\hline
\\
Soundness + Completeness implies 
\\
you can systematically settle \textbf{anything} you can ask your FAS!
\\
\\
\hline
\end{tabular}
\end{center}
In fact, in some limited domains you \textbf{can} do that: Tarski did it for a big chunk of
Euclidean geometry.)

My own approach to incompleteness is similar to Turing's in that I deduce incompleteness from
something deeper, but in my case it'll be from randomness, not from uncomputability, as we shall see.
And random things are also everywhere, they're the rule, not the exception, as I'll explain
in Chapter V.

\begin{center}
\begin{tabular}{|c|}
\hline
\\
\textbf{\emph{\large My Approach:}} Randomness implies Incompleteness!
\\
\\
\hline
\end{tabular}
\end{center}

Anyway, the next major step on this path was taken in 1944 by Emil Post, who happens to have
been a professor at the City College of the City University of New York, which is where I was
a student when I wrote my first major paper on randomness.  They were very impressed by this
paper at City College, and they were nice enough to give me a gold medal, the Belden Mathematical
Prize, and later the Nehemiah Gitelson award for ``the search for truth''. 
It says that on the Gitelson medal in English and Hebrew, which is how I happen to know that
Truth and Law are both Torah in Hebrew; 
that's crucial in order to be able to understand
the Kafka parable at the beginning of this book.  
And as the extremely kind head of the math department, Professor Abraham Schwartz, 
was handing me the
Belden Prize gold medal, he pointed to the photos of former professors hanging on
the wall of his office, one of whom was Emil Post.

Post had the insight and depth to see the key idea in Turing's proof that uncomputability
implies incompleteness.  He extracted that idea from Turing's proof.  It was a jewel that
Turing himself did not sufficiently appreciate.  Post sensed that the essence of a FAS
is that there is an algorithm for generating all the theorems, one that's very slow and
that never stops, but that eventually gets to each and every one of them.  This amazing 
algorithm generates the theorems in size order, but not in order of the size of the statement
of each theorem, in order of the size of each proof.

\begin{center}
\begin{tabular}{|c|}
\hline
\\
\textbf{\emph{\large Hilbert/Turing/Post Formal Axiomatic System}}
\\
\\
Machine for generating all the theorems
\\
one by one, in some arbitrary order.
\\
\\
\hline
\end{tabular}
\end{center}

Of course, as I believe \'Emile Borel pointed out, it would be better to generate just
the interesting theorems, not all the theorems; most of them are totally uninteresting!
But no one knows how to do that. In fact, it's not even clear what an ``interesting''
theorem is.\footnote{Wolfram has some thoughts about what makes theorems interesting in his book.
As usual, he studies a large number of examples and extracts the interesting features.
That's his general \emph{modus operandi}.}

Anyway, it was Post who put his finger on the essential idea, which is the notion of
a set of objects that can be generated one by one by a machine, in some order, any order.
In other words, there is an algorithm, a computer program, for doing this.  
And that's the essential content
of the notion of a FAS, that's the toy model that I'll use to study the limits of the
formal axiomatic method.  I don't care about the details, all that matters to me is that
there's an algorithm for generating all the theorems, that's the key thing.

And since this is such an important notion, it would be nice to give it a name!
It used to be called an r.e.\ or ``recursively enumerable'' set.  The experts seem to
have switched to calling it a ``computably enumerable'' or c.e.\ set.  
I'm tempted to call it something
that Post said in one of his papers, which is that it's a ``generated set'',
one that can be generated by an algorithm, item by item, one by one, in some order or other,
slowly but surely. 
But I'll resist the temptation!

\begin{center}
\begin{tabular}{|c|}
\hline
\\
\textbf{\emph{\large FAS = c.e.\ set of mathematical assertions}}
\\
\\
\hline
\end{tabular}
\end{center}

So what's the bottom line? Well, it's this:
Hilbert's goal of \textbf{one} FAS for all of math
was an impossible goal, because you can't put all of mathematical truth into
just one FAS. Math can't be static, it's got to be dynamic, it's got to evolve.
You've got to keep extending your FAS, adding new principles, new ideas, new axioms,
just as if you were a physicist, without proving them, because they work!
Well, not exactly the way things are done in physics, but more in that spirit.
And this means that the idea of absolute certainty in mathematics becomes untenable.
Math and physics may be different,
but they're not that different, not as
different as people might think.
Neither of them gives you absolute certainty!

We'll discuss all of this at greater length in the concluding chapter, Chapter VII.
Here just let me add that I personally view metamathematics as a \emph{reductio ad absurdum}
of Hilbert's idea of a formal axiomatic system.  
It's an important idea \textbf{precisely because} it can be shot down!
And in the conflict between Hilbert and Poincar\'e over formalism versus intuition,
I am now definitely on the side of intuition.

And as you have no doubt noticed, here we've been doing ``philosophical'' math.
There are no long proofs.
All that counts are the ideas. You just take care of the ideas, and the proofs will take
care of themselves!
That's the kind of math I like to do!
It's also beautiful mathematics, just as beautiful as proving that there are infinitely many primes.
But it's a different kind of beauty!

Armed with these powerful new ideas, let's go back to number theory.
Let's see if we can figure out why number theory is so hard.

We'll use this notion of a computably enumerable set. 
We'll use it to analyze Hilbert's 10th problem, an important problem in number theory.

More precisely, I'll explain what Hilbert's 10th problem is, and then I'll tell you
how Yuri Matiyasevich, Martin Davis, Hilary Putnam and Julia Robinson 
managed to show that it \textbf{can't} be solved. 

Davis, Putnam and Robinson did much of the preparatory work, and then the final steps
were taken by Matiyasevich. Some important additional results were then
obtained much more simply by Matiyasevich together with James Jones. 
Matiyasevich and Jones built on a curious and unappreciated piece of work by \'Edouard Lucas,
a remarkable and extremely talented and unconventional French mathematician of the late 1800's 
who never let fashion stand in his way.
Lucas also invented an extremely fast algorithm for deciding if a Mersenne number
$2^n - 1$ is a prime or not, which is how the largest primes that we currently
know were all discovered. 

By the way, Martin Davis studied with Emil Post at City College. It's a small world!

\section*{Hilbert's 10th Problem \& Diophantine Equations as Computers}

What is a diophantine equation? Well, it's an algebraic equation in which everything, the
constants as well as the unknowns, has got to be an integer.
And what is Hilbert's 10th problem? It's Hilbert's challenge to the world to discover an
algorithm for determining whether a diophantine equation can be solved.  In other words,
an algorithm to decide if there is a bunch of whole numbers that you can plug into 
the unknowns and satisfy the equation.

Note that if there \textbf{is} a solution, it can eventually be found by systematically plugging into
the unknowns all the possible whole-number values, starting with small numbers and gradually
working your way up.  The problem is to decide when to give up. Doesn't that sound familiar?
Yes, that's exactly the same as what happens with Turing's halting problem!

Diophantus was, like Euclid, a Greek scholar in Alexandria, where the famous library was.
And he wrote a book about diophantine equations
that inspired Fermat, who read a Latin translation of the
Greek original.  Fermat's famous ``last theorem'', recently proved by Andrew Wiles, was a marginal
annotation in Fermat's copy of Diophantus.

More precisely, we'll work only with unsigned integers, namely the positive integers and zero:
0, 1, 2, 3, \ldots

For example, $3 \times x \times y = 12$ has the solution $x = y = 2$, 
and $x^2 = 7$ has no whole-number solution.
 
And we'll allow an equation to be assembled from constants and unknowns using only additions
and multiplications, which is called an (ordinary) diophantine equation, or using additions,
multiplications and exponentiations, which is called an exponential diophantine equation.
The idea is to avoid negative integers and minus signs and subtraction by actually using both sides
of the equation. 
So you \textbf{cannot} collect 
everything on the left-hand side and reduce the right-hand side to zero.
 
In an ordinary diophantine equation the exponents will always be constants.
In an exponential diophantine equation, exponents may also contain unknowns.
Fermat-Wiles states that the exponential diophantine equation
\[
   x^n + y^n = z^n
\] 
has no solutions with $x, y$ and $z$ greater than 0 and $n$ greater than 2.
It only took a few hundred years to find a proof!  Too bad that Fermat
only had enough room in the margin of Diophantus to state his result,
but not enough room to say \textbf{how} he did it.
And all of Fermat's personal papers have disappeared, so there is nowhere left to search for clues.
 
So that's just how very bad an exponential diophantine problem can be!
Hilbert wasn't that ambitious, he was only asking for a way to decide if
\textbf{ordinary} diophantine equations have solutions.  That's bad enough!
 
To illustrate these ideas,
let's suppose we have two equations, $u = v$ and $x = y$, and we want to combine them into a single
equation that has precisely the same solutions.  Well, if we could use subtraction, the trick
would be to combine them into
\[
   (u - v)^2 + (x - y)^2 = 0.
\] 
The minus signs can be avoided by converting this into
\[
   u^2 - (2 \times u \times v) + v^2 \;\; + \;\;
   x^2 - (2 \times x \times y) + y^2 \;\; = \;\; 0
\] 
and then to
\[
   u^2 + v^2 \; + \;
   x^2 + y^2 \; = \; 
   (2 \times u \times v) \; + \;
   (2 \times x \times y).
\] 
That does the trick!
 
Okay, now we're ready to state the amazing Matiyasevich/Davis/Put\-nam/Robinson
solution of Hilbert's 10th problem. Diophantine equations come, as we've seen,
from classical Alexandria.  But the solution to Hilbert's problem
is amazingly modern, and could not even have been imagined by Hilbert, since 
the relevant concepts did not yet exist.
 
Here is how you show that there is no solution, that there is no algorithm for
determining whether or not a diophantine equation can be solved, that there will
never be one.
 
It turns out that there is a diophantine equation that is a computer.
It's actually called a universal diophantine equation, because it's like
a universal Turing machine, math-speak for a general-purpose computer, one
that can be programmed to run any algorithm, not just perform special-purpose calculations
like some pocket calculators do.
What do I mean?
Well, here's the equation:
\begin{center}
\begin{tabular}{|c|}
\hline
\\
\textbf{\emph{\large Diophantine Equation Computer:}}
\\
\\
$L(k, n, x, y, z, \ldots) = R(k, n, x, y, z, \ldots)$
\\
\\
\textbf{Program} $k$
\\
\textbf{Output} $n$
\\
\textbf{Time} $x, y, z, \ldots$
\\
\\
\hline
\end{tabular}
\end{center}
(The real thing is too big to write down!)
There's a left-hand side $L$, a right-hand side $R$, a lot of unknowns $x, y, z, \ldots$ and
two special symbols, $k$, called the parameter of the equation, and $n$, which is the unknown
that we really care about. \textbf{$k$ is the program for the computer, $n$ is the output that it
produces, and $x, y, z, \ldots$ are a multi-dimensional time variable!}  In other words,
you put into the equation $L(k, n) = R(k, n)$ a specific value for $k$ which is the computer
program that you are going to run.  Then you focus on the values of $n$ for which there are
solutions to the equation, that is, for which you can find a time $x, y, z, \ldots$ at which the
program $k$ outputs $n$.
That's it! That's all there is to it!
 
So, amazingly enough, this \textbf{one} diophantine equation can perform \textbf{any} calculation.
In particular, you can put in a $k$ that's a program for calculating all the primes, or you
can put in a $k$ that's a program for computably enumerating all the theorems of a particular FAS.
The theorems will come out as numbers, not character strings, but you can just convert the
number to binary, omit the leftmost 1 bit, divide it up into groups of 8 bits, look up each
ASCII character code, and presto chango, you have the theorem!
 
And this diophantine equation means we're in trouble. It means that there \textbf{cannot}
be a way to decide if a diophantine equation has any solutions. Because if we could do that,
we could decide if a computer program gives any output, and if we could do that, we could
solve Turing's halting problem. 

\begin{center}
\begin{tabular}{|c|}
\hline
\\
\textbf{\emph{\large Halting Problem unsolvable implies}}
\\
\textbf{\emph{\large Hilbert's 10th Problem is unsolvable!}}
\\
\\
\hline
\end{tabular}
\end{center}

Why does being able to decide if a computer program has
any output enable us to decide whether or not a program halts? Well, any program that
either halts or doesn't and produces no output, can be converted into one that produces a message
``I'm about to halt!''\ (in binary, as an integer) just before it does.
If it's written in a high-level language, that should be easy to do.
If it's written in a machine language, you just run the program interpretively, not directly,
and catch it just before it halts.
So if you can decide whether the converted program produces any output, you can decide whether the
original program halts.
 
That's all there is to it: diophantine equations are actually computers!
Amazing, isn't it!  What wouldn't I give to be able to explain this to Diophantus and to Hilbert!
I'll bet they would understand right away!
(This book is how I'd do it.)
 
\textbf{And this \emph{\large proves} that number theory is hard! 
This proves that uncomputability and
incompleteness are lurking right at the core, in two thousand-year old diophantine problems!
Later I'll show that randomness is also hiding there\ldots}
 
And how do you actually construct this monster diophantine equation computer? Well, it's a big
job, rather like designing actual computer hardware.  In fact, I once was involved in helping
out on the design of an IBM computer. (I also worked on the operating system and on a compiler,
which was a lot of fun.)  And the job of constructing this diophantine equation reminds me
of what's called the logical design of a CPU.  It's sort of an algebraic version of the design of
a CPU.
 
The original proof by Matiyasevich/Davis/Putnam/Robinson is complicated, and designs an
ordinary diophantine equation computer. The subsequent Matiyasevich/Jones design for
an exponential diophantine equation computer is much, much simpler---fortunately!---so 
I was able to actually program it out and exhibit the equation, 
which nobody else had ever bothered to do. 
I'll tell you about that in the next section, the section on LISP.
 
And 
in designing this computer,
you don't work with individual bits, you work with strings of bits represented in the form
of integers, so you can process a lot of bits at a time.  And the key step was provided by Lucas
a century ago and is a way to compare two bit strings that are the same size
and ensure that every time that a bit is on
in the first string, the corresponding bit is also on in the second string.  Matiyasevich and
Jones show that that's enough, 
that it can be expressed via diophantine equations, and
that it's much more useful to be able to do it than you might think.
In fact, that's all the number theory that Matiyasevich and Jones need to build their
computer; the rest of the job is, so to speak, computer engineering, not number theory.
From that point on,
it's basically just a lot of high-school math!
 
So what's this great idea contributed posthumously by Lucas? Well, it involves ``binomial coefficients''.
Here are some.  They're what you get when you expand powers of $(1 + x)$:
 
\begin{center}
\begin{tabular}{|c|}
\hline
\\
\textbf{\emph{\large Binomial Coefficients:}}
\\
\\
$(1+x)^0 = 1$
\\
$(1+x)^1 = 1 + x$
\\
$(1+x)^2 = 1 + 2x + x^2$
\\
$(1+x)^3 = 1 + 3x + 3x^2 + x^3$
\\
$(1+x)^4 = 1 + 4x + 6x^2 + 4x^3 + x^4$
\\
$(1+x)^5 = 1 + 5x + 10x^2 + 10x^3 + 5x^4 + x^5$
\\
\\
\hline
\end{tabular}
\end{center}

Let's consider the coefficient of $x^k$ in the expansion of $(1+x)^n$.
Lucas's amazing result is that this binomial coefficient is odd if and only if each time that
a bit is on in the number $k$, the corresponding bit is also on in the number $n$! 
 
Let's check some examples. Well, if $k = 0$, then every binomial coefficient is 1 which is odd, no matter what
$n$ is, which is fine, because every bit in $k$ is off.
And if $k = n$, then the binomial coefficient is also 1, which is fine.
How about $n = 5 =$ ``101'' in binary, and $k = 1 =$ ``001''. Then the binomial coefficient is 5 which is odd,
which is correct.   What if we change $k$ to $2 =$ ``010''. Aha! 
Then each bit that's on in $k$ is \textbf{not} on in $n$,
and the binomial coefficient is 10, which is \textbf{even}!
Perfect!
 
Please note that
if the two bit strings aren't the same size, then you have to add 0's \textbf{on the left} to the
smaller one in order to make them the same size.  Then you can look at the corresponding bits.
 
{\footnotesize [\emph{Exercise for budding mathematicians:} Can you check more values? By hand? On a computer?
Can you prove that Lucas was right?!  It shouldn't be too difficult to convince yourself.
As usual, it's much more difficult to \textbf{discover} a new result, 
to \textbf{imagine} the happy possibility,
than to verify that it's correct. Imagination, inspiration and hope are the key!
Verification is routine. \textbf{Anyone} can do that, any competent professional.
I prefer tormented, passionate searchers for Truth!
People who have been seized by a demon!---hopefully a \textbf{good} demon, but a demon nevertheless!]}
 
And how do Matiyasevich and Jones express (if a bit is on in $k$ it also
has to be on in $n$) as a diophantine equation?
 
Well, I'll just show you how they did it, and you can see why it works
with the help of the hints that I'll give.
 
\begin{center}
\begin{tabular}{|c|}
\hline
\\
The coefficient of $x^K$ in the expansion 
\\
of $(1+x)^N$ is odd if and only if there is 
\\
a \textbf{unique} set of seven unsigned whole 
\\
numbers $b, x, y, z, u, v, w$ such that
\\ \\
$b = 2^N$
\\ \\
$(b+1)^N = x b^{K+1} + y b^K + z$ 
\\ \\
$z + u + 1 = b^K$ 
\\
$y + v + 1 = b$
\\
$y = 2w + 1$
\\
\\
\hline
\end{tabular}
\end{center}

Can you figure out how this works? 
\emph{Hint:} $y$ is the binomial coefficient that we're interested in, 
$w$ is used to ensure that $y$ is odd,
$u$ guarantees that $z$ is smaller than $b^k$, $v$ guarantees that $y$ is smaller than $b$,
and when $(b+1)^n$ is written as a base $b$ numeral its
successive digits are precisely the binomial coefficients that you get by expanding
$(1+x)^n$.
The basic idea here is that the powers of the number eleven
give you the binomial coefficients: $11^2 = 121$, $11^3 = 1331$,
$11^4 = 14641$. Then things fall apart because the binomial coefficients get too large
to be individual digits. So you have to work in a large base $b$.  The particular $b$ that
Matiyasevich and Jones pick happens to be the \textbf{sum} of all the binomial coefficients that  
we are trying to pack together into the digits of a single number.
So the base $b$ is larger than all of them,
and there are no carries from one digit to the next to mess things up,
so everything works!
 
And then you combine these five individual equations 
into a \textbf{single} equation the way I showed you before.
There will be ten squares on one side, and five cross-products on the other side\ldots
 
Now that we've done this, we have the fundamental tool that we need in order to be able
to build our computer.
The CPU of that computer is going to have a bunch of hardware registers, just like real
computers do.  But unlike real computers, each register will hold an unsigned integer,
which can be an arbitrarily large number of bits.
And the unknowns in our master equation will give the contents of these registers.
Each digit of one of these unknowns will be the number that that register
happened to have at a particular moment, step by step, as the CPU clock ticks away and
it executes one instruction per clock cycle.
In other words, you have to use the same trick of packing \textbf{a list} of numbers into
a single number that we used before with binomial coefficients. Here it's done with the time
history of CPU register contents, 
a register contents history that ranges
from the start of the computation to the moment that you output an integer.
 
I'm afraid that we're going to have to stop at this point.
After all, constructing a computer is a big job!
And I'm personally much more interested in general ideas than in details.
You just take care of the general ideas, and the details will take care of themselves!
Well, not always: You wouldn't want to buy a computer that was designed that way!
But carrying out this construction in detail is a big job\ldots\
Not much fun! Hard work!  
 
Actually, it \textbf{is} a lot of fun if you do it yourself. 
Like making love,
it's not much fun to hear how somebody else did it.
You have to do it yourself!
And you don't learn much by reading someone else's software. It's excruciating to have to do that,
you're lost in the maze of someone else's train of thought!
But if you yourself write the program and debug it on interesting examples, 
then it's \textbf{your} train of thought and you learn a lot.
Software is frozen thought\ldots
 
So we won't finish working out the design in detail.
But the basic ideas are simple, and you can find them in the short
paper by Matiyasevich and Jones, or 
in a book
I wrote that was published by Cambridge University Press in 1987. It's called \emph{Algorithmic
Information Theory}, and that's where I actually exhibit one of these computer equations.  
It's an equation that runs LISP programs.
 
What kind of a programming language is LISP?
Well, it's a very beautiful language that's highly mathematical in spirit.

\section*{LISP, an Elegant Programming Language}

Utopian social theorists and philosophers have often dreamt of a perfect, ideal
universal human language that would be unambiguous and would promote human understanding
and international cooperation: for example, Esperanto, Leibniz's
\emph{characteristica universalis,} Hilbert's formal axiomatic system for all of mathematics.
These utopian fantasies have never worked.
But they have worked for computers---very, very well indeed!
 
Unfortunately,
as programming languages become increasingly sophisticated, they reflect more and
more the complexity of human society and of the immense world of software applications. 
So they become more and more like giant
tool boxes, like garages and attics stuffed with 30 years of belongings!
On the contrary, LISP is a programming language with considerable
mathematical beauty; it is more like a surgeon's scalpel or a sharp-edged diamond cutting tool
than a two-car garage overflowing with possessions and absolutely no room for a car.  
 
LISP has a few powerful simple basic concepts,
and everything else is built up from that, which is how mathematicians like
to work; it's what mathematical theories look like.
Mathematical theories, the good ones, consist in defining a few new 
key 
concepts,
and then the fireworks begin: they reveal new vistas, they open the door
to entirely new worlds. LISP is like that too; it's more like pure math than most
programming languages are.  At least it is if you strip away the ``useful''
parts that have been added on, the accretions that have made
LISP into a ``practical tool''.   What's left if you do that is the original LISP,
the conceptual heart of LISP, a core which is a jewel of considerable mathematical
beauty, austere intellectual beauty. 
 
So you see, for me this is a real love affair.              
And how did I fall in love with LISP?
In 1970 
I was living in Buenos Aires, and I bought a LISP manual
in Rio de Janeiro.
On page 14 there was a LISP interpreter written in LISP.
I didn't understand a word!
It looked extremely strange.
I programmed the interpreter in FORTRAN and suddenly saw the light.
It was devastatingly beautiful!
I was hooked! (I still am. I've used LISP in six of my books.)
I ran lots of LISP programs on my interpreter.
I wrote dozens of LISP interpreters, probably hundreds of LISP interpreters!
I kept inventing new LISP dialects, one after another\ldots
 
You see, learning a radically
new programming language is a lot of work.
It's like learning a foreign language: you need to pick up the vocabulary,
you need to acquire the culture, 
the foreign world-view, that inevitably comes with a language.
Each culture has a different way of looking at the world!
In the case of a programming language,
that's called its ``programming paradigm'', and to acquire a new programming paradigm
you need to read the manual, you need to see a lot of examples,
and you need to try writing and running programs yourself.
We can't do all that here, no way!  So the idea is just to whet
your appetite and try to suggest how very different LISP is from normal
programming languages.

\section*{How does LISP work?}

First of all, LISP is a non-numerical programming language. 
Instead of computing only with numbers,
you deal
with symbolic expressions, which are called ``S-expressions''.
 
Programs and data in LISP are both S-expressions.
What's an S-expression? Well, it's like an algebraic expression
with full parenthesization. For example, instead of
\[
 a \times b + c \times d
\] 
you write
\[
 ((a \times b) + (c \times d))
\] 
And then you move operators forward, 
in front of their arguments rather than between their arguments. 
\[
 (+ \; (\times \; a \; b) \; (\times \; c \; d))
\] 
That's called prefix as opposed to infix notation.
Actually in LISP it's written this way
\[\mbox{\texttt{
 (+ (* a b) (* c d))
}}\]
because, as in many other programming languages, \texttt{*} is used for multiplication. 
You also have minus \texttt{-} and exponentiation \verb|^| in LISP.
 
Let's get back to S-expressions. In general, an S-expression consists of a nest
of parentheses that balance, like this
\[\mbox{\texttt{
 (  (  ) ((  )) (((  )))  )
}}\]
What can you put inside the parentheses? Words, and unsigned integers, which can both be
arbitrarily big.  And for doing math with integers, which are exact, not approximate numbers, it's
very important to be able to handle very large integers.
 
Also, a word or unsigned integer can appear all by itself as an S-expression, with \textbf{no}
parentheses.   Then it's referred to as an atom (which means indivisible in Greek).
If an S-expression isn't an atom, then it's called a list, and it's considered to be a list
of elements, with a first element, a second, a third, etc.
\[
    \mbox{\texttt{(1 2 3)} and \texttt{((x 1) (y 2) (z 3))}}
\]
are both lists with three elements.
 
And LISP is a functional or expression-based language, not an imperative
or statement-based language.
 
Everything in LISP 
(that's instructions rather than data)
is built up by applying 
functions to arguments like this:
\[\mbox{\texttt{
 (f x y)
}}\]
This indicates that the function $f$ is applied to the arguments $x$ and $y$. 
In pure math that's normally written
\[
 f(x,y)
\] 
And \textbf{everything} is put into this function applied to arguments form.
 
For example,
\[\mbox{\texttt{
 (if condition true-value false-value)
}}\]
is how you choose between two values depending on whether a condition is true or not.
So ``if'' is treated as a three-argument function with the strange property that only
two of its arguments are evaluated.
For instance
\[\mbox{\texttt{
 (if true (+ 1 2) (+ 3 4))
}}\]
gives 3, and
\[\mbox{\texttt{
 (if false (+ 1 2) (+ 3 4))
}}\]
gives 7.
Another pseudo-function is the quote function, which doesn't evaluate its only argument.
For instance
\[\mbox{\texttt{
 (' (a b c))
}}\]
gives \texttt{(a b c)} ``as is''.
In other words, this does \textbf{not} mean that the function $a$ should
be applied to the arguments $b$ and $c$.
The argument of quote is literal data, not an expression to be evaluated.
 
Here are two conditions that may be used in an ``if''.
\[\mbox{\texttt{
 (= x y)
}}\]
gives true if $x$ is equal to $y$.
\[\mbox{\texttt{
 (atom x)
}}\]
gives true if $x$ is an atom, not a list.
 
Next let me tell you about ``let'',
which is very important because it enables you to associate values with variables and 
to define functions.
\[\mbox{\texttt{
 (let x y expression)
}}\]
yields the value of ``expression'' in which $x$ is defined to be $y$.
You can use ``let'' either to define a function, or like an assignment
statement in an ordinary programming language.  However, the function
definition or assignment are only temporarily in effect, inside ``expression''.
In other words, the effect of ``let'' is invariably local.
 
Here are two examples of how to use ``let''.
Let $n$ be $1 + 2$ in $3 \times n$:
\begin{verbatim}
    (let n (+ 1 2)
       (* 3 n)
    )
\end{verbatim} 
This gives 9. And let $f$ of $n$ be $n \times n$ in $f$ of 10:
\begin{verbatim}
    (let (f n) (* n n)
       (f 10)
    )
\end{verbatim} 
This gives 100.
 
And in LISP we can take lists apart and then reassemble them.
To get the first element of a list, ``car'' it:
\[\mbox{\texttt{
 (car (' (a b c)))
}}\]
gives $a$.
To get the rest of a list, ``cdr'' it:
\[\mbox{\texttt{
 (cdr (' (a b c)))
}}\]
gives \texttt{(b c)}.
And to reassemble the pieces, ``cons'' it:
\[\mbox{\texttt{
 (cons (' a) (' (b c)))
}}\]
gives \texttt{(a b c)}.
 
And that's it, that gives you the general idea! LISP is a simple, but powerful formalism.
 
Before showing you two real LISP programs, let me point 
out that in LISP you don't speak of programs, they're called expressions.
And you don't run them or execute them, you evaluate them.
And the result of evaluating an expression is merely a value; there is no side-effect.
The state of the universe is unchanged.
 
Of course, mathematical expressions have always behaved like this.
But normal programming languages are not at all mathematical.
LISP \textbf{is} mathematical.
 
In this chapter we've used the factorial function ``!''\ several times.
So let's program that in a normal language, and then in LISP.
Then we'll write a program to take the factorials of an entire list of numbers, not just one.

\section*{Factorial in a normal language}

Recall that
$N! = 1 \times 2 \times 3 \times \cdots \times (N - 1) \times N$.
So $3! = 3 \times 2 \times 1 = 6$, and $4! = 24$, $5! = 120$.
This program calculates 5 factorial:
\begin{verse}
Set $N$ equal to 5.
\\
Set $K$ equal to 1.
\vspace{\baselineskip}
\\ 
LOOP: Is $N$ equal to 0? If so, stop. The result is in $K$.
\\
   If not, set $K$ equal to $K \times N$.
\\
   Then set $N$ equal to $N - 1$.
\\
   Go to LOOP.  
\end{verse}
When the program stops, $K$ will contain 120,
and $N$ will have been reduced to 0.

\section*{Factorial in LISP}

The LISP code for calculating factorial of 5 looks rather different:
\begin{verbatim}
   (let (factorial N) 
           (if (= N 0) 
               1 
               (* N (factorial (- N 1)))
           ) 
    (factorial 5)
   )
\end{verbatim}
This gives
\begin{verbatim}
   120
\end{verbatim}
 
Here's the definition of the function factorial of $N$ 
in English: 
If $N$ is 0, then
factorial of $N$ is 1, otherwise factorial of $N$ is $N$ times factorial of $N$ minus 1.

\section*{A more sophisticated example: list of factorials}

Here's a second example.
Given a list of numbers,
it produces the corresponding list of factorials:
\begin{verbatim}
   (let (factorial N) 
           (if (= N 0) 
               1 
               (* N (factorial (- N 1)))
           ) 

   (let (map f x)
      (if (atom x) 
          x 
          (cons (f (car x))
                (map f (cdr x))
          ) 
      ) 

    (map factorial (' (4 1 3 2 5)))

   ))
\end{verbatim}
This gives
\begin{verbatim}
   (24 1 6 2 120)
\end{verbatim}
 
The ``map'' function changes $(x \; y \; z\ldots)$ into $(f(x) \; f(y) \; f(z)\ldots)$.
Here's the definition of ``map'' $f$ of $x$
in English: The function ``map'' has two arguments, a function $f$ and a list $x$.
If the list $x$ is empty, then the result is the empty list.
Otherwise it's the list beginning with ($f$ of the
first element of $x$) followed by 
``map'' $f$ of the rest of $x$.
 
{\footnotesize
[For more on my conception of LISP, see my books \emph{The Unknowable}
and \emph{Exploring Randomness,} in which I explain LISP, and my book
\emph{The Limits of Mathematics,} in which I use it. These are technical books.]
}

\section*{My Diophantine Equation for LISP}

So you can see how simple LISP is. But actually, it's still not simple enough.
In order to make a LISP equation, I had to simplify LISP even more.
It's still a powerful language, you can still calculate anything you could before,
but the fact that some things are missing and you have to program them yourself
starts to become really annoying.
It starts to feel more like a machine language and less like a higher-level language.
But I can actually put together the entire equation for running LISP, more precisely,
the equation for evaluating LISP expressions.
Here it is!

\begin{center}
\begin{tabular}{|c|}
\hline
\\
\textbf{\emph{\large Exponential Diophantine Equation Computer:}}
\\ \\
   $L(k, n, x, y, z, \ldots) = R(k, n, x, y, z, \ldots)$
\\ \\
\textbf{LISP expression} $k$
\\
\textbf{Value of LISP expression} $n$
\\
\textbf{Time} $x, y, z, \ldots$
\\ \\
Two hundred page equation, 20,000 unknowns!
\\ \\
If the LISP expression $k$ has no value, then this equation 
\\
will have no solution. If the LISP expression $k$ has a value, 
\\
then this equation will have \textbf{exactly one} solution.
\\
In this unique solution, $n =$ the value of the expression $k$.
\\ \\
(Chaitin, \emph{Algorithmic Information Theory,} 1987.)
\\
\\
\hline
\end{tabular}
\end{center}

Guess what, I'm not going to actually show it to you here!
Anyway, this is a diophantine equation that's a LISP interpreter.
So that's a fairly serious amount of computational ability
packed into one diophantine equation.
 
We're going to put this computer equation to work for us later,
in Chapter VI, in order to show that $\Omega$'s randomness also 
infects number theory.

\section*{Concluding Remarks}

Some concluding remarks.
I've tried to share a lot of beautiful mathematics with you in this chapter.
It's math that I'm very passionate about, math that I've spent my life on!
 
Euclid's proof is classical, Euler's proof is modern, and my proof is post-modern!
Euclid's proof is perfect, but doesn't seem to particularly lead anywhere.
Euler's proof leads to Riemann and analytic number theory and much extremely technical and
difficult modern work.
My proof leads to the study of program-size complexity (which actually preceded this proof).
The work on Hilbert's 10th shows a marvelously unexpected application of the notion of a computably
enumerable set.
 
I'm not really that interested in the primes.  The fact that they now have applications in
cryptography makes me even less interested in the primes.  What is really interesting
about the primes is that there are still lots of things that we don't know about them,
it's the philosophical significance of this bizarre situation.  
What's really interesting is that even
in an area of math as simple as this, you immediately get to questions that nobody knows
how to answer and things begin to look random and haphazard!  
 
\textbf{Is there always some structure or pattern there waiting to be discovered?
Or is it possible that some things 
\textbf{really are} lawless, random, patternless---even in pure math, even
in number theory?}
 
Later I'll give arguments involving information, complexity and the 
extremely fundamental notion of irreducibility that
strongly suggest that the latter is in fact the case.
In fact, that's the \emph{raison d'\^etre} for this book.
 
Based on his own rather individualistic viewpoint---which happens to be 
entirely different from mine---Wolfram gives many interesting examples 
that \textbf{also} strongly suggest
that there are lots of unsolvable problems in mathematics. I think
that his  book on \emph{A New Kind of Science}
is extremely interesting and provides further evidence that lots of things are unknowable,
that the problem is really serious.
Right away in mathematics, in 2000 year-old math, you get into trouble, and there seem
to be limits to what can be known.  Epistemology, which deals with what we can know and why,
is the branch
of philosophy that worries about these questions.  So this is really a book about
epistemology, and so is Wolfram's book.
We're actually working on philosophy as well as math and physics!
 
To summarize, in this chapter we've seen that the computer is a powerful new mathematical
concept illuminating many questions in mathematics and metamathematics.
And we've seen that while formal axiomatic systems are a failure, formalisms for computing
are a brilliant success.
 
To make further progress on the road to $\Omega$, we need to add more to the stew. We need
the idea of digital information---that's measured by the size of computer programs---and
we also need the idea of irreducible digital information, which is a kind of randomness.
In the next chapter we'll discuss the sources of these ideas. We'll see that the
idea of complexity comes from biology, the idea of digital information comes from
computer software, and the idea of irreducibility---that's my particular contribution---can
be traced back to Leibniz in 1686.

\chapter*{Chapter III---\\Digital Information: DNA/Software/Leibniz}
\addcontentsline{toc}{chapter}{Chapter III---Digital Information: DNA/Software/Leibniz}
\markright {Chapter III---Digital Information: DNA/Software/Leibniz}

\begin{quote}
``Nothing is more important than to see the sources of invention which are,
in my opinion, more interesting than the inventions themselves.''---\emph{Leibniz,}
as quoted in P\'olya, \emph{How to Solve It.}
\end{quote}

In this chapter I'm going to show you the ``sources of invention''
of the ideas of digital information, program-size complexity, and algorithmic
irreducibility or randomness.
In fact, the genesis of these ideas can be traced to DNA, to software, and to Leibniz himself.
 
Basically the idea is just to measure the size of software in bits, that's it!
But in this chapter I want to explain to you just why this is so important, 
why this is such a universally applicable idea.
We'll see that it plays a key role in biology, where DNA is that software,
and in computer technology as well of course, and even in Leibniz's metaphysics
where he analyzes what's a law of nature, and how can we decide if something
is lawless or not. And later, in Chapter VI, we'll see that it plays a key role in
metamathematics, in discussing what a FAS can or cannot achieve.
 
Let's start with Leibniz.

\section*{Who was Leibniz?}

Let me tell you about Leibniz.
 
Leibniz invented the calculus, invented binary arithmetic, a superb mechanical
calculator, clearly envisioned symbolic logic, gave the name to topology
(analysis situs) and combinatorics, discovered Wilson's theorem (a primality test; see 
Dantzig, \emph{Number, The Language of Science}),
etc.\ etc.\ etc.
 
Newton was a great physicist, but he was definitely inferior to Leibniz both as a mathematician
and as a philosopher. And Newton was a rotten human being---so much so that 
Djerassi and Pinner call their recent book \emph{Newton's Darkness.} 
 
Leibniz invented the calculus, published it, wrote letter after letter to
continental mathematicians to explain it to them, initially received all the credit for this
from his contemporaries, and then was astonished to learn that Newton, who had
never published a word on the subject, claimed that Leibniz had stolen it all from
him. Leibniz could hardly take Newton seriously!
 
But it was Newton who won, not Leibniz.
 
Newton bragged that he had destroyed Leibniz
and rejoiced in Leibniz's death after Leibniz was abandoned by his royal patron,
whom Leibniz had helped to become the king of England. 
It's extremely ironic that Newton's incomprehensible \emph{Principia}---written
in the style of Euclid's \emph{Elements}---was only appreciated
by continental mathematicians \textbf{after} they succeeded in translating
it into that effective tool, the infinitesimal calculus that Leibniz had taught them!
 
Morally, what a contrast!
Leibniz was such an elevated soul that he found good in
all philosophies: Catholic, Protestant, Cabala, medieval scholastics, the ancients,
the Chinese\ldots\ It pains me to say that Newton enjoyed witnessing the executions
of counterfeiters he pursued as Master of the Mint.
 
{\footnotesize [The science-fiction writer Neal Stephenson has recently published the first 
volume, \emph{Quicksilver,} 
of a trilogy about Newton versus Leibniz, and comes out strongly
on Leibniz's side.
See also Isabelle Stengers, \emph{La Guerre des sciences aura-t-elle lieu?,} a
play about Newton vs.\ Leibniz, and the above mentioned book, consisting
of two plays and a long essay, called \emph{Newton's Darkness.}]}
 
Leibniz was also a fine physicist. In fact, Leibniz was good at everything.
For example, there's his remark that ``music is the unconscious joy that the
soul experiences on counting without realising that it is counting.''
Or his effort to discern pre-historic human migration patterns through similarities
between languages---something that is now done with DNA!
 
So you begin to see the problem: Leibniz is at too high an intellectual level.
He's too difficult to understand and appreciate.  In fact, you can only 
really appreciate Leibniz if you are at \textbf{his} level.
You can only 
realize that Leibniz has anticipated you
\textbf{after} you've invented a new field by yourself---which as C. MacDonald Ross
says in his little Oxford University Press book \emph{Leibniz,} has happened to
many people.
 
In fact, that's what happened to me. I invented and developed my
theory of algorithmic information, and one day not so long ago, when asked
to write a paper for a philosophy congress in Bonn, I went back to a little
1932 book by Hermann Weyl, \emph{The Open World,} Yale University Press.
I had put it aside after being surprised to read in it that Leibniz in
1686 in his \emph{Discours de M\'etaphysique}---that's the original French,
in English, \emph{Discourse on Metaphysics}---had made a key observation
about complexity and randomness, \textbf{the} key observation that started
me on all this at age 15!
 
{\footnotesize [Actually, Weyl himself was a very unusual mathematician, Hilbert's mathematical heir and
philosophical opponent, who wrote a beautiful book on philosophy that I had
read as a teenager: \emph{Philosophy of Mathematics and Natural Science,}
Princeton University Press, 1949. In that book Weyl also recounts Leibniz's idea,
but the wording is not as sharp, it's not formulated as clearly and as dramatically as in Weyl's
1932 Yale book. Among his other works: an important book on relativity, \emph{Raum-Zeit-Materie}
(space, time, matter).]}
 
``Could Leibniz have done it?!'', I asked myself. I put the matter aside until such time as
I could check what he had actually said. Years passed\ldots\  Well, for my Bonn paper I finally
took the trouble to obtain an English translation of the \emph{Discours,}
and then the original French. 
And I tried to find out more about Leibniz.
 
It turns out that Newton wasn't the only important opponent that Leibniz had had.
You didn't realize that math and philosophy were such dangerous professions, did you?!
 
The satirical novel \emph{Candide} by Voltaire, which was made into a musical comedy
when I was a child, is actually a caricature of Leibniz.
Voltaire was Leibniz's implacable opponent, and a terrific Newton
booster---his mistress 
la Marquise du Ch\^atelet
translated Newton's \emph{Principia} into French.
Voltaire was against one and in favor of the other, not based on an understanding of
their work, but simply because Leibniz constantly mentions God, whereas Newton's work
seems to fit in perfectly with an atheist, mechanistic world view. 
This was leading up to the French revolution, which was against the Church just as much as 
it was against the Monarchy.
 
Poor Voltaire---if he had read Newton's private papers, he would have realised that
he had backed the wrong man!  Newton's beliefs were primitive and literal---Newton
computed the age of the world based on the Bible. Whereas Leibniz was never seen
to enter a church, and his notion of God 
was sophisticated and subtle.
Leibniz's God was a logical necessity
to provide the initial complexity to create the world, and is required because
nothing is necessarily simpler than something. That's Leibniz's answer to the 
question, 
``Why is there something rather than nothing? For nothing is \textbf{simpler} 
and easier than something.''
(\emph{Principles of Nature and Grace,} 1714, Section 7)
 
In modern language, this is like saying that the initial complexity of the universe
comes from the choice of laws of physics and the initial conditions to which these laws apply.
And if the initial conditions are simple, for example an empty universe or an exploding singularity, 
then all the initial complexity comes from the laws of physics.
 
The question of where all the complexity in the world comes from continues to fascinate
scientists to this day. For example, it's the focus of Wolfram's book \emph{A New Kind of Science.}
Wolfram solves
the problem of complexity by claiming that it's only an illusion, that the world is actually
very simple. For example, according to Wolfram
all the randomness in the world is only pseudo-randomness
generated by simple algorithms!  That's certainly a philosophical possibility, but
does it apply to \textbf{this} world?  Here it seems that quantum-mechanical randomness
provides an inexhaustible source of potential complexity, 
for example via ``frozen accidents'', 
such as biological mutations that change the course of evolutionary history.
 
Before explaining to what extent and how Leibniz anticipated the starting point for my theory,
let me recommend some good sources of information about Leibniz, which are not
easy to find.  On Leibniz's mathematical work, see the chapter on him
in E. T. Bell's \emph{Men of Mathematics} and Tobias Dantzig's \emph{Number, The Language of
Science.} 
On Leibniz the philosopher, see C. MacDonald Ross, \emph{Leibniz,}   
and
\emph{The Cambridge Companion to Leibniz,}
edited by Nicholas Jolley. 
For works \textbf{by} Leibniz, including his \emph{Discourse on Metaphysics}
and \emph{Principles of Nature and Grace,}
see
G. W. Leibniz, \emph{Philosophical Essays,} 
edited and translated by Roger Ariew and Daniel Garber.
 
Actually, let me start by telling you about Leibniz's discovery of binary
arithmetic, which in a sense marks the very beginning of information theory,
and then I'll tell you what I discovered in the \emph{Discours de M\'etaphysique.}

\section*{Leibniz on Binary Arithmetic and Information}

As I said, Leibniz discovered base-two arithmetic, and he was extremely enthusiastic
about it.  About the fact that
\[
 10110 \; = \; 1 \times 2^4
 + 0 \times 2^3
 + 1 \times 2^2
 + 1 \times 2^1
 + 0 \times 2^0
\]
\[
 \; = \; 16 + 4 + 2
 \; = \; 22 \; ???
\] 
Yes indeed, he sensed in the 0 and 1 bits the 
combinatorial power to create the entire universe, which is
exactly what happens in modern digital electronic computers and the rest of our
digital information technology: CD's, DVD's, digital cameras, PC's, \ldots\
It's all 0's and 1's, that's our entire image of the world!
You just combine 0's and 1's, and you get everything.
 
And Leibniz was very proud to note how easy it is to perform calculations with binary
numbers, in line with what you might expect if you have reached the logical bedrock
of reality.  Of course, that observation was also made by computer engineers three
centuries later; the idea is the same, even though the language and the cultural 
context in which it is formulated is rather different.
 
Here is what Dantzig has to say about this in his \emph{Number, The Language of Science:}
\begin{quotation}
``It is the mystic elegance of the binary system that made Leibniz exclaim:
\emph{Omnibus ex nihil ducendis sufficit unum.} (One suffices to derive all out of nothing.)
[In German: ``Einer hat alles aus nichts gemacht.'' 
Word for word: ``One has all from nothing made.'']
Says Laplace:
\begin{quotation}
`Leibniz saw in his binary arithmetic the image of Creation\ldots\
He imagined that Unity represented God, and Zero the void; that the Supreme Being drew
all beings from the void, just as unity and zero express all numbers in his system
of numeration\ldots\ I mention this merely to show how the prejudices of childhood may cloud the
vision even of the greatest men!\ldots'
\end{quotation}
Alas! What was once hailed as a monument to monotheism
ended in the bowels of a robot. For most of the high-speed calculating machines [computers]
operate on the \emph{binary} principle.''
\end{quotation}
In spite of the criticism by Laplace, Leibniz's vision of creating the world from 0's and 1's
refuses to go away.  
In fact, it has begun to inspire some contemporary physicists,
who probably have never even heard of Leibniz.
 
Opening the 1 January 2004 issue
of the prestigious journal \emph{Nature,}
I discovered a book review entitled ``The bits that make up the Universe.''
This turned out to be a review of von Baeyer, \emph{Information: The New Language of Science}
by Michael Nielsen, himself coauthor of the large and authoritative 
Nielsen and Chuang, \emph{Quantum Computation and Quantum Information.}
Here's what Nielsen had to say:
\begin{quotation}
``What is the Universe made of? A growing number of scientists suspect that information
plays a fundamental role in answering this question. Some even go as far as to suggest
that information-based concepts may eventually fuse with or replace traditional notions
such as particles, fields and forces. The Universe may literally be made of information,
they say, an idea neatly encapsulated in physicist John Wheeler's slogan: `It from bit'
[Matter from information!]\ldots\
These are speculative ideas, still in the early days of development\ldots\
Von Baeyer has provided an accessible and engaging overview of the emerging role of
information as a fundamental building block in science.''
\end{quotation}
So perhaps Leibniz was right after all! At any rate, it certainly is a grand vision!

\section*{Leibniz on Complexity and Randomness}

Okay, it's time to look at Leibniz's \emph{Discourse on Metaphysics!} Modern science was really just
beginning then.  And the question that Leibniz raises is, how can we tell the difference
between a world that is governed by law, one in which science can apply, and a lawless world?
How can we decide if science actually works?!
In other words how can we distinguish between a set of observations that obeys a mathematical
law, and one that doesn't?
 
And to make the question sharper, Leibniz asks us to think about scattering points at random
on a sheet of paper, by closing our eyes and stabbing at it with a pen many times.
Even then, observes Leibniz, there will always be a mathematical law that passes through
those very points!
 
Yes, this is certainly the case. For example, if the three points are $(x,y) = (a,A), (b,B), (c,C)$,
then the following curve will pass through those very points:
\[
         y \; = \; A(x-b)(x-c)/(a-b)(a-c)
           \; + \; B(x-a)(x-c)/(b-a)(b-c)
\]
\[
           \; + \; C(x-a)(x-b)/(c-a)(c-b).
\] 
Just plug in $a$ every time you see an $x$, 
and you'll see that the entire right-hand side reduces to $A$.  
And this also works if you set $x = b$ or $x = c$.
Can you figure out how this works and write the corresponding equation for four points? This is
called Lagrangian interpolation, by the way.
 
So there will always be a mathematical law, no matter what the points are, no matter how they
were placed at random!
 
That seems very discouraging. How can we decide if the universe is capricious or if science
actually works?
 
And here is Leibniz's answer: If the law has to be extremely complicated 
(``fort compos\'ee'') then the points are placed
at random, 
they're ``irr\'egulier'',
not in accordance with a scientific law.  But if the law is simple, then it's a genuine
law of nature, we're not fooling ourselves!
 
See for yourself, look at Sections V and VI of the \emph{Discours.}
 
The way that Leibniz summarizes his view of what precisely it is that makes the
scientific enterprise possible is this: 
``God has chosen that which is the most perfect,
that is to say, in which at the same time the hypotheses are as simple as possible, and
the phenomena are as rich as possible.''
The job of the scientist, of course, is then to figure out these simplest possible hypotheses.
 
(Please note that these ideas of Leibniz are much stronger than Occam's razor, 
because they tell us \textbf{why} Occam's razor works, why it is necessary.
Occam's razor merely asserts that the simplest theory is best.)
 
But how can we \textbf{measure} the complexity of a law and compare it with the complexity of
the data that it attempts to explain?  
Because it's only a valid law if it's simpler than the data, hopefully much simpler.
Leibniz does not answer that question, which bothered
Weyl immensely.
But Leibniz had all the pieces, he only had to put them together. For he worshiped 0 and 1,
and appreciated the importance of calculating machines.
 
\textbf{
The way I would put it is like this: 
I think of a scientific theory as a binary computer program for calculating the observations,
which are also written in binary.
And you have a law of nature if there is compression, 
if the experimental data is compressed into a computer program that has a smaller 
number of bits 
than are in the data that it explains.
The greater the degree of compression, the better the law, the more you understand the data.
}
 
\textbf{
But if the experimental data cannot be compressed, if the smallest program for calculating
it is just as large as it is (and such a theory can always be found, that can always be done,
that's so to speak our ``Lagrangian interpolation''), 
then the data is lawless,
unstructured, patternless, not amenable to scientific study, incomprehensible.
In a word, random, irreducible!
}
 
And \textbf{that} was the idea that burned in my brain like a red-hot ember when I was 15 years old! 
Leibniz would have understood it instantaneously!
 
But working out the details, and showing that mathematics contains such randomness---that's 
where my $\Omega$ number comes in---that's the
hard part, that took me the rest of my life. As they say, genius is 1\% inspiration and 99\%
perspiration.  And the devil is in the details.  And the detail that cost me the most effort
was the idea of self-delimiting information, that I'll explain later in this chapter, and without
which, as we'll see in Chapter VI, there would in fact be no $\Omega$ number.
 
The ideas that I've been discussing are summarized in the following diagrams:
 
\begin{center}
\begin{tabular}{|c|}
\hline
\\
\textbf{\emph{\large Leibniz's Metaphysics:}}
\\ \\
ideas $\longrightarrow$ \textbf{\large Mind of God} $\longrightarrow$ universe
\\ \\
The ideas are simple, but the 
\\
universe is very complicated!
\\
(If science applies!)
\\
God minimizes the left-hand side,
\\
and maximizes the right-hand side.
\\
\\
\hline
\end{tabular}
\end{center}
 
\begin{center}
\begin{tabular}{|c|}
\hline
\\
\textbf{\emph{\large Scientific Method:}}
\\ \\
theory $\longrightarrow$ \textbf{\large Computer} $\longrightarrow$ data
\\ \\
The theory is concise, the data isn't.
\\
The data is compressed into the theory!
\\ \\
\textbf{Understanding is compression!}
\\
\textbf{To comprehend is to compress!}
\\
\\
\hline
\end{tabular}
\end{center}
 
\begin{center}
\begin{tabular}{|c|}
\hline
\\
\textbf{\emph{\large Algorithmic Information Theory:}}
\\ \\
binary program $\longrightarrow$ \textbf{\large Computer} $\longrightarrow$ binary output 
\\ \\
What's the smallest program that produces a given output?
\\
If the program is concise and the output isn't, we have a theory.
\\ \\
If the output is random, then no compression is possible,
\\
and the input has to be the same size as the output.
\\
\\
\hline
\end{tabular}
\end{center}
 
And here's a diagram from Chapter II, which I think of
in exactly the same way:
 
\begin{center}
\begin{tabular}{|c|}
\hline
\\
\textbf{\emph{\large Mathematics (FAS):}}
\\ \\
axioms $\longrightarrow$ \textbf{\large Computer} $\longrightarrow$ theorems
\\ \\
The theorems are compressed into the axioms!
\\
\\
\hline
\end{tabular}
\end{center}
 
I think of the axioms as a computer program for generating all theorems.
I measure the amount of information in the axioms by the size of this program.
Again, you want the smallest program, the most concise set of axioms,
the minimum number of assumptions (hypotheses) that give you that set of theorems.
One difference with the previous diagrams:  This program never halts, it keeps
generating theorems forever; it produces an infinite amount of output, not a finite
amount of output.
 
And the same ideas sort of apply to biology:
 
\begin{center}
\begin{tabular}{|c|}
\hline
\\
\textbf{\emph{\large Biology:}}
\\ \\
DNA $\longrightarrow$ \textbf{\large Pregnancy} $\longrightarrow$ organism
\\ \\
The organism is determined by its DNA!
\\
The less DNA, the simpler the organism.
\\ \\
Pregnancy is decompression 
\\
of a compressed message.
\\
\\
\hline
\end{tabular}
\end{center}
 
So I think of the DNA as a computer program for calculating the organism.
Viruses have a very small program, single-celled organisms have a larger program,
and multi-celled organisms need an even larger program.

\section*{The Software of Life: Love, Sex and Biological Information}

A passionate young couple can make love several times a night, every night.
Nature does not care about recreational sex; the reason that men and women
romanticize about each other and are attracted and fall in love, is so that they
will be fooled into having children, even if they think they are trying to avoid conception!
They're actually trying very hard to transmit information from the male to the female.
Every time a man ejaculates inside a woman he loves, an enormous number of sperm cells,
each with half the DNA software for a complete individual, try to fertilize
an egg.
 
Jacob Schwartz once surprised a computer science class by calculating the
bandwidth of human sexual intercourse, the rate of information transmission achieved in human
love-making.  I'm too much of a theoretician to care about the exact answer, which
anyway depends on details like how you measure the amount of time that's involved, but
his class was impressed that the bandwidth that's achieved is quite respectable!
 
What is this software like?  It isn't written in 0/1 binary like computer software. 
Instead DNA is written 
in a 4-letter alphabet, the 4 bases that can be each rung of the twisted double-helix ladder that is
a DNA molecule.
Adenine, A, thymine, T, guanine, G, and cytosine, C, are those four letters.
Individual genes, which code for a single protein, are kilobases of information.
And an entire human genome is measured in gigabases, so that's sort of like gigabytes
of computer software.
 
\textbf{Each cell in the body has the same DNA software, the complete genome, but
depending on the kind of tissue or the organ that it's in,
it runs different
portions of this software, while using many basic subroutines that are common to all cells.}
 
\begin{center}
\begin{tabular}{|c|}
\hline
\\
Program (10011\ldots) $\longrightarrow$ \textbf{\large Computer} $\longrightarrow$ Output
\\ \\
DNA (GCTATAGC\ldots) $\longrightarrow$ \textbf{\large Development} $\longrightarrow$ Organism
\\
\\
\hline
\end{tabular}
\end{center}
 
And this software is highly conservative, much of it is quite ancient:
Many common subroutines are shared among 
fruitflies, invertebrates, mice and humans,
so they have to have originated in an ancient common ancestor.
In fact, there is surprisingly little difference between a chimp and a human, or even
between a mouse and human. 
 
We are not that
unique; Nature likes to re-use good ideas. Instead of starting afresh each time,
Nature ``solves'' new problems by
patching---that is, slightly modifying or mutating---the solutions to old problems, as the
need arises. Nature is a cobbler, a tinkerer.  It's much too much work,
it's much too expensive, to start over again each time.
Our DNA software accumulates by accretion, it's a beautiful patch-work quilt!
And our DNA software 
also includes all those frozen accidents, those mutations due to DNA copying errors
or ionizing radiation, which is a possible pathway for quantum uncertainty to be incorporated
in the evolutionary record.
In a sense this is an amplification mechanism, one that magnifies quantum uncertainty into an
effect that is macroscopically visible.
 
Another example of this constant re-use of biological ideas
is the fact that the human embryo is briefly a fish, or, more generally, the fact that the
development of the embryo tends to recapitulate 
the evolutionary history leading to that particular organism.
The same thing happens in human software, where old code is hard to get rid of, since much
is built on it, and it quickly becomes full of add-ons and patches.
 
When a couple falls madly in love, what is happening in information-theoretic terms is that they
are saying, what nice subroutines the other person has, let's try combining some of her
subroutines with some of his subroutines, let's do that right away!  
That's what it's really all about!
Just like a child cannot have a first name
followed by the entire father's and mother's names, because then names would 
double in size each generation and soon
become too long to remember (although the Spanish
nobility attempted to achieve this!), a child cannot have
all the software, all the DNA, from both parents.  So a child's DNA consists of a random
selection of subroutines, half from one parent, half from the other.
You only pass half of your software to each child, just as you can't include your whole name
as part of your child's name.
 
And in my great-grandmother's generation in the old country, women would have a dozen children,
most of whom would die before puberty.
So you were trying a dozen mixes of DNA subroutines from both parents. 
(In the middle ages, babies weren't even named til they were a year old, since so many of them
would die the first year.)  Now instead of trying to keep women pregnant all the time, we
depend on massive amounts of expensive medical care to keep alive one or two children, no matter
how unhealthy they are.
While such medical care is wonderful for the individual, 
the quality of the human gene pool inevitably deteriorates
to match the amount of medical care that is available.
The more medical care there is, the sicker people become!
The massive amounts of medical care become part of the ecology,
and people come to depend on it to survive til child-bearing and for their
children to survive in turn\ldots
 
Actually, many times a woman will not even realize that she was briefly pregnant, because
the embryo was not at all viable and the pregnancy quickly aborted itself. 
So parents still actually do
try lots of different combinations of their subroutines, even if they only have a few children.
 
We are just beginning to discover just how powerful biological software is.
For example, take human life expectancy.  Ageing isn't just a finally fatal accumulation
of wear and tear. No, death is programming, the life-expectancy of an individual is determined
by an internal alarm clock.  Just like the way that ``apoptosis'', the process by which individual cells
are ordered to self-destruct, is an intrinsic part of new cell growth and the plasticity
in the organism (it must continually tear itself down in order to continually rebuild itself),
a human organism systematically starts to self-destruct according to a pre-programmed schedule,
no doubt designed to get him or her out of the way and no longer competing for food with
children and
younger child-bearing individuals. (Although grandparents seem to serve a role to help caring
for children.)  So this is just software!  Well, then it should be easy to modify, you just
change a few key parameters in your code!  And, in fact, people are beginning to think this can
be done. For one invertebrate, the nematode worm \emph{C. elegans,} it has \textbf{already} been done.
(This is work of Cynthia Kenyon.) For humans, maybe in 50 or a hundred
years such modification will be common practice, or permanently incorporated into the human genome.
 
By the way, bacteria may not have different sexes like we do, but
they certainly do pass around useful DNA subroutines, which is how antibiotics quickly
create superbugs in hospitals.  This is called horizontal gene transfer, because the genes
go to contemporaries, not to descendents (that's vertical). And this process results in a 
kind of ``bacterial intelligence'', in their ability to quickly adapt to a new environment.

\section*{What are bits and what are they good for?}

The Hindus were fascinated by large numbers.  There is a parable of a mountain
of diamond, the hardest substance in the world, one-mile high and one-mile wide.
Every thousand years a beautiful golden bird flies over this diamond mountain and
sheds a single tear, which wears away a tiny, tiny bit of the diamond.  And the
time that it will take for these tears to completely wear away the diamond mountain,
that immense time, it is but the blink of an eye to one of the Gods!
 
And when I was a small child, like many future mathematicians, I was also fascinated
by large numbers. I would sit in a stair-well of the apartment building where we lived
with a piece of paper and a pencil, and I would write down the number 1 and then double it,
and then double it again, and again, until I ran out of paper or of patience:
\begin{center}
\begin{tabular}{c}
$2^0 = 1$ \\
$2^1 = 2$ \\
$2^2 = 4$ \\
$2^3 = 8$ \\
$2^4 = 16$ \\
$2^5 = 32$ \\
$2^6 = 64$ \\
$2^7 = 128$\\
$2^8 = 256$ \\
$2^9 = 512$\\
$2^{10} = 1024$ \\
$2^{11} = 2048$ \\
$2^{12} = 4096$ \\
$2^{13} = 8192$ \\
$2^{14} = 16384$ \\
$2^{15} = 32768$ \\
\ldots
\end{tabular}
\end{center}
 
A similar attempt to reach infinity is to take a piece of paper and fold it in half, and
then fold that in half, and then do it again and again, until the paper is too thick to
fold, which happens rather quickly\ldots
 
What does all this have to do with information? 
Well, the doubling has a lot to do with bits of information, with messages in binary,
and with the bits or \textbf{b}inary dig\textbf{its} 0 and 1, 
out of which we can build the entire universe!
 
One bit of information, a string 1-bit long, can distinguish two different possibilities.
Two bits of information, a string 2-bits long, can distinguish four different possibilities.
Each time you add another bit to the message, you double the number of possible messages.
There are 256 possible different 8-bit messages, 
and there are 1024 possible different 10-bit messages.
So ten bits of information is roughly equivalent to three digits of information,
because ten bits can represent 1024 possibilities, while three digits can represent 1000 possibilities.
 
Two possibilities:
\[
0, 1.
\]
Four possibilities:
\[
00, 01, 10, 11.
\]
Eight possibilities:
\[
000, 001, 010, 011, 100, 101, 110, 111.
\]
Sixteen possibilities:
\[
0000, 0001, 0010, 0011, 0100, 0101, 0110, 0111,
\]
\[
1000, 1001, 1010, 1011, 1100, 1101, 1110, 1111.
\]
 
So that's what raw information is, that's what raw binary information looks like inside a computer.
Everything is in binary, everything is built up out of 0's and 1's.
And computers have clocks to synchronize everything, and in a given clock cycle an
electrical signal represents a 1, and no signal represents a 0.  These two-state systems
can be built very reliably, and then everything in a computer is built up out of them.
 
What do these binary strings, these strings of bits, these binary messages, what do they
represent?  Well, the bit strings can represent many things. For instance, numbers:
\begin{center}
\begin{tabular}{lr}
0: & 0 \\
1: & 1 \\
2: & 10 \\
3: & 11 \\
4: & 100 \\
5: & 101 \\
6: & 110 \\
7: & 111 \\
8: & 1000 \\
9: & 1001 \\
10: & 1010
\end{tabular}
\end{center}
 
Or you can use an 8-bit ``byte'' to represent a single character, using the ASCII character-code,
which is used by most computers.
For example:
\begin{center}
\begin{tabular}{lr}
A: & 0100 0001 \\
B: & 0100 0010 \\
C: & 0100 0011 \\
a: & 0110 0001 \\
b: & 0110 0010 \\
c: & 0110 0011 \\
0: & 0011 0000 \\
1: & 0011 0001 \\
2: & 0011 0010 \\
(: & 0010 1000 \\
): & 0010 1001 \\
blank: & 0010 0000
\end{tabular}
\end{center}
 
And you can use a string of 8-bit bytes to represent character strings.
For example, in the LISP programming language 
how do you write the arithmetic expression for 12 plus 24?
This is written in LISP as
\[\mbox{\texttt{
  (+ 12 24) 
}}\]
This 9-character LISP expression
indicates the result of adding the numbers 12 and 24, and its
value is therefore the number 36. 
This LISP expression has 9 characters because there are seven visible characters plus two blanks, 
and inside a computer it becomes
a $9 \times 8 = 72$-bit bit string.
And you can easily build up more complicated objects like this.
For example, the LISP expression for adding the product of 3 and 4 to the product of 5 and 6 is this:
\[\mbox{\texttt{
  (+ (* 3 4) (* 5 6))
}}\]
This has 19 characters and is therefore represented as a $19 \times 8 = 152$-bit bit string.
But this 152-bit bit string could also be interpreted as an extremely large number, one that is
152 bits or, recalling that 10 bits is about 3 digits worth of information, about
$(152/10 = 15) \times 3 = 45$ digits long.
 
So bit strings are neutral, they are pure syntax, but it is a convention regarding semantics
that gives them meaning. 
In other words, bit strings can be used to represent many things: numbers,
character strings, LISP expressions, and even, as we all know, color pictures, since these
are also in our computers.  In each case it is the bit string plus a convention for interpreting
them that enables us to determine its meaning, be it a number, a character string, a LISP expression,
or a color picture.  For pictures one must in fact give the height and width of the picture in pixels,
or picture elements (points on a display) and then the red/green/blue intensity at each point, each
with 8-bit precision, which is a total of 24 bits per pixel.

\section*{What is biological information?}

Here is a specific example that is of great interest to us as human beings.
Our genetic information (DNA) is written using an alphabet of 4 symbols:
\[
    \mbox{A, C, G, T}
\]
These are the symbols for each of the possible bases at each rung of a DNA double-helix.
So each of these bases is exactly 2 bits of information, since two bits enables us to
specify exactly $2 \times 2 = 4$ possibilities.
 
In turn, a triple (codon) of 3 bases ``codes for'' (specifies) a specific amino acid.
So each DNA codon is $3 \times 2 = 6$ bits of information.
An individual protein is determined by giving the linear sequence of its amino acids 
in a portion of the DNA that is called a gene.
The gene determines the amino-acid ``spinal cord'' of each protein.
Once synthesized in a ribosome, the protein immediately folds into a complicated three-dimensional
shape. This folding process isn't well understood and in the current state of the art requires
massive amounts of computation to simulate.  It is the complicated geometrical form of the protein
that determines its biological activity.  For example, enzymes are  proteins that catalyze
(greatly facilitate and speed-up specific chemical reactions) by holding the reagents close
to each other in exactly the right way for them to react with each other.
 
That's the story, roughly speaking, but
DNA is actually much more sophisticated than that.
For example, 
some proteins turn other genes on and off; in other words, they control gene ``expression''.
We are dealing here with a programming language that can perform complicated
calculations and run through sophisticated sequences of gene expression in response to changes in
environmental conditions!
 
And as I said before, the DNA
software of some of our cousin apes and near-relative mammals is surprisingly similar to our own.
DNA subroutines are strongly ``conserved'', they are reused constantly
across many different species. Many of our basic subroutines are present in much more primitive living
beings. They haven't changed much; Nature likes to re-use good ideas.
 
By the way, a human being has 3 giga-bases, that's 3Gb, of DNA: 
 
\begin{center}
\begin{tabular}{|c|}
\hline
\\
\textbf{\emph{\large Human}} = 3 giga-bases = 6 giga-bits !!!
\\
\\
\hline
\end{tabular}
\end{center}
 
The conventional units of biological information
are kilo-bases (genes), mega-bases (chromosomes) and giga-bases (entire genomes for organisms).
That's thousands, millions and billions of bases, respectively.

\section*{Compressing digital pictures and video}

It is said that
a picture is worth a thousand words. Indeed, a thousand-by-thousand-point picture has a million
pixels, and each pixel requires 3 bytes = 24 bits to specify the admixture of primary colors.
So a thousand by a thousand picture is 3 megabytes 
or 24 megabits
of information. We have all
experienced the frustration of having to wait for a large picture to take shape in a web browser.  
 
Digital video requires transmitting many pictures per second  
to create the illusion
of smooth, continuous motion, 
and this requires extremely high internet bandwidth to work well,  
say, thirty frames per second.
This kind of high-bandwidth
connection is now usually only available in-house within organizations, but not across the internet
between organizations.
It would be nice to be able to distribute HDTV, high-resolution digital video, across the web,
but this is unfortunately not practical at this time.
 
But there isn't really that much information there---you don't 
really have to transmit all those bits.
Compression techniques are widely used to speed-up the transmission of digital pictures
and video.  The idea is to take advantage of the fact that large regions of a photo may
be more or less the same, and that successive frames of a video may not differ that much from
each other.  So instead of transmitting the pictures directly, you just send compact descriptions
of the pictures, and then at the receiving end you de-compress them and 
re-create the original pictures.
It's sort of like the freeze-dried desiccated foods that you can take with you up a mountain.
Just add water and heat, and you recreate the original foods, which are mostly water and much
heavier than their freeze-dried ``compressions''.
 
Why does this compression and decompression work so well?  
It works so well because pictures are far from random.  If each
pixel had absolutely no connection with any of its neighboring pixels, and if successive frames
of a digital video were totally unconnected, then no compression technique would work.
Compression techniques are useless if they are applied to noise, to the mad jumble that you get
if an antenna is disconnected, because there is absolutely no pattern there to compress away.
 
Another way to put it is that the most informative picture is one in which each pixel
is a complete surprise. Fortunately real pictures are almost never like that.  But we shall
re-encounter this important theme later on, in Chapter V, when we'll take incompressibility and 
use it as the basis for
a mathematical definition of a random real number.  
And then, in Chapter VI, we'll discover such a random number in pure 
mathematics: the halting probability $\Omega$. 
$\Omega$ is an infinite sequence of bits in which there is no pattern, and there are no correlations. 
Its bits are
mathematical facts that cannot be compressed into axioms that are
more concise than they are.
So, surprisingly enough, digital TV and DVD compression/decompression technology 
may actually have something to do with
more important questions, namely with philosophy and the limits of knowledge!
 
Indeed, the philosophical, mathematical ideas anti-dated these practical applications 
of compression and decompression.
There were no DVD's when I started working on these ideas in the early 1960's.
At that time audio and video recording was analog, not digital. And the cheap computing
power necessary to compress and decompress audio and video was simply not available back then.
Vinyl records used to use the depth of a groove to record a sound, not 0's and 1's like a CD does.
I even remember my grandparents' totally mechanical record player, absolutely without
electricity, that still worked perfectly when I was a child on their old, heavy 78 rpm records.
What astonishing progress!
 
Too bad that human society hasn't progressed at the same pace that this technology has!
Indeed, in some ways we've retrogressed, we haven't progressed at all.
Technology can be improved easily enough, but the human soul, that's immensely hard to improve!
 
But according to John Maynard Smith and E\"ors Szathm\'ary
(see their books \emph{The Major Transitions in Evolution} and \emph{The Origins of Life}), 
this kind of technological progress isn't just
important in consumer electronics, it's also responsible for the major steps forward in the
evolution of life on this planet. It played an extremely important role in the creation of 
vastly new and improved life-forms, which are in a sense multiple \textbf{origins} of life, moments
when the life on this planet successfully re-invented itself. 
To be more specific, they see the major steps forward in evolution
as occurring
at each point in time when biological organisms were able to invent and take advantage of radically 
new and improved techniques for representing and storing biological information.
 
And this even extends to human society and to human history! 
In his books \emph{The New Renaissance} and \emph{Phase Change},
Douglas Robertson makes a convincing case that the major steps forward in human social evolution
were because of the invention of language, which separates us from the apes, 
because of written language, which makes
civilization possible, because of the printing press and cheap paper and cheap books, 
which provoked the Renaissance and
the Enlightenment, and because of the PC and the web, 
which are the current motor for drastic social change.
Each major transition occurred because it became possible for human society 
to record and transmit a vastly increased amount of information.
That may not be the \textbf{whole} story, but the amount of information that could be recorded and
remembered, the availability of information, certainly played a major role.

\section*{What is self-delimiting information?}

Okay, so far we've discussed what information is, and the important role that it plays in 
biology---and even in determining our form of human social organization!---and 
the idea of compression and decompression.  But there's another important idea that
I'd like to tell you about right away, which is how to make information be ``self-delimiting''.
What the heck is that?!
 
The problem is very simple: How can we tell where one binary message ends and another begins, 
so that we can have \textbf{several} messages in a row and not get confused?  
(This is also, it will turn out, an essential step
in the direction of the number $\Omega$!)
Well, if we know exactly how long each message is, then there is no problem.  But how
can we get messages to indicate how long they are?  
 
Well, there are many ways to do this, with increasing degrees of sophistication.
 
For example,
one way to do it is to double each bit of the original message, and then tack on a pair of unequal
bits at the end.
For example, the plain ordinary binary string
\[
   \mbox{\textbf{Original string}:} \;\; 011100
\] 
becomes the self-delimiting binary string
\[
   \mbox{\textbf{With bits doubled}:} \;\; 00\, 11 \, 11\, 11\, 00\, 00\, 01.
\]
The only problem with this technique, is that we've just doubled the length of every string! 
True, if we
see two such strings in a row, if we are given two such strings in a row,
we know just where to separate them.
For example,
\[
   \mbox{\textbf{Two self-delimiting strings}:} \;\; 00\, 11\, 11\, 11\, 00\, 00\, 01 \;\; 11\, 00\, 11\, 00\, 10
\]
gives us the two separate strings
\[
    011100, \;\; 1010.
\]
But we had to double the size of everything!  
That's too big a price to pay!
Is there a better way to make information self-delimiting?
Yes, there certainly is!
 
Let's start with
\[
   011100
\]
again.  But this time, let's put a prefix in front of this string, a header, that tells
us precisely how many bits there are in this message (6).  To do this, let's write 6 in
binary, which gives us 110 ($4+2$), and then let's use the bit-doubling trick on the prefix
so that we know where the header ends and the actual information begins:
\[
   \mbox{\textbf{With a header}:} \;\; 11\, 11\, 00\, 01 \;\; 011100 
\]
In this particular case,
using a bit-doubled header like this 
is no better than doubling the entire original message.
But if our original message were very long, doubling the header would be a lot more economical
than doubling everything.
 
Can we do even better?! Yes, we can, there's no need to double the entire header.
We can just put a prefix in front of the header that tells us how long the header is,
and just double \textbf{that}. So now we have \textbf{two} headers, 
and only the first one is doubled, nothing else is.
How many bits are there in our original header? Well, the header is  6 = 110, 
which is only 3 = 11 bits long,
which doubled gives us 11 11 01.  So now we've got this:
\[
   \mbox{\textbf{With two headers}:} \;\; 11\, 11\, 01 \;\; 110 \;\; 011100
\]
In this case using two headers is
actually longer than just using one header.  But if our original message were very, very long,
then this would save us a lot of bits.
 
And you can keep going on and on in this way, 
by adding more and more headers, each of which tells us the size
of the \textbf{next} header, and only the very \textbf{first} header is doubled!
 
So I think you get the idea.  There are lots and lots of ways to make information self-delimiting!
And the number of bits that you have to add in order to do this isn't very large. 
If you start with an $N$-bit message, you can make it self-delimiting by only adding the number
of bits you need to specify $N$ in a self-delimiting manner. That way you know where all the header
information, the part that tells us $N$, ends, and where the actual information, all $N$ bits of it,
begins. Okay? Is that clear?  
 
This is more important than it seems.  First of all, self-delimiting information is \textbf{additive},
which means that the number of bits it takes to convey two messages in a row is just the sum of
the number of bits it takes to convey each message separately.  In other words, if you think of
individual messages as subroutines, you can combine many subroutines into a larger program without having
to add any bits---as long as you know just how many subroutines there are, which you usually do.
 
And that's not all.
We'll come back to self-delimiting information later when we discuss how the halting probability
$\Omega$ is defined, and why the definition works and actually defines a probability.
In a word, it works because we stipulate that these programs, the ones that
the halting probability $\Omega$ counts, have to be self-delimiting binary information.
Otherwise $\Omega$ turns out to be nonsense, because you can't use a single number to count the
programs of \textbf{any} size that halt, you can only do it for programs of a particular size,
which isn't very interesting.
 
We'll discuss this again later.
It's probably the most difficult thing to understand about $\Omega$.
But even if I can't explain it to you well enough for you to understand it, 
just remember that we're always
dealing with self-delimiting information, 
and that all the computer programs that we've been and will be considering are self-delimiting binary,
that'll give you the general idea of how things work.
Okay?
 
\begin{center}
\begin{tabular}{|c|}
\hline
\\
\textbf{\emph{\large Algorithmic Information Theory:}}
\\ \\
self-delimiting information $\longrightarrow$ \textbf{\large Computer} $\longrightarrow$ output 
\\ \\
{\small What's the smallest self-delimiting program that produces a given output?}
\\
{\small The size of that program in bits is the complexity $H$(output) of the output.}
\\
\\
\hline
\end{tabular}
\end{center}
 
After all, you don't want to become an expert in this field, (at least not yet!), you just
want to get a general idea of what's going on, and of what the important ideas are.

\section*{More on Information Theory \& Biology}

On the one hand, different kinds of organisms need different amounts of DNA.
And the DNA specifies how to construct the organism, and how it works.
 
\begin{center}
\begin{tabular}{|c|}
\hline
\\
DNA $\longrightarrow$ \textbf{Development/Embryogenesis/Pregnancy} $\longrightarrow$ Organism 
\\ \\
\emph{\textbf{Complexity/Size Measure:}}
\\
\emph{kilo/mega/gigabases of DNA.}
\\ \\
And one can more or less classify organisms into a complexity hierarchy 
\\
based on the amount of DNA they need: viruses, cells without nuclei, 
\\
cells with nuclei, multi-celled organisms, humans\ldots
\\
\\
\hline
\end{tabular}
\end{center}
 
That's the biological reality. Now let's abstract a highly-simplified mathematical model of this:
 
\begin{center}
\begin{tabular}{|c|}
\hline
\\
Program $\longrightarrow$ \textbf{\large Computer} $\longrightarrow$ Output 
\\ \\
\emph{\textbf{Complexity/Size Measure:}}
\\
\emph{bits/bytes kilo/mega/gigabytes of software.}
\\ \\
$H$(Output) = size of smallest program for computing it.
\\
\\
\hline
\end{tabular}
\end{center}
 
Actually, a byte = 8 bits, and kilo/mega/gigabytes are the conventional units of software size, 
leading to confusion with bases (= 2 bits), which also starts with the letter ``b''!
One way to avoid the confusion is to use capital ``B'' for ``bytes'', and lowercase ``b'' for ``bases''.
So KB, MB, GB, TB = kilo/mega/giga/terabytes, and Kb, Mb, Gb, Tb = kilo/mega/giga/terabases.
And kilo = one thousand ($10^3$), mega = one million ($10^6$), 
giga = one billion ($10^9$), tera = a million millions ($10^{12}$).
However as a theoretician, I'll think in just plain bits, usually, I won't use any of these
practical information size measures.
 
Now let me be provocative:
Science is searching for the DNA for the Universe!
 
\begin{center}
\begin{tabular}{|c|}
\hline
\\
\textbf{\emph{\large A Mixed-Up Model:}}
\\ \\
DNA $\longrightarrow$ \textbf{\large Computer} $\longrightarrow$ Universe!
\\ \\
What is the complexity of the entire universe?
\\
What's the smallest amount of DNA/software to construct the world?
\\
\\
\hline
\end{tabular}
\end{center}
 
So, in a way, AIT is inspired by biology!  
Certainly biology is the domain of the complex, 
it's the most obvious example of complexity in the world,
not physics, where there are simple unifying principles, 
where there are simple equations that explain everything!
But can AIT contribute to biology?!
Can ideas flow from AIT to biology, rather than from biology to AIT?
One has to be very careful!
 
First of all, the 
\begin{center}
\begin{tabular}{|c|}
\hline
\\
\textbf{\emph{\large Algorithmic Information Theory:}}
\\ \\
program $\longrightarrow$ \textbf{\large Computer} $\longrightarrow$ output 
\\ \\
What's the smallest program for a given output?
\\
\\
\hline
\end{tabular}
\end{center}
model of AIT imposes \textbf{no} time limit
on the computation. But in the world of biology, 9 months is already a long time to produce
a new organism, that's a long pregnancy. 
And some genes are repeated because one needs large quantities of the protein
that they code for.  So AIT is much too simple-minded a model to apply to biology.
DNA is \textbf{not} the minimum-size for the organism that it produces.  
 
In my theory a program should have no redundancy. But in DNA and the real world,
redundancy is good. If there were none, any change in a message yields \textbf{another}
valid message. But redundancy makes it possible to do error correction and error detection,
which is very important for DNA (indeed DNA
actually contains \textbf{complementary base pairs}, not individual bases, so that is already
a form of redundancy).
 
But unfortunately, in order for us to
be able to prove theorems,
we need to use a less complicated model, a toy model, one that applies well to meta-mathematics, 
but that does not apply well to biology.
Remember:
Pure mathematics is much easier to understand, much simpler, than the messy real world!
 
\begin{center}
\begin{tabular}{|c|}
\hline
\\
\textbf{\emph{\large Complexity!}}
\\ \\
\emph{It's an idea from Biology imported into Mathematics.}
\\ \\
\emph{It's \textbf{not} an idea from Mathematics imported into Biology.}
\\
\\
\hline
\end{tabular}
\end{center}
 
AIT also fails in biology in another crucial respect.
Look at a crystal and at a gas.  One has high program-size complexity, the other has low
program-size complexity, but \textbf{neither} is organized, neither is of any biological interest!
 
\begin{center}
\begin{tabular}{|c|}
\hline
\\
\textbf{\large \emph{Gas, Crystal}}
\\ \\
$H$(Crystal) is very low because 
\\
it's a regular array of motionless atoms.
\\ \\
$H$(Gas) is very high because 
you have to specify 
\\
where each atom is and where it is going
and how fast.
\\
\\
\hline
\end{tabular}
\end{center}
 
Here, however, is a (very theoretical) biological application of AIT, using mutual information,
the extent to which two things considered together are simpler than considered apart.
In other words,
this is the extent to which the smallest program that computes \textbf{both} simultaneously 
is smaller than the sum of the size of the individual smallest programs for each separately.
 
\begin{center}
\begin{tabular}{|c|}
\hline
\\
\textbf{\large \emph{Mutual Information Between X and Y}}
\\ \\
The mutual information is equal to $H(X) + H(Y) - H(X, Y)$.
\\ \\
It's \textbf{small} if 
$H(X, Y)$ is approximately equal to $H(X) + H(Y)$.
\\ \\
It's \textbf{large} if 
$H(X, Y)$ is much smaller than $H(X) + H(Y)$.
\\ \\
Note that $H(X, Y)$ \textbf{cannot} be larger than $H(X) + H(Y)$
\\
because programs are self-delimiting binary information.
\\
\\
\hline
\end{tabular}
\end{center}
 
If two things have very little in common, it makes no difference if we consider them together
or separately.  But if they have a lot in common, then we can eliminate common subroutines
when we calculate them at the same time, and so the mutual information will be large.
How about biology?  Well, there is the problem of the whole versus the parts.  By what right
can we partition the world of our experience into parts, instead of considering it as a unified whole?
I'll discuss this in a moment.
Mutual information also has a very theoretical application in music theory, 
not just in theoretical biology.
How?
 
In music, we can use this measure of mutual information 
to see how close two compositions are, to see how much they have in common.
Presumably, two works by Bach will have higher mutual information than a work by Bach and
one by Shostakovich. Or compare the entire body of work of two composers.  Two Baroque
composers should have more in common than a Baroque and a Romantic composer.
 
\begin{center}
\begin{tabular}{|c|}
\hline
\\
\textbf{\large \emph{Bach, Shostakovich}}
\\ \\
$H$(Bach, Shostakovich) is approximately equal to 
\\
$H$(Bach) + $H$(Shostakovich).
\\ \\
\textbf{\large \emph{Bach$_1$ = Brandenburg Concertos}}
\\
\textbf{\large \emph{Bach$_2$ = Art of the Fugue}}
\\ \\
$H$(Bach$_1$, Bach$_2$) is much less than $H$(Bach$_1$) + \emph{H}(Bach$_2$).
\\
\\
\hline
\end{tabular}
\end{center}
 
How about biology?  As I said, there is the problem of the whole versus the parts.  By what right
can we partition the universe into parts, rather than consider it as an interacting whole?
Well, the parts of an organism have high mutual information.  
(Because they all contain the full genome for that organism, they all have the same DNA, even
though different cells in different tissues turn on different parts of that DNA.)
Also, it is natural to divide
a whole into pieces if the program-size complexity decomposes additively, that is, if the
complexity of the whole is approximately equal to the sum of the complexity of those parts, 
which means
that their mutual interactions are not as important as their internal interactions.  This
certainly applies to living beings!
 
\begin{center}
\begin{tabular}{|c|}
\hline
\\
\textbf{\large \emph{Fred, Alice}}
\\ \\
\emph{H}(Fred, Alice) is approximately equal to $H$(Fred) + $H$(Alice).
\\ \\
\textbf{\large \emph{Fred's left arm, Fred's right arm}}
\\ \\
$H$(left arm, right arm) is much less than $H$(left arm) + $H$(right arm).
\\
\\
\hline
\end{tabular}
\end{center}
 
But do we really have the right to talk about the complexity of a physical object
like $H$(Alice), $H$(Fred)?!
Complexity of digital objects, yes, they're all just finite strings of 0's and 1's.
But physicists usually use \textbf{real numbers} to describe the physical world, 
and they assume that space and time and many measurable quantities can change continuously.
And a real number isn't a digital object, it's an analog object.
And since it varies continuously, not in discrete jumps, if you convert such a number
into bits you get an
infinite number of bits.
But computers can't perform computations with numbers that have an infinite number of bits!
And my theory of algorithmic information is based on what computers can do.
 
Does this pull the rug out from under everything? Fortunately it doesn't.
 
In Chapter IV we'll discuss whether the physical world is really continuous or whether
it might actually be discrete, like some rebels are beginning to suspect.
We'll see that there are in fact plenty of physical reasons to reject real numbers. 
And in 
Chapter V we'll discuss mathematical and philosophical reasons to reject real
numbers.
There are lots of arguments against real numbers, 
it's just that people usually aren't willing to listen to any of them!
 
So our digital approach may in fact have a truly broad range of applicability.
After this reassuring vindication of our digital viewpoint, the road will finally be clear
for us to get to $\Omega$, which will happen in Chapter VI.

\chapter*{Chapter IV---Intermezzo}
\addcontentsline{toc}{chapter}{Chapter IV---Intermezzo}
\markright{Chapter IV---Intermezzo}

\section*{The Parable of the Rose}

Let's take a break.  Please take a look at this magnificent rose:
\vspace{7cm}
 
In his wonderfully philosophical short story 
\emph{The Rose of Paracelsus} (1983), Borges depicts the
medieval alchemist Paracelsus, who had the fame of being able to recreate
a rose from its ashes:
\begin{quotation}
The young man raised the rose into the air.
 
``You are famed,'' he said, ``for being able to burn a rose to ashes
and make it emerge again, by the magic of your art. Let me 
witness that prodigy. I ask that of you, and in return I will offer up
my entire life.''  \ldots
 
``There is still some fire there,'' said Paracelsus, pointing toward 
the hearth. ``If you cast this rose into the embers, you would
believe that it has been consumed, and that its ashes are real. I
tell you that the rose is eternal, and that only its appearances may
change. At a word from me, you would see it again.''
 
``A word?'' the disciple asked, puzzled. ``The furnace is cold, and
the retorts are covered with dust. What is it you would do to bring
it back again?''\footnote
{\textbf{Source}: Jorge Luis Borges, 
\emph{Collected Fictions,} translated by Andrew Hurley,
Penguin Books, 1999, pp.\ 504--507.
For the original \emph{La rosa de Paracelso,}
see Borges, 
\emph{Obras Completas,} Tomo III, Emec\'e, 1996, pp.\ 387--390.}
\end{quotation}
 
Now let's translate this into the language of our theory.
 
The \emph{algorithmic information content} (program-size complexity) $H$(rose) of a rose 
is defined to be
the size in bits of the smallest computer program (algorithm) $p_{\mbox{\scriptsize rose}}$ that 
produces the digital image of this rose.
 
This minimum-size algorithmic description $p_{\mbox{\scriptsize rose}}$ 
captures the \textbf{irreducible essence}
of the rose,
and is the number of bits that you need to preserve in order to be able to recover the rose:

\begin{center}
\fbox{
ashes $\longrightarrow$
\textbf{\large Alchemy/Paracelsus}
$\longrightarrow$ rose
}
\vspace{\baselineskip}
\\ 
\fbox{
minimum-size program $p_{\mbox{\scriptsize rose}}$ $\longrightarrow$
\textbf{\large Computer}
$\longrightarrow$ rose 
}
\vspace{\baselineskip}
\\
size in bits of $p_{\mbox{\scriptsize rose}} =$ 
\vspace{.25\baselineskip}
\\
\fbox{
algorithmic information
content $H$(rose) of the rose 
}
\end{center}
 
$H$(rose) measures the \textbf{conceptual complexity} of the rose,
that is, the difficulty (in bits) of conceiving of the rose,
the number of \textbf{independent yes/no choices} that must be made in order to do this.
The larger $H$(rose) is, the \textbf{less inevitable} it is for God to
create this particular rose, as an independent act of creation.
 
In fact, the work of Edward Fredkin, Tommaso Toffoli and Norman Margolus
on reversible cellular automata, has shown that there are
discrete toy universes in which \textbf{no information is
ever lost}, that is, speaking in medieval terms, \textbf{the soul is immortal}.  
Of course, this is not a terribly personal kind of immortality; it only
means that enough information is always present to be able to reverse time and recover 
any previous state of the universe.
 
Furthermore
the rose that we have been considering is only a jpeg image, it's
only a digital rose. What about a real rose? Is that analog or digital?
Is physical reality analog or digital?
In other words,
in this universe
are there infinite precision real numbers, or is everything built up out of a finite number
of 0's and 1's?
 
Well, let's hear what Richard Feynman has to say about this.
\begin{quotation}
``It always bothers me that, according to the laws as we understand them today,
it takes a computing machine an infinite number of logical operations to
figure out what goes on in no matter how tiny a region of space, and no
matter how tiny a region of time. How can all that be going on in that tiny space?
Why should it take an infinite amount of logic to figure out what one tiny
piece of space/time is going to do?
So I have often made the hypothesis that ultimately physics will not require
a mathematical statement, that in the end the machinery will be revealed,
and the laws will turn out to be simple like the chequer board\ldots''---Feynman,
\emph{The Character of Physical Law}.
\end{quotation}
Fredkin told me that he spent years trying to make Feynman
take digital physics seriously, 
and was very pleased to see his ideas reproduced in this passage in Feynman's book!
For an idea to be successful, you have to give it away, you have to 
be willing to 
let other people
think that it was theirs!
You can't be possessive, you can't be jealous!\ldots
 
I used to be fascinated by physics, and I wanted very badly to be able to know if physical
reality is discrete or continuous.
 
But what I really want to know now is:
``What is life?'', ``What is mind?'', ``What is intelligence?'', ``What is consciousness?'',
``What is creativity?''.
By working in a toy world and using digital physics, we can hopefully avoid
getting bogged down in the
physics and we can concentrate on the biology instead.
 
I call this theoretical theoretical physics, because even if these models do not
apply to \textbf{this particular} world, they are \textbf{interesting possible} worlds!
 
And my ultimate ambition, which hopefully somebody will someday achieve, would be
to \textbf{prove} that life, intelligence and consciousness must with high probability
evolve in one of these toy worlds.  In playing God like this, it's important to
get out more than you put in, because otherwise the whole thing can just be a cheat.  
So what would really be nice is to be able to obtain \textbf{more} intelligence, or
a higher degree of consciousness, than you yourself possess!
You, who set up the rules of the toy universe in the first place.
 
This is certainly not an easy task to achieve, but I think that this 
kind of higher-level understanding
is probably more fundamental than just disentangling the microstructure of the physical
laws of this one world.
     
\section*{Theoretical Physics \& Digital Philosophy}

Nevertheless, there are some intriguing hints that this particular universe
may in fact be a discrete digital universe, not a continuous analog universe
the way that most people would expect.
 
In fact these ideas actually go back to Democritus, who argued that matter must be
discrete, and to Zeno, who even had the audacity to suggest that continuous
space and time were self-contradictory impossibilities.
 
Through the years I've noticed many times, as an armchair physicist, places where
physical calculations diverge to infinity at extremely small distances. Physicists are adept
at not asking the wrong question, one that gives an infinite answer.  
But I'm a mathematician, and each time I would
wonder if Nature wasn't really trying to tell something, namely that real numbers and continuity are
a sham, and that infinitesimally small distances \textbf{do not exist!}
 
\emph{Two Examples:} the infinite amount of energy stored in the field around 
a point electron according to Maxwell's theory of electromagnetism,
and the infinite energy content of the vacuum according to quantum field theory.
 
In fact,
according to quantum mechanics, infinitely precise measurements require infinite energies (infinitely
big and expensive atom smashers), but long before that you get gravitational collapse
into a black hole, if you believe in general relativity theory.
 
For more on ideas like this, see the discussion of the 
Bekenstein bound and the holographic principle
in Lee Smolin, \emph{Three Roads to Quantum Gravity.}
 
And although it's not usually presented that way,
string theory was invented to eliminate these divergences as distances
become arbitrarily small. String theory does this by eliminating point field sources
like the electron, which as I pointed out 
causes difficulties for Maxwell's theory of electromagnetism.
In string theory elementary particles are changed from points into loops of string.
And string theory
provides a crucial \textbf{minimum distance scale}, which is the length of these loops of string.
(They're in a tightly curled-up extra dimension.)
So you can no longer make distances go below that minimum length scale.
 
That's one of the reasons that string theory had to be invented, in order to get rid of
arbitrarily small distances!
 
Last but not least, I should share with you the following experience:
For many years, whenever I would bump into my esteemed late colleague Rolf Landauer
in the halls of the IBM Watson Research Center in Yorktown Heights, he would often gently
assure me that
no physical quantity has ever been measured 
with more than, say, 20 digits of precision.
And the experiments that achieve that degree of precision are masterpieces of laboratory
art and the sympathetic search for physical quantities stable enough and sharply defined enough
\textbf{to allow themselves} to be measured with that degree of precision, not to mention eliminating
all possible perturbing sources of noise!
So why, Rolf would ask, should \textbf{anyone} believe in arbitrarily precise measurements 
or in the real numbers used to record the results of such measurements?!
 
So, you see, there are lots of reasons to suspect that we might live
in a digital universe, that God prefers to be able to copy things exactly
when he has to, rather than to get the inevitable increase in noise that
accompanies analog copying!
 
And in the next chapter I'd like to continue in Zeno's footsteps as well as Rolf's and argue
that \textbf{a number with infinite precision, a so-called real number, is actually
rather unreal!}

\chapter*{Chapter V---The Labyrinth of the Continuum}
\addcontentsline{toc}{chapter}{Chapter V---The Labyrinth of the Continuum}
\markright{Chapter V---The Labyrinth of the Continuum}

The ``labyrinth of the continuum'' is how Leibniz referred to the philosophical problems
associated with real numbers, which we shall discuss in this chapter.
So the emphasis here will be on philosophy and mathematics, rather than on physics as
in the last chapter.
 
What is a real number? Well, in geometry
it's the length of a line, measured exactly, with infinite
precision, for example 1.2749591\ldots, which doesn't sound too problematical, 
at least at first.
And in analytic
geometry you need \textbf{two} real numbers to locate a point (in two dimensions), its
distance from the $x$ axis, and its distance from the $y$ axis.
\textbf{One} real will locate a point on a line,
and the line that we will normally consider will be the so-called ``unit interval''
consisting of all the real numbers from zero to one.
Mathematicians write this interval as $[0,1)$,  to indicate that 0 is included but 1 is not,
so that all the real numbers corresponding to these points have no integer part,
only a decimal fraction.
Actually, $[0,1]$ works too, as long as you write one as 0.99999\ldots\ instead of as 1.00000\ldots\
But not to worry, we're going to ignore all these subtle details. You get the general 
idea, and that's enough for reading this chapter.
 
{\footnotesize
[By the way,
why is it called ``real''? To distinguish it from so-called ``imaginary'' numbers like $\sqrt{-1}$.
Imaginary numbers are neither more nor less imaginary than real numbers, but there was initially,
several centuries ago, including at the time of Leibniz, much resistance to placing them on 
an equal footing with real numbers.  In a letter to Huygens, Leibniz points out that calculations
that temporarily traverse this imaginary world can in fact start and end with real numbers.
The usefulness of such a procedure was, he argued, an argument in favor of such numbers.
By the time of Euler, imaginaries were extremely useful. For example, Euler's famous result
that 
\[
   e^{i x} = \cos x + i \sin x 
\] 
totally tamed trigonometry.
And the statement (Gauss)
that an algebraic equation of degree $n$ has exactly $n$ roots only works with the 
aid of imaginaries. Furthermore, the theory of functions of a complex variable (Cauchy)
shows that the calculus and in particular so-called power series
\[
   a_0 + a_1 x + a_2 x^2 + a_3 x^3 + \ldots
\] 
make much more sense with imaginaries than without.
The final argument in their favor, if any was needed, was provided by Schr\"odinger's equation
in quantum mechanics, in which imaginaries are absolutely essential, since quantum probabilities
(so-called ``probability amplitudes'') have to have direction as well as magnitude.]
}
 
As is discussed in Burbage and Chouchan, \emph{Leibniz et l'infini,} PUF, 1993,
Leibniz referred to what we call the infinitesimal calculus as ``the calculus of
transcendentals''. And he called curves ``transcendental'' if they cannot be obtained via
an algebraic equation, the way that the circles, ellipses, parabolas and hyperbolas
of analytic geometry most certainly can.  
 
Leibniz was extremely proud of his quadrature of the circle, a problem that had eluded
the ancient Greeks, but that he could solve with \emph{transcendental} methods:
\[
 \pi/4 \; = \; 1 - 1/3 + 1/5 - 1/7 + 1/9 - 1/11 + 1/13 + \ldots
\] 
What is the quadrature of the circle? 
The problem is to geometrically construct a square having the same area as a given circle, that
is, to determine the area of the circle. Well, that's $\pi r^2$, $r$ being the radius
of the circle, which converts the problem into determining $\pi$, precisely what Leibniz accomplished
so elegantly with the infinite series displayed above.
 
\textbf{
Leibniz could not have failed to be aware that in using this term he was evoking
the notion of God's transcendence of all things human, of human limitations, of
human finiteness.
}
 
As often happens, history has thrown away the philosophical ideas that inspired
the creators and kept only a dry technical husk of what they thought that they had achieved.
What remains of Leibniz's idea of transcendental methods is merely the distinction between
algebraic numbers and transcendental numbers.
A real number $x$ is algebraic if it is the solution of an equation of the form
\[
   a x^n + b x^{n-1} + \cdots + p x + q = 0
\] 
where the constants $a, b, \ldots$ are all integers. Otherwise $x$ is said to be transcendental.
The history of proofs of the existence of transcendental numbers is rich in intellectual
drama, and is one of the themes of this chapter.
 
Similarly, it was Cantor's obsession with God's infiniteness and transcendence 
that led him to create his spectacularly successful 
but extremely controversial
theory of infinite sets and infinite numbers.
What began, at least in Cantor's mind, as a kind of madness, 
as a kind of mathematical theology full---necessarily full---of paradoxes, 
such as the one discovered by Bertrand Russell, 
since any attempt by a finite mind to apprehend God is inherently paradoxical,
has now been
condensed and desiccated into an extremely technical and untheological field of math,
modern axiomatic set theory.
 
Nevertheless, the intellectual history of the proofs of the existence of transcendental numbers
is quite fascinating. New ideas totally transformed
our way of viewing the problem, not once, but in fact \textbf{four} times!
Here is an outline of these developments:
\begin{itemize}
\item
Liouville, Hermite and Lindemann, with great effort, were the first to exhibit individual
real numbers that could be proved to be transcendental. 
Summary: \textbf{individual transcendentals.}
\item
Then Cantor's theory of infinite sets revealed that the transcendental
reals had the same cardinality as the set of all reals, while the algebraic reals were merely
as numerous as the integers, a smaller infinity.
Summary: \textbf{most reals are transcendental.}
\item
Next Turing pointed out that all algebraic reals are computable, but again, the uncomputable
reals are as numerous as the set of all reals, while the computable reals are only as
numerous as the integers. The existence of transcendentals is
an immediate corollary.
Summary: \textbf{most reals are uncomputable and \emph{therefore} transcendental.}
\item
The next great leap forward involves probabilistic ideas: the set of random
reals was defined, and it turns out that with probability one, 
a real number is random and therefore necessarily
uncomputable and transcendental. Non-random, computable and algebraic reals
all have probability zero.  So now you can get a transcendental real merely by picking a real
number at random with an infinitely sharp pin, or, 
alternatively, by using independent tosses of a fair
coin to get its binary expansion.
Summary: \textbf{reals are transcendental/uncomputable/random with probability one.}
And in the next chapter 
we'll exhibit a natural construction that picks out an individual random real, namely
the halting probability $\Omega$,
without the need for an infinitely sharp pin.
\item
Finally, and perhaps even more devastatingly, 
it turns out that the set of all reals that can be individually named or specified
or even defined or referred to---constructively or not---within a formal language 
or within an individual FAS, has probability zero.
Summary: \textbf{reals are un-nameable with probability one.}
\end{itemize}
So the set of real numbers, 
while natural---indeed, immediately given---geometrically, nevertheless remains quite elusive:
 
\textbf{
Why should I believe in a real number if I can't calculate it, if I can't prove
what its bits are, and if I can't even refer to it? 
And each of these things happens with probability one!
The real line from 0 to 1 looks more and more like a Swiss cheese, 
more and more
like a stunningly black high-mountain
sky studded with pin-pricks of light.
}
 
Let's now set to work to explore these ideas in more detail.

\section*{The ``Unutterable'' and the Pythagorean School}

This intellectual journey actually begins, as is often the case, with the ancient Greeks. 
Pythagoras is credited with naming both mathematics and philosophy.
And the Pythagoreans believed that number---whole numbers---rule the universe, 
and that God is a mathematician,
a point of view largely vindicated by modern science, especially quantum mechanics,  
in which the hydrogen atom is modeled as a musical instrument that produces a discrete scale of notes.
Although, as we saw in Chapter III, perhaps God is actually a computer programmer!
 
Be that as it may, these early efforts to understand the universe suffered
a serious setback when the 
Pythagoreans discovered geometrical lengths that cannot be expressed as the
ratio of two whole numbers. Such lengths are called irrational or incommensurable.
In other words, they discovered real numbers that cannot be expressed as a ratio of 
two whole numbers.
 
How did this happen?
 
The Pythagoreans considered the unit square, 
a square one unit in length on each side, and they discovered that
the size of both of the two diagonals, $\sqrt{2}$,
isn't a rational number $n/m$. That is to say, it cannot
be expressed as the ratio of two integers. 
In other words, there are no integers $n$ and $m$ such that
\[
   (n/m)^2 = 2 \;\; \mbox{or} \;\; n^2 = 2 m^2
\] 
An elementary proof of this from first principles is given in Hardy's 
well-known
\emph{A Mathematician's
Apology.} 
He presents it there because
he believes that it's
a mathematical argument whose beauty anyone should be able to appreciate.
However, the proof that Hardy gives, which is actually from Euclid's \emph{Elements,} 
does not give as much insight as a more advanced proof using
unique factorization into primes. I \textbf{did not} prove unique factorization in Chapter II.
Nevertheless, I'll use it here.
It is relevant because the two sides of the equation 
\[
   n^{2} = 2 m^{2}
\] 
would give us \textbf{two different}
factorizations of the same number. How?
 
Well, factor $n$ into primes, and factor $m$ into primes. 
By doubling the exponent of each prime in the 
factorizations of $n$ and $m$,
\[
 2^{\alpha} \, 3^{\beta}\,  5^{\gamma} \ldots
 \; \longrightarrow \;
 2^{2\alpha} \, 3^{2\beta} \, 5^{2\gamma} \ldots,
\] 
we get factorizations of
$n^2$ and $m^2$. 
This gives us a factorization of
$n^2$ in which
the exponent of 2 is even, and 
a factorization of
$2 m^2$ 
in which the exponent of 2 is odd.
So we have two different factorizations of the same number into primes, which is impossible.
 
According to Dantzig, \emph{Number, The Language of Science,}
the  discovery of irrational or incommensurable numbers like $\sqrt{2}$ 
\begin{quote}
``caused great
consternation in the ranks of the Pythagoreans. The very name given to these entities
testifies to that. \emph{Algon,} the \emph{unutterable,} these incommensurables were called\ldots\
How can number dominate the universe when it fails to account even for the most
immediate aspect of the universe, namely \emph{geometry?}
So ended the first attempt to exhaust nature by number.''
\end{quote}
This 
intellectual history also left its traces in the English language:
In English such irrationals are referred to as ``surds'', which comes from
the French ``sourd-muet'', meaning deaf-mute, one who cannot hear or speak. 
So the English word ``surd'' comes from the French word for ``deaf-mute'', and algon = mute.
In Spanish it's ``sordomudo'', deaf-mute.
 
In this chapter we'll retrace this history, and we'll see that real numbers
not only confound the philosophy of Pythagoras, they confound as well Hilbert's belief in the 
notion of a FAS,
and they provide us with many additional reasons for doubting their existence, and for
remaining quite skeptical. To put it bluntly, our purpose here is to review and discuss
the mathematical arguments \textbf{against real numbers}.

\section*{The 1800's: Individual transcendentals (Liouville, Hermite, Lindemann)}

Although Leibniz was extremely proud of the fact that 
he been able to square the circle using transcendental methods,
the 1800's wanted to be \textbf{sure} that they were really required. In other words, they demanded
proofs that $\pi$ and other individual numbers defined via the sums of infinite series
\textbf{were not} the solution of any algebraic equation.
 
Finding a natural specific example of a transcendental real turned out to be much harder
than expected. It took great ingenuity and cleverness to exhibit provably transcendental numbers!
 
The first such number was found by Liouville:
\[
   \mbox{Liouville number} =  1/10^{1!} + 1/10^{2!} + \cdots + 1/10^{n!} + \ldots
\]
He showed that algebraic numbers cannot be approximated that well by rational numbers.
In other words, he showed that his number cannot be algebraic, because 
there are rational approximations
that work too well for it: they can get too close, too fast.
But Liouville's number isn't a natural example, because no one had ever
 been interested in this particular
number before Liouville.  It was constructed precisely so that Liouville could prove its
transcendence. 
What about
$\pi$ and Euler's number $e$?
 
Euler's number
\[
   e = 1/1! + 1/2! + 1/3! + \cdots + 1/n! + \ldots
\] 
was finally proved transcendental by Hermite.  Here at last was a natural example!
This was an important number that people really cared about!
 
But what about the number that Leibniz was so proud of conquering?  He had squared the circle
by transcendental methods:
\[
   \pi/4 = 1/1 - 1/3 + 1/5 - 1/7 + 1/9 - 1/11 + \ldots
\] 
But can you prove that transcendental methods are really necessary?  
This question attracted a great deal of attention after Hermite's result, since $\pi$ seemed
to be the obvious next candidate for a transcendence proof.
This feat was finally accomplished by
Lindemann, provoking the famous remark by Kronecker that ``Of what use is your beautiful
proof, since \textbf{$\pi$ does not exist!''}  Kronecker was a follower of Pythagoras;
Kronecker's best known statement is,
``God created the integers; all the rest is the work of man!''
 
These were the first steps
on the long road to understanding transcendence,
but they were difficult complicated proofs that were specially
tailored for each of these specific numbers, and gave no general insight into what was going on.

\section*{Cantor: The number of transcendentals is a higher order infinity than the number of algebraic reals}

As I've said, a real number is one that can be determined with arbitrary precision, such as
$\pi = 3.1415926\ldots$
Nevertheless,
in the late 1800's two mathematicians, Cantor and Dedekind, 
were moved to come up with much more careful
definitions of a real number.
Dedekind did it via ``cuts'', thinking of an irrational real $r$ as a way to partition all the rational
numbers $n/m$ into those less than $r$ and those greater than $r$.
In Cantor's case a real was defined as an infinite sequence of rational numbers $n/m$
that approach $r$ more and more closely.\footnote
{Earlier versions of the work of Dedekind and of Cantor
on the reals are due to Eudoxus and to Cauchy, respectively. 
History repeats itself, even in mathematics.}

History did not take any more kindly to their work than it has to any other attempt at a
``final solution''.
 
But first, let me tell you about Cantor's theory of infinite sets and his
invention of new, infinite numbers for the purpose of measuring the sizes of all infinite sets.
A very bold theory, indeed!
 
Cantor's starting point is his notion of comparing two sets, finite or infinite,
by asking whether or not there is a one-to-one correspondence, a pairing between the
elements of the two sets that exhausts both
sets and leaves no element of either set unpaired, and no element of one of the sets
paired with more than one partner in the other set.
If this can be done, then Cantor declares that the two sets are equally big.
 
Actually, Galileo had mentioned this idea in one of his dialogues, the one published
at the end of his life when he was under house arrest. Galileo points out that there
are precisely as many positive integers $1, 2, 3, 4, 5, \ldots$ as there are square numbers
$1, 4, 9, 16, 25, \ldots$  Up to that point, history has decided that Galileo was right on
target. 
 
However, he then declares that the fact that the squares are just a tiny fraction of
all the positive integers contradicts his previous observation that they are equally
numerous, and that this paradox
precludes making any sense of the notion of the size of an infinite set.
 
The paradox of the whole being equivalent to one of its parts, may have deterred
Galileo, but Cantor and Dedekind took it entirely in stride.
It did not deter them at all. In fact, Dedekind even 
put it to work for him, he used it.
Dedekind \textbf{defined} an infinite set to be one having
the property that
a proper subset of it is just as numerous as it is!
In other words, according to Dedekind,
a set is infinite if and only if it can be put in a one-to-one correspondence
with a part of itself, one that excludes some of the elements of the original set!
 
Meanwhile, Dedekind's friend Cantor
was starting to apply this new way of comparing the size of
two infinite sets to common everyday mathematical objects: integers, rational numbers, algebraic
numbers, reals, points on a line, points in the plane, etc.
 
Most of the well-known mathematical objects broke into two classes: 1) sets like
the algebraic real numbers and the rational numbers, which were exactly as numerous
as the positive integers, and are therefore called ``countable'' or ``denumerable'' infinities,
and 2) sets like the points in a finite or infinite line or in the plane or in space,
which turned out all to be exactly as numerous as each other, and which are said to 
``have the power of the continuum''.
And this gave rise to two new infinite numbers,
$\aleph_0$ (aleph-nought) and \textbf{\emph{c}}, both invented by Cantor,
that are, respectively, the size (or as Cantor called it, the ``power''
or the ``cardinality'') of the positive integers and of the continuum of real numbers.
 
\begin{center}
\begin{tabular}{|c|}
\hline
\\
\textbf{\emph{\large Comparing Infinities!}}
\\ \\
\#\{reals\} = \#\{points in line\} = \#\{points in plane\} = \textbf{\emph{c}}
\\ \\
\#\{positive integers\} = \#\{rational numbers\} = 
\\
\#\{algebraic real numbers\} = $\aleph_0$
\\
\\
\hline
\end{tabular}
\end{center}
 
Regarding his proof that there were precisely as many points in a plane as 
there are in a solid or in a line,
Cantor remarked in a letter to Dedekind, ``Je le vois, mais je ne le crois pas!'',
which means ``I see it, but I don't believe it!'', and which happens to have a pleasant melody 
in French.
 
And then Cantor was able to prove the extremely important and basic theorem that
\textbf{\emph{c}} is larger than $\aleph_0$, 
that is to say, that the continuum is a nondenumerable infinity, 
an uncountable infinity, in other words,
that there are
more real numbers than there are positive integers, infinitely more.
This he did by using Cantor's well-known diagonal method, 
explained in Wallace, \emph{Everything and More,}
which is all about Cantor and his theory.
 
In fact, it turns out that
the infinity of transcendental reals is exactly as large as the infinity of all reals,
and the smaller infinity of algebraic reals is exactly as large as the infinity of positive integers.
Immediate corollary: most reals are transcendental, not algebraic, infinitely more so.
 
Well, this is like stealing candy from a baby!
It's
much less work than struggling with individual real numbers and trying to prove that they are
transcendental!
Cantor gives us a much more general perspective from which to view this particular problem.
And it's much easier to see that \textbf{most} reals are transcendental 
than to decide if a \textbf{particular}
real number happens to be transcendental!
 
So that's the first of what I would call the ``philosophical'' proofs that transcendentals exist.
Philosophical as opposed to highly-technical, like flying by helicopter to the top of the
Eiger instead of reaching the summit by 
climbing up its infamous snow-covered north face.
 
Is it really that easy?  Yes, but this set-theoretic approach 
created as many problems as it solved.  
The most famous is called Cantor's continuum problem.
 
What is Cantor's continuum problem?
 
Well, it's the question of whether or not there happens to be any set that has
more elements than there are positive integers, 
and that has fewer elements than there are real numbers.
In other words, is there an infinite set whose cardinality or power 
is bigger than $\aleph_0$ 
and smaller than \textbf{\emph{c}}?  In other words, is \textbf{\emph{c}} the next
infinite number after $\aleph_0$, which has the name 
$\aleph_1$
(aleph-one) reserved for it in Cantor's theory,
or are there a lot of other aleph numbers in between?
 
\begin{center}
\begin{tabular}{|c|}
\hline
\\
\textbf{\emph{\large Cantor's Continuum Problem}}
\\ \\
Is there a set $S$ such that $\aleph_0 < \#S <$ \textbf{\emph{c}}?
\\ \\
In other words, is \textbf{\emph{c}} = $\aleph_1$, which is
\\
the first cardinal number after $\aleph_0$?
\\
\\
\hline
\end{tabular}
\end{center}
 
A century of work has not sufficed to solve this problem!
 
An important milestone was the proof by the combined efforts of
G\"odel and Paul Cohen that the
usual axioms of axiomatic set theory (as opposed to the ``naive''
paradoxical original Cantorian set theory), do 
not suffice to decide one way or another. You can add
a new axiom asserting there is a set with intermediate power, or that there is no
such set, and the resulting system of axioms will not lead to a contradiction
(unless there was already one there, without even having to use this new axiom, 
which everyone fervently hopes is not the case).
 
Since then there has been a great deal of work to see if there might be new axioms
that set theorists can agree on that might enable them to settle Cantor's continuum problem.
And indeed, something called the axiom of projective determinacy has become quite
popular among set theorists, since it permits them to solve many open problems that interest them.
However, it doesn't suffice to settle the continuum problem!
 
So you see, the continuum refuses to be tamed!
 
And now we'll see how the real numbers, annoyed at being ``defined''
by Cantor and Dedekind, got their revenge in the century after Cantor,
the 20th century.

\section*{Borel's amazing know-it-all real number}

The first intimation that there might be something wrong, something terribly wrong,
with the notion of a real number comes from a small paper published by \'Emile Borel in 1927.
 
Borel pointed out that if you really believe in the notion of a real number as
an infinite sequence of digits 3.1415926\ldots, then you could put all of human knowledge
into a single real number.  Well, that's not too difficult to do, that's only a finite
amount of information. You just take your favorite encyclopedia, for example, the 
\emph{Encyclopedia Britannica}, which I used to use when I was in high-school---we had a nice
library at the Bronx High School of Science---and you digitize it, you convert it into
binary, and you use that binary as the base-two expansion of a real number in the unit
interval between zero and one!
 
So that's pretty straight-forward, especially now that most information, including books,
is prepared in digital form before being printed.
 
But what's more amazing is that there's nothing to stop us from putting an infinite amount
of information into a real number.
In fact, there's a single real number, I'll call it Borel's number, since he imagined it, in 1927,
that can serve as an oracle and answer any yes/no question that we could ever pose to it.
How? Well, you just number all the possible questions, and then the $N$th digit or $N$th
bit of Borel's number tells you whether the answer is yes or no!
 
If you could come up with a list of all possible yes/no questions and only valid yes/no questions, 
then Borel's number
could give us the answer in its binary digits.
But it's hard to do that.  It's much easier to simply list all possible texts in the
English language (and Borel did it using the French language), all possible finite strings of
characters that you can form using the English alphabet, including a blank for use between
words.  You start with all the one-character strings, then all the two-character strings, etc.
And you number them all like that\ldots
 
Then you can use the $N$th digit of Borel's number to tell you whether the $N$th string of characters is
a valid text in English, then whether it's a yes/no question, then whether it has an answer, 
then whether the answer is yes or no.
For example, ``Is the answer to this question `No'?'' looks like a valid yes/no question, but in
fact has no answer.
 
So we can use $N$th digit 0 to mean bad English, 1 to mean not a yes/no question,
2 to mean unanswerable, and 3 and 4 to mean ``yes'' and ``no'' are the answers, respectively.
Then 0 will be the most common digit, 
then 1, then there'll be about as many 3's as 4's, and, I expect,
a smattering of 2's.
 
Now Borel raises the extremely troubling question, ``Why should we believe in this real number
that answers every possible yes/no question?''
And his answer is that he doesn't see any reason to believe in it, none at all!
According to Borel, this number is merely a mathematical fantasy, a joke, a 
\emph{reductio ad absurdum}
of the concept of a real number!
 
You see, some mathematicians have what's called a ``constructive'' attitude. 
This means that they only believe in mathematical objects that can be constructed, that,
given enough time, in theory one could actually calculate.
They think that there ought to be some
way to \textbf{calculate} a real number, 
to calculate it digit by digit, otherwise in what sense can it be said to have some kind of
mathematical existence?
 
And this is precisely the question discussed by Alan Turing in his famous 1936 paper
that invented the computer as a mathematical concept.  
He showed that there were lots and lots of computable real numbers.
That's the positive part of his paper.
The negative part, is that he also showed that there were lots and lots of uncomputable real numbers.
And that gives us another philosophical proof that there are transcendental numbers, because
it turns out that all algebraic reals are in fact computable.

\section*{Turing: Uncomputable reals are transcendental}

Turing's argument is very simple, very Cantorian in flavor.
First he invents a computer (on paper, as a mathematical idea, a model computer).
Then he points out that the set of all possible computer programs is a countable set,
just like the set of all possible English texts.
Therefore the set of all possible computable real numbers must also be countable.
But the set of all reals is uncountable, it has the power of the continuum.
Therefore the set of all uncomputable reals is also uncountable and has the power of
the continuum.  Therefore most reals are uncomputable, infinitely more than are computable.
 
That's remarkably simple, if you believe in the idea of a general-purpose digital
computer. Now we are all very familiar with that idea.  Turing's paper is long precisely
because that was not at all the case in 1936. So he had to work out a simple computer
on paper and argue that it could compute anything that can ever be computed, before giving
the above argument that most real numbers will then be uncomputable, in the sense
that there cannot be a program for computing them digit by digit forever.
 
The other difficult thing is to work out in detail precisely why algebraic reals can
be computed digit by digit.
Well, it's sort of intuitively obvious that this has to be the case; after all, what
could possibly go wrong?!
In fact, this is now well-known technology using something called Sturm sequences,
that's the slickest way to do this; I'm sure that
it comes built into \emph{Mathematica} and \emph{Maple,}
two symbolic computing software packages.
So you can use these software packages to calculate as many digits as you want.
And you need to be able to calculate hundreds of digits in
order to do research
the way described by Jonathan Borwein and David Bailey 
in their book \emph{Mathematics by Experiment.}
 
But in his 1936 paper Turing mentions a way to calculate algebraic reals that will
work for a lot of them, and since it's a nice idea, I thought I'd tell you about it.
It's a technique for root-solving by successive interval halving.
 
Let's write the algebraic equation that determines an individual algebraic real $r$ that 
we are interested in as $\phi(x) = 0$; $\phi(x)$ is a polynomial in $x$.
So $\phi(r) = 0$, and let's suppose we know two rational numbers $\alpha, \beta$ such that
$\alpha < r < \beta$ and $\phi(\alpha) < \phi(r) < \phi(\beta)$
and we also know that there is no other root of the equation  $\phi(x) = 0$ in that interval.
So the signs of 
$\phi(\alpha)$ and $\phi(\beta)$
have to be different, neither of them is zero, and precisely one of them is greater than
zero and one of them is less than zero, that's key.
Because if $\phi$ changes from positive to negative it must pass through zero somewhere in between.
 
Then you just bisect this interval $[\alpha,\beta]$.  You look at the midpoint
$(\alpha+\beta)/2$ which is also a rational number, and you plug that into $\phi$ and
you see whether or not $\phi((\alpha+\beta)/2)$ is equal to zero, less than zero,
or greater than zero. It's easy to see which, since you're only dealing with rational numbers,
not with real numbers, which have an infinite number of digits.
 
Then if $\phi$ of the midpoint gives zero, we have found $r$ and we're finished.
If not, we choose the left half or the right half of our original interval in such
a way that the sign of $\phi$ at both ends is different, and this new interval
replaces our original interval, $r$ must be there, and we keep on going like that forever.
And that gives us better and better approximations to the algebraic number $r$,
which is what
we wanted to show was possible,
because at each stage the interval containing $r$ is half the size it was before.
 
And this will work if $r$ is what is called a ``simple'' root of its defining equation $\phi(r) = 0$,
because in that case the curve for $\phi(x)$ will in fact cross zero at $x = r$.
But if $r$ is what is called a ``multiple'' root, then the curve may just graze zero, not cross it,
and the Sturm sequence approach is the slickest way to proceed.
 
Now let's stand back and take a look at Turing's proof that there are transcendental reals.
On the one hand, it's philosophical like Cantor's proof; on the other hand, it is some
work to verify in detail that all algebraic reals are computable, although to me that
seems obvious in some sense that I would be hard-pressed to justify/explain.
 
At any rate, now I'd like to take another big step, and show you that there are
uncomputable reals in a very different way from the way that Turing did it, which is very 
much in the spirit of Cantor.  Instead I'd like to use probabilistic ideas, ideas
from what's called measure theory, which was developed by Lebesgue, Borel, and Hausdorff,
among others, and which immediately shows that there are uncomputable reals in a totally
un-Cantorian manner. 

\section*{Reals are uncomputable with probability one!}

I got this idea from reading Courant and Robbins, 
\emph{What is Mathematics?,}
where they give a measure-theoretic
proof that the reals are non-denumerable (more numerous than the integers).
 
Let's look at all the reals in the unit interval between zero and one.
The total length of that interval is of course exactly one.
But it turns out that all of the computable reals in it can be covered with
intervals having total length exactly $\epsilon$, and we can make $\epsilon$ as small as
we want. How can we do that?
 
Well, remember that Turing points out that all the possible computer programs can be
put in a list and numbered one by one, so there's a first program, a second program,
and so forth and so on\ldots\
Some of these programs don't compute computable reals digit by digit; let's just forget
about them and focus on the others.
So there's a first computable real, a second computable real, etc.
And you just take the first computable real and cover it with an interval of size $\epsilon/2$,
and you take  the second computable real and you 
cover it with an interval of size $\epsilon/4$, and you keep going that way, halving the
size of the covering interval each time.
So the total size of all the covering intervals is going to be exactly
\[
   \epsilon/2 + \epsilon/4 + \epsilon/8 + \epsilon/16 + \epsilon/32 + \ldots \; = \; \epsilon
\] 
which can be made as small as you like.
 
And it doesn't matter if some of these covering intervals fall partially outside of the
unit interval, that doesn't change anything.
 
So all the computable reals can be covered this way, using an arbitrarily small part 
$\epsilon$
of
the unit interval, which has length exactly equal to one.
 
So if you close your eyes and pick a real number from the unit interval at random,
in such a way that any one of them is equally likely, the probability is zero that
you get a computable real.  And that's also the case if you get the successive binary
digits of your real number using independent tosses of a fair coin.  It's possible
that you get a computable real, but it's infinitely unlikely.  So with probability one
you get an uncomputable real, and that has also got to be a transcendental number, what
do you think of that!
 
Liouville, Hermite and Lindemann worked so hard to exhibit individual transcendentals,
and now we can do it, almost certainly, by just picking a real number out of a hat!
That's progress for you!
 
So let's suppose that you do that and get a specific uncomputable real that I'm
going to call $R*$.  What if you try to prove what some of its bits are when you
write $R*$ in base-two binary?
 
Well, we've got a problem if we try to do that\ldots

\section*{A FAS cannot determine infinitely many bits of an uncomputable real}

The problem is that if we are using a Hilbert/Turing/Post FAS, as we saw in Chapter II
there has got to be an algorithm for computably enumerating all the theorems, and
so if we can prove what all the bits are, then we can just go ahead and calculate $R*$
bit by bit, which is impossible.  
To do that, you would just go through all the theorems one by one until you find the value
of any particular bit of $R*$.
 
In fact, the FAS has got to fail infinitely many
times to determine a bit of $R*$, otherwise we could just keep a little (finite) table
on the side telling us which bits the FAS misses, and what their values are, and again
we'd be able to compute $R*$, by combining the table with the FAS, which is impossible.
 
So our study of transcendentals now has a new ingredient, which is that we're using
probabilistic methods.  This is how I get a specific uncomputable real $R*$, not the
way that Turing originally did it, which was using Cantor's diagonal method.
 
And now let's really start using probabilistic methods, let's start talking about
algorithmically incompressible, irreducible random real numbers.  What do I mean by that?

\section*{Random reals are uncomputable and have probability one}

Well, a random real is a real number with the property that its bits are irreducible, incompressible
information, as much as is possible.
Technically, the way you guarantee that, is to require that the smallest self-delimiting
binary program that calculates the first $N$ bits of the binary expansion of $R$ is always greater
than $N - c$ bits long, for all $N$, where $c$ is a constant that depends on $R$ but not on $N$.
 
\begin{center}
\begin{tabular}{|c|}
\hline
\\
\textbf{\emph{\large What's a ``Random'' Real Number R?}}
\\ \\
There's a constant $c$ such that 
\\
$H$(the first $N$ bits of $R) > N - c$ for all $N$.
\\ \\
The program-size complexity of the first $N$ bits of $R$ 
\\
can never drop too far below $N$.
\\
\\
\hline
\end{tabular}
\end{center}
 
But I don't want to get into the details. The general idea is that you just require that
the program-size  complexity that we defined in Chapter III of the first $N$ bits of $R$ has
got to be as large as possible, as long as most reals can satisfy the lower bound that
you put on this complexity. That is to say, you demand that the complexity be as high
as possible, as long as the reals that can satisfy this requirement continue to have
probability one, in other words, as long as the probability that a real fails to have
complexity this high is zero.
 
So this is just a simple application of the ideas that we discussed in Chapter III.
 
So of course this random incompressible real won't be computable, for computable reals
only have a finite amount of information.
But what if we try to use a particular FAS to determine particular bits of a particular
random real $R*$?
Now what?

\section*{A FAS can determine only finitely many bits of a random real!}

Well, it turns out that things have gotten worse, much worse than before. Now we can
only determine finitely many bits of that random real $R*$.
 
Why?
 
Well, because if we could determine infinitely many bits, then those bits aren't really
there in the binary expansion of $R*$, we get them for free, you see, by generating all
the theorems in the FAS.   You just see how many bits there are in the program for
generating all the theorems, you see how complex the FAS is.  Then you use the FAS to
determine that many bits of $R*$, and a few more (just a little more than the $c$ in the
definition of the random real $R*$).   And then you fill in the holes up to the last bit
that you got using the FAS. 
(This can
actually be done in a self-delimiting manner, 
because we already know exactly how many bits we need, we know that in advance.)
And so that
gives you a lot of bits from the beginning of the binary expansion of $R*$,  but the program
for doing it that I described is just a little bit too small, its size is substantially
smaller than the number of bits of $R*$ that it got us, which contradicts the definition
of a random real.
 
So random reals are bad news from the point of view of what can be proved.
Most of their bits are unknowable; a given FAS can only determine about as many bits of the
random real as the number of bits needed to generate the set of all its theorems.
In other words, the bits of a random real, any finite set of them,
cannot be compressed into a FAS with a smaller number of bits.
 
So by using random reals we get a much worse incompleteness result than by merely using
uncomputable reals.  You get at most a finite number of bits using any FAS. In fact,
essentially the only way to prove what a bit of a particular random real is using a FAS
is if you put that information directly into the axioms! 
 
The bits of a random real are maximally unknowable!
 
So is there any reason to believe in such a random real?   Plus I need to pick out a specific
one, otherwise this incompleteness result isn't very interesting.  It just says, pick
the real $R*$ out of a hat, and then any particular FAS can prove what at most a finite number of bits
are.  But there is no way to even refer to that specific real number $R*$ within the FAS: it doesn't
have a name!  
 
Well, we'll solve that problem in the next chapter by picking out one random real, the
halting probability $\Omega$.
 
Meanwhile, we've realised that naming individual real numbers can be a problem.
In fact, most of them \textbf{can't even be named}.

\section*{Reals are un-nameable with probability one!}

The proof is just like the one that computable reals have probability zero.
The set of all names for real numbers, if you fix your formal language or FAS, is just a countable
infinity of names: because there's a first name, a second name, etc.
 
So you can cover them all using intervals that get smaller and smaller, and
the total size of all the covering intervals is going to be exactly as before
\[
   \epsilon/2 + \epsilon/4 + \epsilon/8 + \epsilon/16 + \epsilon/32 + \ldots \; = \; \epsilon
\] 
which, as before, can be made as small as you like.
 
So, with probability one, a specific real number chosen at random cannot even be
named uniquely, we can't specify it somehow, constructively or not, we can't define it
or even refer to it!
 
So why should we believe that such an un-nameable real even exists?!
 
\textbf{I claim that this makes incompleteness obvious:
a FAS cannot even name all reals!}
 
Instantaneous proof! 
We did it in three little paragraphs!
Of course, this is a somewhat different kind of incompleteness than the one that G\"odel
originally exhibited.
 
People were very impressed by the technical difficulty of G\"odel's proof.
It had to be that difficult, because it involves constructing an assertion about whole
numbers that can't be proved within a particular, popular FAS (called Peano arithmetic).
But if you change whole numbers to real numbers, and if you talk about what you can
name, rather than about what you can prove, then incompleteness is, as we've just seen, immediate!
 
So why climb
the north face of the Eiger, when you can take a helicopter and have
a picnic lunch on the summit in the noonday sun?
Of course, there are actually plenty of reasons to climb that north face.
I was recently proud to shake hands with someone who's tried to do it several times.
 
But in my opinion, this \textbf{is} the best proof of incompleteness!
As P\'olya says in \emph{How to Solve It}, after you solve a problem, 
if you're a future mathematician, that's just the beginning.
You should look back, reflect on what you did,
and how you did it, and what the alternatives were.
What other possibilities are there?
How general is the method that you used?
What was the key idea?
What else is it good for?
And can you do it without any computation,
or can you see it at a glance?
 
This way, the way we just did it here, yes, you can certainly see incompleteness at a glance!
 
The only problem is, that most people are not going to be too impressed by this particular
kind of incompleteness, because it seems too darn philosophical.
G\"odel made it much more down to earth by talking about the
positive integers, instead of exploiting problematical aspects of the reals.
 
And in spite of what mathematicians may brag, they have always had a slightly
queasy attitude about reals, it's only with whole numbers that they feel really and 
absolutely confident.
 
And in the next chapter I'll solve that, at least for $\Omega$.
I'll take my paradoxical real, $\Omega$, and dress
it up as a diophantine problem, which just talks about whole numbers.
This shows that even though $\Omega$ is a real number, you've got to take it seriously!

\section*{Chapter summary}

\begin{center}
\begin{tabular}{|c|}
\hline
\\
\textbf{\emph{\large Against Real Numbers!}}
\\ \\
\textbf{Prob}\{algebraic reals\} = \textbf{Prob}\{computable reals\} = 
\\
\textbf{Prob}\{nameable reals\} = 0
\\ \\
\textbf{Prob}\{transcendental reals\} = \textbf{Prob}\{uncomputable reals\} = 
\\
\textbf{Prob}\{random reals\} = \textbf{Prob}\{un-nameable reals\} = 1
\\
\\
\hline
\end{tabular}
\end{center}
 
In summary:
\textbf{
Why should I believe in a real number if I can't calculate it, if I can't prove
what its bits are, and if I can't even refer to it? 
And each of these things happens with probability one!
}
 
We started this chapter with $\sqrt{2}$, which the ancient Greeks referred to as ``unutterable'',
and we ended it by showing that with probability one there is no way to 
even name or to specify or define or refer to, no matter how non-constructively,
individual real numbers. We have come full circle, from the unutterable to the un-namable!
There is no escape, these issues will not go away!
 
In the previous chapter we saw \textbf{physical arguments} against real numbers.
In this chapter we've seen that reals are also problematic from a \textbf{mathematical} point
of view, mainly because they contain an \textbf{infinite} amount of information, and
infinity is something we can imagine but rarely touch.
So I view these two chapters as validating the discrete, digital information approach of AIT,
which does not apply comfortably in a physical or mathematical world made up out of real numbers.
And I feel that this gives us the right to go ahead and see what looking at the size of
computer programs can buy us, now that we feel reasonably comfortable with this new digital, discrete
viewpoint, now that we've examined the philosophical underpinnings 
and the tacit assumptions
that allowed us to posit this
new concept.
 
\textbf{In the next chapter I'll finally settle my two debts to you, dear reader, by proving that
Turing's halting problem cannot be solved---my way, not the way that Turing originally did.
And I'll pick out an individual random real, my $\Omega$ number, which as we saw in this
chapter must have the property that any FAS can determine at most finitely many bits of its
base-two binary expansion.
And we'll discuss what the heck it all means,
what it says about how we should do mathematics\ldots}

\chapter*{Chapter VI---Complexity, Randomness \& Incompleteness}
\addcontentsline{toc}{chapter}{Chapter VI---Complexity, Randomness \& Incompleteness}
\markright{Chapter VI---Complexity, Randomness \& Incompleteness}

In Chapter II I showed you Turing's approach to incompleteness.
Now let me show you how I do it\ldots
 
I am very proud of my two incompleteness results in this chapter!
These are the jewels in the AIT crown, the best (or worst) incompleteness results,
the most shocking ones, the most devastating ones, the most enlightening ones,
that I've been able to come up with!
Plus they are a consequence of the digital philosophy viewpoint that goes back to
Leibniz and that I described in Chapter III. That's why these results are so astonishingly
different from the classical incompleteness results of G\"odel (1931) and Turing (1936).

\section*{Irreducible Truths and the Greek Ideal of Reason}

I want to start by telling you about the very dangerous idea of
``logical irreducibility''\ldots
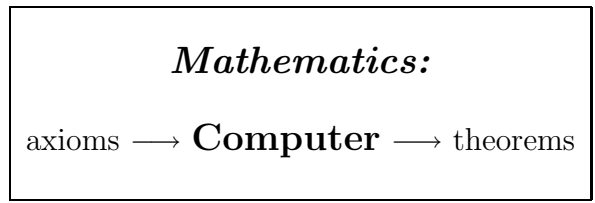
\begin{center}
\begin{tabular}{|c|}
\hline
\\
\textbf{\emph{\large Mathematics:}}
\\ \\
axioms $\longrightarrow$ \textbf{\large Computer} $\longrightarrow$ theorems
\\
\\
\hline
\end{tabular}
\end{center}
We'll see in this chapter that
the traditional notion of what math is about is all wrong: reduce things to axioms,
compression. Nope, sometimes this doesn't work at all. The irreducible mathematical facts
exhibited here in this chapter---the bits of $\Omega$---\textbf{cannot} be derived from 
any principles simpler than they themselves are.
 
So the normal notion of the utility of proof fails for them---proof doesn't help at all
in these cases.
Simple axioms with complicated results is where proof helps.
But here the axioms have to be as complicated as the result. So what's the point
of using reasoning at all?!
 
Put another way: The normal notion of math is to look for structure and law in the world
of math, for a theory.  But theory implies compression, and here there cannot be any---there
is no structure or law at all in this particular corner of the world of math.
 
And since there can be no compression, there can be no understanding of these mathematical facts!
 
In summary\ldots
\begin{center}
\begin{tabular}{|c|}
\hline
\\
\textbf{\emph{\large When is Reasoning Useful?}}
\\ \\
``Axioms = Theorems'' implies reasoning is useless!
\\ \\
``Axioms $\ll$ Theorems'' implies compression \& comprehension!
\\
\\
\hline
\end{tabular}
\end{center}
If the axioms are \textbf{exactly equal} in size to the
body of interesting theorems, then reasoning was absolutely useless.
But if the axioms are \textbf{much smaller} than the body of interesting theorems,
then we have a substantial amount of compression, and so a substantial amount of understanding!
 
Hilbert, taking this tradition to its extreme, 
believed that a single FAS of finite complexity, a finite
number of bits of information, must suffice to generate \textbf{all} of mathematical truth.
He believed in a final theory of everything, at least for the world of pure math.
The rich, infinite, imaginative, open-ended world of all of math, 
all of that compressed into a finite number of bits!
What a magnificent compression that would have been!
What a monument to the power of human reason!

\section*{Coin Tosses, Randomness vs.\ Reason, True for no Reason, Unconnected Facts}

And now, violently opposed to the Greek ideal of pure reason:
Independent tosses of a fair coin, an idea from physics!
 
A ``fair'' coin means that it is as likely for heads to turn up as for tails.
``Independent'' means that the outcome of one coin toss doesn't influence the next outcome.
 
So each outcome of a coin toss is a unique, atomic fact that has no connection with any other
fact: not with any previous outcome, not with any future outcome.
 
And knowing the outcome of the first million coin tosses, if we're dealing with independent
tosses of a fair coin, gives us absolutely no help in predicting the very next outcome.
Similarly, if we could know all the even outcomes (2nd coin toss, 4th toss, 6th toss), that
would be no help in predicting any of the odd outcomes (1st toss, 3rd toss, 5th toss).
 
This idea of an infinite series of
independent tosses of a fair coin may sound like a simple idea, a toy physical
model, but it is a serious challenge, indeed a horrible nightmare,
for any attempt to formulate a rational world view!
Because each outcome is a fact that is true for \textbf{no reason}, that's true only by accident!
 
\begin{center}
\begin{tabular}{|c|}
\hline
\\
\textbf{\emph{\large Rationalist World View}}
\\ \\
In the physical world, everything happens for a reason.
\\ \\
In the world of math, everything is true for a reason.
\\ \\
The universe is comprehensible, logical!
\\ \\
Kurt G\"odel subscribed to this philosophical position.
\\
\\
\hline
\end{tabular}
\end{center}
 
So rationalists like Leibniz and Wolfram have always rejected physical randomness, or
``contingent events'', as Leibniz called them, because they cannot be understood using reason,
they utterly refute the power of reason. 
Leibniz's solution to the problem is to claim that contingent
events are also true for a reason, but in such cases 
there is in fact an infinite series of reasons, an infinite chain of cause and effect, that
while utterly beyond the power of human comprehension, 
is not at all beyond the power of comprehension of the divine mind.  Wolfram's solution to the
problem is to say that all of the \textbf{apparent} randomness 
that we see 
in the world is actually only \textbf{pseudo}-randomness.
It \textbf{looks like} randomness, but it is actually the result of simple laws, in the same
way that the digits of $\pi = 3.1415926\ldots$ \textbf{seem} to be random.
 
Nevertheless, at this point in time quantum mechanics demands real intrinsic randomness
in the physical world, real unpredictability, and chaos theory even shows that a somewhat
milder form of randomness is actually present in classical, deterministic physics, if you
believe in infinite precision real numbers and in the power of acute sensitivity to initial
conditions to quickly amplify random bits in initial conditions into the macroscopic domain\ldots 
 
The physicist Karl Svozil has the following interesting position on these
questions.  Svozil has classical, deterministic leanings 
and sympathizes with Einstein's assertion that
``God doesn't play dice''.  Svozil admits that in its 
\textbf{current} state quantum theory contains randomness.
But he thinks that this is only temporary, and that some new, deeper, hidden-variable theory will
eventually
restore determinacy and law to physics.  On the other hand, Svozil believes that, as he puts it,
$\Omega$ shows that
there is \textbf{real} randomness in the (unreal) mental mindscape fantasy world of pure math!

\section*{A Conversation on Randomness at City College, NY, 1965}

As background information for this story, I should start by telling
you that a century ago Borel proposed a
mathematical definition of a random real number, in fact, infinitely many variant definitions.
He called such reals ``normal'' numbers.
This is the same Borel who in 1927 invented the know-it-all-real that we discussed in Chapter V.
And in 1909 he was able to show that most real numbers must satisfy
\textbf{all} of his variant definitions of normality.
The probability is zero of failing to satisfy any of them.
 
What is Borel's definition of a normal real? Well, it's a real number with the property
that every possible digit occurs with equal limiting frequency: in the long run, exactly 10\%
of the time. That's called ``simple'' 10-normal.  
And 10-normal \emph{tout court} means that for each $k$, the base-ten
expansion of the real has each of the $10^k$ possible blocks of $k$ successive digits
with exactly the same limiting relative frequency $1/10^k$. All the possible blocks
of $k$ successive digits are exactly equally likely to appear in the long run, and this is the case for
\textbf{every} $k$, for $k = 1, 2, 3, \ldots$
Finally there is just plain ``normal'', which is a real that has this property when written
in \textbf{any} base, not just base-ten.
In other words, normal means that it's 2-normal, 3-normal, 4-normal, and so forth and so on.
 
So most reals are normal, with probability one (Borel, 1909).
But what if we want to exhibit a \textbf{specific} normal number?
And there seems to be no reason to doubt that $\pi$ and $e$ are normal, but no one,
to this day, has the slightest idea how to prove this.
 
Okay, that's the background information on Borel normality.
Now let's fast forward to 1965.
 
I'm an undergraduate at City College, just starting my second year.
I was writing and revising my very first paper on randomness.
My definition of randomness was not at all like Borel's.
It was
the much more demanding definition that I've discussed in this book,
that's based on the idea of algorithmic incompressibility, on the idea of looking at the size
of computer programs. 
In fact, it's actually based on a remark that Leibniz made in 1686, though I wasn't aware of
that at the time.
This is my paper that was published in 1966 and 1969 in the \emph{ACM Journal}.
 
The Dean had excused me from attending classes so that I could prepare all this for publication,
and the word quickly got around at City College that I was working on a new definition of randomness.
And it turned out that 
there was
a professor at City College, Richard Stoneham, whom I had never had in any of my courses,
who was very interested in normal numbers.  
 
We met in his office in wonderful old, stone pseudo-gothic
Shepard Hall, 
and I explained to him that I was working on a definition of randomness, one that would
imply Borel normality,
and that it was a definition that most reals
would satisfy, with probability one.
  
He explained that he was interested
in proving that \textbf{specific} numbers like $\pi$ or $e$ were normal. 
He wanted to
show that some well-known
mathematical object already contained randomness, Borel normality, not some new kind of randomness.
I countered that no number like $\pi$ or $e$ could satisfy
my much more demanding definition of randomness, 
because they were computable reals and therefore compressible.
 
He gave me some of his papers. One paper was computational and was 
on the distribution of the different
digits in the decimal expansions of $\pi$ and $e$. 
They seemed to be normal\ldots\
And another paper was theoretical, on digits in rational
numbers, in periodic decimal expansions.  Stoneham was able to show that for some rationals
$m/n$ the digits were sort of equi-distributed,
subject to conditions on $m$ and $n$ that I no longer recall. 
 
That was it, that was the only time that we met.
 
And the years went by, and I came up with the halting probability $\Omega$. And I never
heard of Stoneham again, until I learned in Borwein and Bailey's chapter on normal
numbers in \emph{Mathematics by Experiment} 
that Stoneham had actually managed to do it!  
Amazingly enough, we had \textbf{both} achieved our goals!
 
Bailey and Crandall, during the course of
their own work in 2003 that produced a much stronger result, had discovered that 30
years before, Stoneham, now deceased (that's information from Wolfram's book), had succeeded in
finding what as far as I know was the first ``natural'' example of a normal number.
 
Looking like a replay of Liouville, Stoneham had not been able to show that
$\pi$ or $e$ were normal.  But he had been able to show that the sum of a natural-looking
infinite series was 2-normal, that is, normal for blocks of bits of every possible size in base-2!
 
He and I \textbf{both} found what we were looking for!  
How delightful!
 
And in this chapter I'll tell you how we did it.  But first, as a warm-up
exercise, to get in the right mood, I want to prove that Turing's halting problem
is unsolvable.  We should do that before looking at my halting probability $\Omega$.
 
And to do that, I want to show that you can't prove that a program is elegant,
except finitely often.  Believe it or not, the idea for that proof actually comes
from \textbf{\'Emile Borel}, the same Borel as before.
Although I wasn't aware of this when I originally found the proof on my own\ldots
 
So let me start by telling you Borel's beautiful idea.

\section*{Borel's Undefinability-of-Randomness Paradox}
  
I have a very clear and distinct recollection of reading the idea that I'm about
to explain to you,
in an English translation of a book by Borel on the basic ideas of probability theory.
Unfortunately, 
I have never been able to discover which book it was, 
nor to find again the following discussion by Borel, 
about a paradox regarding any attempt to give a definitive notion of randomness.
 
So it is possible that this is a false memory, perhaps a dream that I once had, which
sometimes seem real to me, or a ``reprocessed'' memory that I have somehow fabricated
over the years.  
 
However, Borel was quite prolific, and was very interested in these issues,
so this is probably somewhere in his oeuvre!
If you find it, please let me know!
 
That caution aside, let me share with you my recollection of Borel's discussion
of a problem that must inevitably arise with \textbf{any} attempt to define the notion
of randomness.
 
Let's say that somehow you can distinguish between whole numbers whose decimal digits
form a particularly random sequence of digits, and those that don't.
 
Now think about the first $N$-digit number that satisfies your definition of randomness.
But this particular number is rather atypical, because it happens to be precisely
the first $N$-digit number that has a particular property!
 
The problem is that random means ``typical, doesn't stand out from the crowd, no distinguishing
features''.  But if you can define randomness, then the property of being random becomes 
just one more feature that you can use to show that certain numbers are atypical and
stand out from the crowd!
 
So in this way you get a hierarchy of notions of randomness, the one you started out with, then
the next one, then one derived from that, and so forth and so on\ldots\  And each of these
is derived using the previous definition of randomness as one more characteristic
just like any other that can be used to classify numbers!
 
Borel's conclusion is that there can be no one definitive definition of randomness.
You can't define an all-inclusive notion of randomness.
Randomness is a slippery concept, there's something paradoxical about it, it's hard to grasp.
It's all a matter of deciding how much we want to demand. You have to decide on a cut-off, 
you have to say ``enough'', let's take \textbf{that} to be random.
 
Let me try to explain this by using some images that are definitely not in my 
possibly false Borel recollection.
 
The moment that you fix in your mind a notion of randomness, that very mental act invalidates
that notion and creates a new more demanding notion of randomness\ldots\
So fixing randomness in your mind is like trying to stare at something without blinking
and without moving your eyes.  If you do that, the scene starts to disappear
in pieces from your visual field.  To see something, you have to keep moving your eyes,
changing your focus of attention\ldots
 
The harder you stare at randomness, the less you see it!  
It's sort of like trying to see faint objects
in a telescope at night, 
which you do by not staring directly at them, but instead looking aside, where the
resolution of the retina is lower and the color sensitivity is lower, but you can see
much fainter objects\ldots
 
This discussion that I shall attribute to Borel in fact contains the germ of the idea of my proof
that I'll now give that you can't prove that a computer program is ``elegant'', that
is to say, that it's the smallest program that produces the output that it does.  More precisely,
a program is elegant if no program smaller than it produces the same output.  You see,
there may be a tie, there may be several different programs that have exactly the same
minimum-possible size and produce the same output.
 
Viewed as a theory, as discussed in Chapter III, an elegant program is the optimal
compression of its output, it's the simplest scientific theory for that output, considered
as experimental data. 
 
So if the output of this computer program is the entire universe, then an elegant program
for it would be the optimal TOE, the optimal Theory of Everything, the one with no redundant elements,
the best TOE, the one, Leibniz would say, that a perfect God would use to produce that particular
universe.

\section*{Why Can't You Prove that a Program is ``Elegant''?}

So let's consider a Hilbert/Turing/Post formal axiomatic system FAS, which (as discussed in
Chapter II) for us
is just a program for generating all the theorems.  And we'll assume that this program
is the smallest possible one that produces that particular set of theorems.  So the size of this
program is precisely the program-size complexity of that theory. There's absolutely no redundancy!
 
Okay, so that's how we
can measure the power of FAS by how many bits of information it contains. As you'll see
now this really works, it really gives us new insight.
 
Building on Borel's paradox discussed in the previous section, consider this 
computer program\ldots
 
\begin{center}
\begin{tabular}{|c|}
\hline
\\
\textbf{\emph{\large Paradoxical Program $P$:}}
\\ \\
The output of $P$ is the same as the output of   
\\
the first provably elegant program $Q$ that you encounter
\\
(as you generate all the theorems of your chosen FAS) 
\\
that is larger than $P$.
\\
\\
\hline
\end{tabular}
\end{center}
 
\textbf{Why are we building on Borel's paradox?  Because the idea is, that if we could prove
that a program is elegant, then that would enable us to find a smaller program that
produces the same output, contradiction!  Borel's paradox, is that if we could define
randomness, then that would enable us to pick out a random number that isn't at all random,
contradiction. So I view these two proofs, Borel's informal paradox, and this actual
theorem, or, more precisely, meta-theorem, as being in the same spirit.}
 
In other words, $P$ generates all the theorems of the FAS until it finds a proof that
a particular program $Q$ is elegant, and what's more the size of $Q$ has got to be larger than
the size of $P$.  If $P$ can find such a $Q$, then it runs $Q$ and produces $Q$'s output as its own output.
 
Contradiction, because
$P$ is too small to produce the same output as $Q$, since $P$ is smaller than $Q$ and
$Q$ is supposed to be elegant (assuming that all theorems proved in the FAS are correct, are true)!
The only way to avoid this contradiction is that $P$ never finds $Q$, because no program $Q$ that's
larger than $P$ is ever shown to be elegant in this FAS.
 
So using this FAS, you can't prove that a program $Q$ is elegant if it's larger than $P$ is.
So in this FAS, you can't prove that more than finitely many specific programs are elegant. 
(So the halting problem is unsolvable, as we'll see in the next section.)
 
And $P$ is just a fixed number of bits larger than the FAS that we're using.
Mostly $P$ contains the instructions for generating all the theorems, that's the variable part,
plus a little
more, a fixed part, to filter them and carry out our proof as above.
 
So the basic idea is that you can't prove that a program is elegant if it's
larger than the size of the program for generating all the theorems in your FAS.
In other words, if the program is larger than the program-size complexity of your FAS,
then you can't prove that that program is elegant!
No way!
 
You see how useful it is to be able to measure the complexity of a FAS, to be able to
measure how many bits of information it contains?
 
\textbf{Tacit assumption}: In this discussion we've assumed the soundness of the FAS, which is to
say, we've assumed that all the theorems that it proves are true.
 
\textbf{Caution}: We don't really have complete irreducibility yet, because the different cases
``$P$ is elegant'', ``$Q$ is elegant'', ``$R$ is elegant'' are \textbf{not} independent of one another.
In fact, you can determine all elegant programs less than $N$ bits in size 
\textbf{from the same} $N$-bit axiom, one that tells you which program with less than $N$
bits takes longest to halt.
Can you see how to do this?
 
But before we solve this problem and achieve total irreducibility with the bits of
the halting probability $\Omega$, let's pause and
deduce a useful corollary.

\section*{Back to Turing's Halting Problem}
 
Here's an immediate consequence, it's a corollary, of the fact that we've just established
that you can't mechanically find more than finitely many elegant programs:
 
There is no algorithm to solve the halting problem, to decide whether or not a given
program ever halts. 
Proof by \emph{reductio ad absurdum}:
Because if there \textbf{were such an algorithm}, we could use it to find all the elegant programs.
You'd do this by checking in turn every program to see if it halts, and then running
the ones that halt to see what they produce, and then keeping only the first program
you find that produces a given output.  If you look at all the programs in size order,
this will give you precisely all the elegant programs (barring ties, which are unimportant
and we can forget about).
 
In fact, we've 
actually just
shown that if 
our halting problem algorithm is
$N$ bits in size,
then there's got to be a program that never halts, 
and that's at most just a few bits larger than $N$ bits in size,
but we can't decide that this program never halts using our $N$-bit halting problem algorithm.
(This assumes that the halting problem algorithm prefers never to give an answer rather than to give
the wrong answer. Such an algorithm can be reinterpreted as a FAS.)
 
So we've just proved Turing's famous 1936 result completely differently from the way that
he originally proved it, 
using the idea of program size, of algorithmic information, of software complexity.
Of all the proofs that I've found of Turing's result, this one is my favorite.
 
And now I must confess something else.  Which is that I don't think that you can
really understand a mathematical result until you find your own proof.
Reading somebody else's proof is not as good as finding your own proof.
In fact, one fine mathematician that I know, Robert Solovay, never let me explain
a proof to him.  He would always insist on just being told the statement of the result,
and then he would think it through on his own!  I was very impressed!
 
This is the most straight-forward, the most direct, the most basic proof of Turing's
result that I've been able to find.
I've tried to get to the heart of the matter, to remove all the clutter, all the
peripheral details that get in the way of understanding.
 
And this also means that this missing piece in our Chapter II discussion of Hilbert's 10th problem
has now been filled in.
 
Why is the halting problem interesting?  
Well, because in Chapter II we showed that if the halting problem is unsolvable,
then Hilbert's 10th problem cannot be solved, and there is no algorithm to
decide whether a diophantine equation has a solution or not.
In fact, Turing's halting problem is equivalent to Hilbert's 10th problem, in the sense
that a solution to either problem would automatically provide/entail a solution for the other one.

\section*{The Halting Probability $\Omega$ and Self-Delimiting Programs}

Now we're really going to get irreducible mathematical facts, mathematical facts that
``are true for no reason'', and which simulate in pure math, as much as is possible, independent
tosses of a fair coin: It's
the bits of the base-two expansion of the halting probability
$\Omega$. The beautiful, illuminating
fact that you can't prove that a program is elegant was just a warm-up
exercise!
 
Instead of looking at \textbf{one} program, like Turing does, and asking whether
or not it halts, let's put \textbf{all possible} programs in a bag, shake it up,
close our eyes, and pick out a program.
What's the probability that this program that we've just chosen at random will eventually halt?
Let's express that probability as an infinite precision binary real between zero and one.
And \emph{voila!,} its bits are our independent mathematical facts.
 
\begin{center}
\begin{tabular}{|c|}
\hline
\\
\textbf{\emph{\large The Halting Probability $\Omega$}}
\\ \\
You run a program chosen by chance on a fixed computer.
\\
Each time the computer requests the next bit of the program,
\\
flip a coin to generate it, using independent tosses of a fair coin.
\\
The computer must decide 
\textbf{by itself} when to stop reading the program.
\\
This forces the program to be self-delimiting binary information.
\\ \\
You sum for each program that halts 
\\
the probability of getting precisely
\\
that program by chance:
\\ \\
{\large $\Omega \; = \; \sum_{\mbox{ program $p$ halts}} 1/2^{\mbox{(size in bits of $p$)}}$}
\\ \\
Each $k$-bit self-delimiting program $p$ that halts
\\
contributes $1/2^k$ to the value of $\Omega$.
\\
\\
\hline
\end{tabular}
\end{center}
 
The self-delimiting program proviso is crucial: Otherwise the halting probability
has to be defined for programs of \textbf{each particular size}, but it cannot be defined
over \textbf{all} programs of \textbf{arbitrary size}.
 
To make $\Omega$ seem more real, let me point out that you can compute it in the limit
from below:
 
\begin{center}
\begin{tabular}{|c|}
\hline
\\
\textbf{\large \emph{Nth Approximation to $\Omega$}}
\\ \\
Run each program up to $N$ bits in size for $N$ seconds.
\\ \\
Then each $k$-bit program you discover that halts   
\\
contributes $1/2^k$ to this approximate value for $\Omega$.
\\ \\
These approximate values get bigger and bigger (slowly!)
\\
and they approach $\Omega$ 
more and more closely, from below.
\\ \\
1st approx.\ $\leq$ 2nd approx.\ $\leq$ 3rd approx.\ $\leq \; \ldots \; \leq \Omega$
\\
\\
\hline
\end{tabular}
\end{center}
 
This process is written in LISP in my book \emph{The Limits of Mathematics.}
The LISP function that gives this approximate value for $\Omega$ as a function of $N$
is about half a page of code using the special LISP dialect that I present in that book.
 
However, this process converges very, very slowly to $\Omega$. In fact, you can never
know how close you are, which makes this a rather weak kind of convergence.  
 
Normally, for approximate values of a real number to be useful, you have to know how
close they are to what they're approximating, you have to know what's called ``the rate
of convergence'', or have what's called ``a computable regulator of convergence''.  But here
we don't have any of that, just a sequence of rational numbers that very slowly creeps
closer and closer to $\Omega$, but without enabling us to ever know precisely how close we are 
to $\Omega$ at a given point in this unending computation.
 
In spite of all of this,
these approximations are extremely useful:
In the next section,
we'll use them in order to show that $\Omega$ is an algorithmically
``random'' or ``irreducible'' real number.
And later in this chapter,
we'll use them to construct diophantine equations for the bits of $\Omega$.

\section*{$\Omega$ as an Oracle for the Halting Problem}

Why are the bits of $\Omega$ irreducible mathematical facts? Well, it's because
we can use the first $N$ bits of $\Omega$ to settle the halting problem for all
programs up to $N$ bits in size.
That's $N$ bits of information, and the first $N$ bits of $\Omega$ are therefore an irredundant
representation of this information.

\begin{center}
\begin{tabular}{|c|}
\hline
\\
\textbf{\large \emph{How much information is there}}
\\
\textbf{\large \emph{in the first N bits of $\Omega$?}}
\\ \\
Given the first $N$ bits of $\Omega$, get better and better
\\
approximations for $\Omega$ as indicated in the previous section,
\\
until the first $N$ bits of the approximate value are correct.
\\ \\
At that point you've seen every program up to $N$ bits long
\\
that ever halts. Output something not included in any of the 
\\
output produced by all these programs that halt.  It cannot have been
\\
produced using any program having less than or equal to $N$ bits.
\\ \\
Therefore the first $N$ bits of $\Omega$ cannot be produced with 
\\
any program having substantially less than $N$ bits, and $\Omega$
\\
satisfies the definition of a ``random'' or ``irreducible'' real
\\
number given in Chapter V:
\\ \\
$H$(the first $N$ bits of $\Omega) \; > \; N - c$
\\
\\
\hline
\end{tabular}
\end{center}
 
This process is written in LISP in my book \emph{The Limits of Mathematics.}
The LISP function that produces something with complexity greater than $N$ bits if it is
given any program that calculates the first $N$ bits of $\Omega$,
is about a page of code using the special LISP dialect that I present in that book.
The size in bits of this one-page LISP function is precisely the value of that
constant $c$ with the property that 
$H$(the first $N$ bits of $\Omega$) is greater than $N - c$ for all $N$.
So, the program-size complexity of the first $N$ bits of $\Omega$ never drops
too far below $N$.
 
Now that we know that $\Omega$ is an algorithmically ``random'' or ``irreducible'' real number,
the argument that a FAS can determine only finitely many bits of such a number given in
Chapter V immediately applies to $\Omega$.  
The basic idea is that if $K$ bits of $\Omega$ could be ``compressed'' into a substantially
less than $K$-bit FAS, then $\Omega$ wouldn't really be irreducible.
In fact, using the argument given in Chapter V,
we can say exactly how many bits of
$\Omega$ a given FAS can determine.
Here's the final result\ldots
 
\begin{center}
\begin{tabular}{|c|}
\hline
\\
\textbf{\large \emph{An FAS can only determine}}
\\
\textbf{\large \emph{as many bits of $\Omega$ as its complexity}}
\\ \\
As we showed in Chapter V, there is (another) constant $c$ such that
\\
a formal axiomatic system FAS with program-size complexity $H$(FAS)
\\
can never determine more than $H(\mbox{FAS}) + c$ bits of the value for $\Omega$.
\\ \\
These are theorems of the form ``The 39th bit of $\Omega$ is 0''
\\
or ``The 64th bit of $\Omega$ is 1''.  
\\ \\
(This assumes that the FAS only enables you to prove such theorems 
\\
if they are true.)
\\
\\
\hline
\end{tabular}
\end{center}
 
This is \textbf{an extremely strong incompleteness result}, it's the very best I can do,
because it says that essentially
the only way to determine bits of $\Omega$ is to put that information directly
into the axioms of our FAS, without using any reasoning at all, only, so to speak,
table look-up to determine these finite sets of bits.
 
In other words, the bits of $\Omega$ are logically irreducible, they cannot be obtained
from axioms simpler than they are.  Finally!  We've found a way to simulate independent
tosses of a fair coin, we've found ``atomic'' mathematical facts, 
an infinite series of math facts that have no connection with each other 
and that are, so to speak, ``true for no reason'' (\textbf{no reason simpler than they are}).
 
Also, this result can be interpreted informally as saying that math is random,
or more precisely, contains randomness, namely the bits of $\Omega$.
What a dramatic conclusion!
But a number of serious caveats are in order!
 
Math isn't random in the sense of being arbitrary, not at all---it
most definitely is not the case that
$2 + 2$ is occasionally equal to 5 instead of 4!
But math does contain \textbf{irreducible information}, of which $\Omega$
is the prime example.
 
To say that $\Omega$ is random may be a little confusing. It's a specific
well-determined real number, and technically it satisfies the definition
of what I've been calling a ``random real''. But math often uses familiar
words in unfamiliar ways.  Perhaps a better way to put it is to say
that $\Omega$ is algorithmically incompressible.  Actually, I much prefer
the word ``irreducible''; I'm coming to use it more and more, although for
historical reasons the word ``random'' is unavoidable.
 
So perhaps it's best to avoid misunderstandings and to say that $\Omega$ is irreducible,
which is true both algorithmically or computationally, and \textbf{logically}, by means
of proofs.
And \emph{which happens to imply} that $\Omega$ has many of the characteristics of the typical
outcome of a random process, in the physical sense of an unpredictable process
amenable to statistical treatment.
 
For example, as we'll discuss in the next section, in base two, in $\Omega$'s infinite
binary expansion, each of the $2^k$ possible $k$-bit blocks will appear with
limiting relative frequency exactly $1/2^k$, and this is provably the case
for the specific real number $\Omega$, even though it's only true with probability
one, but not with certainty, for the outcome of an infinite series of independent
tosses of a fair coin.  So perhaps, in retrospect, the choice of the word ``random'' wasn't
so bad after all!
 
Also, a random real may be meaningless, or it may be extremely meaningful; my
theory cannot distinguish between these two possibilities, it cannot
say anything about that.  If the real was produced using an independent
toss of a fair coin for each bit, it'll be irreducible and it'll be meaningless.
On the other hand, $\Omega$ is a random real with lots of meaning, since it contains
a lot of information about the halting problem, and this information is stored in $\Omega$
in an irreducible fashion, with no redundancy.
You see, once you compress out all the redundancy from anything meaningful, the result
necessarily \textbf{looks} meaningless, even though it isn't, not at all, it's just jam-packed
with meaning!

\section*{Borel Normal Numbers Again}

In fact, it's not difficult to see that $\Omega$ is normal, that is to say, $b$-normal
for any base $b$, not just 2-normal.  And the lovely thing about this, is that the
definition of $\Omega$ as the halting probability doesn't seem to have anything to do
with normality.  $\Omega$ wasn't constructed especially with normality in mind; that just dropped out,
for free, so to speak!
 
So for any base $b$ and any fixed number $k$ of base-$b$ ``digits'', the limiting relative
frequency of each of the $b^k$ possible $k$-digit sequences in $\Omega$'s $b$-ary
expansion will be exactly $1/b^k$.  In the limit, they are all equally likely
to appear\ldots
 
How can you prove that this has got to be the case? Well, if this were \textbf{not} the
case, then $\Omega$'s bits would be highly compressible, by a fixed multiplicative factor
depending on just how unequal the relative frequencies happen to be.  In other words, the bits of
$\Omega$ could be compressed by a fixed percentage, which is a lot of compression\ldots\
These are ideas that go back to a famous paper published by Claude Shannon in 
\emph{The Bell System Technical Journal} in
the 1940's, although
the context in which we are working is rather different from his.
Anyway that's where I got the idea, by reading Shannon.
Shannon and I both worked for industrial labs: in his case, it was the phone
company's, in my case, IBM's.
 
So that's how I succeeded in finding a normal number!  It was because I wasn't interested
in normality per se, but in deeper, philosophical issues, and normality just dropped out
as an application of these ideas.
Sort of similar to Turing's proof that transcendentals exist, because all algebraic numbers
have to be computable\ldots
 
And how did the City College professor Richard Stoneham succeed in \textbf{his} quest for
randomness?  Here it is, here is Stoneham's provably 2-normal number:
\[
   1/(3 \times 2^3) + 1/(3^2 \times 2^{3^2})
                    + 1/(3^3 \times 2^{3^3})
             + \ldots 1/(3^k \times 2^{3^k}) + \ldots
\]
This is normal base 2 for blocks of all size.
 
And what is the more general result obtained by
David Bailey and Richard Crandall in 2003? 
\[
   1/(c \times b^c) + 1/(c^2 \times b^{c^2})
                    + 1/(c^3 \times b^{c^3})
             + \ldots 1/(c^k \times b^{c^k}) + \ldots
\]
This is $b$-normal as long as $b$ is greater than 1 and $b$ and $c$ have no common factors.
For example, $b = 5$ and $c = 7$ will do; that gives you a 5-normal number.
For the details, see the chapter on normal numbers, Chapter 4, in Borwein and Bailey,
\emph{Mathematics by Experiment.}
And there's lots of other interesting stuff in that book, for example, amazing new ways to
calculate $\pi$---which actually happens to be connected with these normality proofs!
(See their Chapter 3).

\section*{Getting Bits of $\Omega$ using Diophantine Equations}

In Chapter V, I expressed \textbf{mathematical} skepticism about real numbers.
And in Chapter IV, I expressed \textbf{physical} skepticism about real numbers.
So why should we take the real number $\Omega$ seriously?
Well, it's not just \textbf{any} real number, you can get it from a diophantine equation!
In fact, you can do this in two rather different ways.  
 
One approach, that I discovered in 1987,
makes the number of solutions of an equation jump from finite to infinite in a way that
mimics the bits of $\Omega$. The other approach, 
discovered by Toby Ord and Tien D. Kieu in Australia 2003,
makes the number of solutions of the equation jump from even to odd in a way that mimics
the bits of $\Omega$. So take your choice; it can be done either way. There's a Northern
and there's a Southern Hemisphere approach to this problem---and no doubt many other 
interesting ways to do it!
 
Remember Kronecker's credo that ``God created the integers; all else is
the work of man''?
If you prefer, $\Omega$ isn't a real number at all, it's a fact about certain diophantine
equations; it has to do only with whole numbers, with positive integers!
So you cannot shrug away the fact that the bits of the halting probability $\Omega$ are
irreducible mathematical truths, for this can be reinterpreted 
as a statement about diophantine equations.
 
\begin{center}
\begin{tabular}{|c|}
\hline
\\
\textbf{\large \emph{Chaitin (1987):}}
\\
\textbf{\large \emph{Exponential Diophantine Equation \#1}}
\\ \\
In this equation $n$ is a \textbf{parameter}, 
\\
and $k, x, y, z, \ldots$ are the \textbf{unknowns}:
\\ \\
$L(n, k, x, y, z, \ldots) \; = \; R(n, k, x, y, z, \ldots)$.
\\ \\
It has \textbf{infinitely many} positive-integer solutions
\\
if the $n$th bit of $\Omega$ is a 1.
\\ \\
It has \textbf{only finitely many} positive-integer solutions
\\
if the $n$th bit of $\Omega$ is a 0.
\\
\\
\hline
\end{tabular}
\end{center}
 
\begin{center}
\begin{tabular}{|c|}
\hline
\\
\textbf{\large \emph{Ord, Kieu (2003):}}
\\
\textbf{\large \emph{Exponential Diophantine Equation \#2}}
\\ \\
In this equation $n$ is a \textbf{parameter}, 
\\
and $k, x, y, z, \ldots$ are the \textbf{unknowns}:
\\ \\
$L(n, k, x, y, z, \ldots) \; = \; R(n, k, x, y, z, \ldots)$.
\\ \\
For any given value of the parameter $n$,
\\
it only has finitely many positive-integer solutions. 
\\ \\
\textbf{For each particular value of $n$:}
\\ \\
the number of solutions of this equation will be \textbf{odd}
\\
if the $n$th bit of $\Omega$ is a 1, and 
\\ \\
the number of solutions of this equation will be \textbf{even} 
\\
if the $n$th bit of $\Omega$ is a 0.
\\
\\
\hline
\end{tabular}
\end{center}
 
How do you construct these two diophantine equations?
Well, the details get a little messy. The boxes below give you
the general idea; they summarize what needs to be done.
As you'll see, that computable sequence of approximate values of $\Omega$
that we discussed before plays a key role.  It's also important, particularly
for Ord and Kieu (2003),
to recall that these approximate values are a \textbf{non-decreasing} sequence
of rational numbers that get closer and closer to $\Omega$ but that always
remain \textbf{less than} the value of $\Omega$.
 
{\footnotesize
\begin{center}
\begin{tabular}{|c|}
\hline
\\
\textbf{\large \emph{Chaitin (1987):}}
\\
\textbf{\large \emph{Exponential Diophantine Equation \#1}}
\\ \\
Program$(n,k)$ calculates the $k$th approximation to $\Omega$, 
\\
in the manner explained in a previous section.
\\
Then Program$(n,k)$ looks at the $n$th bit of this approximate value for $\Omega$.
\\
If this bit is a 1, then Program$(n,k)$ immediately halts; otherwise it loops forever.
\\
So Program$(n,k)$ halts if and only if 
\\
(the $n$th bit in the $k$th approximation to $\Omega$) is a 1.
\\ \\
As $k$ gets larger and larger, 
the $n$th bit of the $k$th approximation to $\Omega$
\\
will eventually settle down to
the correct value.
Therefore for all sufficiently large $k$:
\\ \\
Program$(n,k)$ will halt if the $n$th bit of $\Omega$ is a 1,
\\ \\
and Program$(n,k)$ will fail to halt if the $n$th bit of $\Omega$ is a 0.
\\ \\
Using all the work on Hilbert's 10th problem that we explained in Chapter II, 
\\
we immediately get an exponential diophantine equation
\\ \\
$L(n, k, x, y, z, \ldots) \; = \; R(n, k, x, y, z, \ldots)$
\\ \\
that has \textbf{exactly one} positive-integer solution
if Program$(n,k)$ eventually halts,  
\\ \\
and that has \textbf{no} positive-integer solution if Program$(n,k)$ never halts.
\\ \\
\textbf{Therefore}, fixing $n$ and considering $k$ to be an unknown, this exact same equation
\\ \\
$L(n, k, x, y, z, \ldots) \; = \; R(n, k, x, y, z, \ldots)$
\\ \\
has \textbf{infinitely many solutions}
if the $n$th bit of $\Omega$ is a 1,  
\\ \\
and it has \textbf{only finitely many} solutions
if the $n$th bit of $\Omega$ is a 0.
\\
\\
\hline
\end{tabular}
\end{center}
}

{\footnotesize
\begin{center}
\begin{tabular}{|c|}
\hline
\\
\textbf{\large \emph{Ord, Kieu (2003):}}
\\
\textbf{\large \emph{Exponential Diophantine Equation \#2}}
\\ \\
Program$(n,k)$ halts if and only if $k > 0$ and
\\
$2^n \times (j$th approximation to $\Omega) \; > \; k$
\\
for some $j = 1, 2, 3, \ldots$
\\ \\
So Program$(n,k)$ halts if and only if 
$2^n \times \Omega \; > \; k \; > \; 0$.
\\ \\
Using all the work on Hilbert's 10th problem that we explained in Chapter II,
\\
we immediately get an exponential diophantine equation
\\ \\
$L(n, k, x, y, z, \ldots) \; = \; R(n, k, x, y, z, \ldots)$
\\ \\
that has \textbf{exactly one} positive-integer solution
if Program$(n,k)$ eventually halts,
\\ \\
and that has \textbf{no} positive-integer solution if Program$(n,k)$ never halts.
\\ \\
Let's fix $n$ and ask for which $k$ does this equation have a solution.
\\ \\
\textbf{Answer:}
$L(n,k) = R(n,k)$ is solvable precisely for 
\\
$k = 1, 2, 3, \ldots$, 
up to the integer part of $2^n \times \Omega$.
\\ \\
Therefore, $L(n) = R(n)$ has exactly integer part of $2^n \times \Omega$ solutions,
\\
which is the integer you get by shifting the binary expansion for $\Omega$ left $n$ bits.
\\
And
the right-most bit of the integer part of $2^n \times \Omega$ will be the $n$th bit of $\Omega$.
\\ \\
\textbf{Therefore}, fixing $n$ and considering $k$ to be an unknown, this exact same equation
\\ \\
$L(n, k, x, y, z, \ldots) \; = \; R(n, k, x, y, z, \ldots)$
\\ \\
has \textbf{an odd number} of solutions
if the $n$th bit of $\Omega$ is a 1,
\\ \\
and it has \textbf{an even number} of solutions
if the $n$th bit of $\Omega$ is a 0.
\\
\\
\hline
\end{tabular}
\end{center}
}
\newpage

\section*{Why is the Particular Random Real $\Omega$ Interesting?}

This is a good place to discuss a significant issue.  In the previous
chapter we pointed out that with probability one a real number is
algorithmically irreducible.  
Algorithmically compressible reals have probability zero.
So why is the \textbf{particular} random real 
$\Omega$ of any interest?! There are certainly plenty of them!
 
Well, for a number of reasons.
 
First of all, $\Omega$ links us to Turing's famous result; the halting 
problem is unsolvable and the halting probability is random!
Algorithmic unsolvability in one case, algorithmic randomness or incompressibility in the other.
Also, the first $N$ bits of $\Omega$ give us a lot of information about
particular, individual cases of the halting problem.
 
But the main reason that $\Omega$ is interesting is this:
We have reached into the infinitely dark blackness of random reals
and picked out \textbf{a single} random real!  I wouldn't say that
we can touch it, but we can certainly point straight at it.
And it is important to make this obscurity tangible by exhibiting a
specific example.
After all, why should we believe that most things have a certain property
if we can't exhibit a specific thing that has this property?
 
(Please note, however, that 
in the case of the un-nameable reals, 
which also have probability one, we're never ever going to be able to
pick out an individual un-namable real!)
 
Here's another way to put it:
$\Omega$ is violently, maximally uncomputable, but it \textbf{almost} looks
computable.  It's just across the border between what we can deal with
and things that transcend our abilities as mathematicians.  So it serves
to establish a sharp boundary, it draws a thin line in the sand that we dare not
cross, that we \textbf{can not} cross!
 
And that's also connected with the fact that we can compute better and better
lower bounds on $\Omega$, we just can't ever know how close we've gotten.
 
In other words, for an incompleteness result to be really shocking, the
situation has got to be that just as we were about to reach out and touch something,
we get our fingers slapped. We can never ever have it, even though it's attractively
laid out on the dining room table, next to a lot of other inviting food!
That's much more frustrating and much more interesting than being told
that we can't have something or do something that never seemed concrete enough
or down to earth enough that this ever looked like a legitimate possibility in the first place!

\section*{What is the Moral of the Story?}

So the world of mathematical truth has infinite complexity, even though any given FAS
only has finite complexity.  
In fact, even the world of diophantine problems has infinite complexity, no finite FAS will do.
 
I therefore believe that we cannot stick with a single FAS, as Hilbert wanted, we've
got to keep adding new axioms, new rules of inference, or some other kind of new mathematical
information to the foundations of our theory.  And where can we get new stuff that cannot
be deduced from what we already know?  Well, I'm not sure, but I think that it may come
from the same place that physicists get their new equations: based on inspiration, imagination and
on---in the case of math, computer, not laboratory---experiments.
 
So this is a ``quasi-empirical'' view of how to do mathematics, which is a term coined
by Lakatos in an article in Thomas Tymoczko's interesting collection \emph{New Directions in the
Philosophy of Mathematics.}  And this is closely connected 
with the idea of
so-called ``experimental mathematics'', 
which uses computational evidence rather than conventional proof to ``establish'' new truths.
This research methodology, whose benefits are 
argued for in a two-volume
work by Borwein, Bailey and Girgensohn,
may not only sometimes be \textbf{extremely convenient}, as they argue, but in fact 
it may sometimes even be \textbf{absolutely necessary} 
in order for mathematics to be able to progress
in spite of the incompleteness phenomenon\ldots
 
Okay, so that's \textbf{my} approach to incompleteness, and it's rather
different from G\"odel's and Turing's. The main idea is to measure
the complexity or the information content of a FAS by the size in bits of the smallest program
for generating all the theorems. Once you do that, everything else
follows. From that initial insight, it's more or less straight-forward, 
the developments are more or less systematic.
 
Yes, but as P\'olya pointedly asks, ``Can you 
see it all at a glance?''  Yes, in fact I think that we can:
 
\textbf{No mechanical process (rules of the game) can be really creative,
because in a sense anything that ever comes out was already
contained in your starting point. Does this mean that physical randomness,
coin tossing, something non-mechanical, is the only possible source of creativity?!
At least from this (enormously over-simplified!)\ point of view, it is.}
    
\section*{Against Egotism}
     
The history of ideas offers many surprises.
In this book we've seen that:
\begin{itemize}
\item
Digital philosophy can be traced back to Leibniz,
and digital physics can be traced back to Zeno. 
\item
My definition of randomness via complexity goes back to Leibniz.
\item
My $\Omega$ number can be traced back to Borel's know-it-all real number.
Borel's number is also an example of what Turing would later call an uncomputable real.
\item
The main idea of my proof that you can't prove that a program
is elegant---in fact it's the basic idea of \textbf{all} my incompleteness results---can 
be traced back to Borel's undefinability-of-randomness paradox.
\end{itemize}
I think that these observations should serve as an antidote to the excessive
egotism, competitiveness and the foolish fights over priority that poison science.
No scientific idea has only one name on it; they are the joint production
of the best minds in the human race, building on each other's insights
over the course of history.
 
And the fact that these ideas can be traced that far back, 
that the threads are that long, doesn't weaken them in any way. 
On the contrary, it gives them even greater significance.
 
As my friend Jacob T. Schwartz once told me, the medieval cathedrals
were the work of many hands, anonymous hands, and took lifetimes to build.
And Schwartz delighted in quoting a celebrated doctor from that period, who said about a patient,
``I treated him and God cured him!''  I think that this is also the right attitude to have
in science and mathematics.

\chapter*{Chapter VII---Conclusion}
\addcontentsline{toc}{chapter}{Chapter VII---Conclusion}
\markright{Chapter VII---Conclusion}

As you have no doubt noticed, this is really a book on philosophy, not just a
math book.
And as Leibniz says in the quote at the beginning of the book, 
math and philosophy are inseparable.
 
Among the great philosophers, only Pythagoras and Leibniz were great 
mathematicians.  Indeed, Pythagoras, though lost in the mists of time,
is credited with inventing philosophy and mathematics and also inventing
the \textbf{words} ``philosophy'' and ``mathematics''.
As for Leibniz, he is what the Germans call a ``universalgenie'', a universal
genius, interested in everything; and if he was interested in something, 
he always came up with an important new idea, 
he always made
a good suggestion.
 
Since Leibniz, perhaps only Poincar\'e was a little like this.
He was enough of a philosopher to come up with a version of relativity theory before Einstein did,
and his popular books of essays were full of philosophical observations and are
still in print.
 
And Leibniz used the word ``transcendental'' as in transcendental curves, numbers, and methods
etc.\ deliberately, thinking of God's transcendence of all things human, 
which also inspired Cantor to develop his theory of infinite magnitudes.
After all, math deals with the world of ideas, which transcends the real world.
And for ``God'' you can understand the laws of the universe, as Einstein did, or the entire world,
as Spinoza did, that doesn't change the message.
 
The topic of monotheism and polytheism is also germane here. I love
the complexity and sophistication of South Indian vegetarian cuisine, and my
home is decorated with Indian sculpture and fabrics---and much more.
I admired
Peter Brook's \emph{Mahabharata,} which clearly brings out the immense philosophical depth
of this epic.
 
But my intellectual personality is resolutely monotheistic.
Why do I say that my personality is monotheistic? I mean it only in this rather abstract sense:
that I am always searching for simple, unifying ideas, rather than glorying intellectually
in ``polytheistic'' subjects like biology, in which there 
is a rich tapestry of extremely complicated facts
that resists being reduced to a few simple ideas. 
 
And let's look at Hilbert's FAS's. They failed miserably. Nevertheless formalism has been
a brilliant success this past century, but not in math, not in philosophy, but as computer
technology, as software, as programming languages.
It works for machines, but not for us!\footnote
{Mathematicians this century have nevertheless made the fatal mistake of 
gradually eliminating all words
from their papers. That's another story, which I prefer not to discuss here, the Bourbakisation
of 20th century mathematics.}
 
Just look at Bertrand Russell's three-volume magnum opus (with Whitehead) 
\emph{Principia Mathematica,} which I keep in my office at IBM.  
One entire formula-filled volume to prove that $1 + 1 = 2$!
And by modern standards, the formal axiomatic system used there is inadequate; it's not formal
enough!  
No wonder that Poincar\'e decried this as a kind of madness!
 
But I think that it was an interesting philosophical/intellectual exercise to refute this madness
as convincingly as possible; somehow I ended up spending my life on that!
 
Of course, formalism fits the 20th century zeitgeist so well: Everything is meaningless,
technical papers should never discuss ideas, only present
the facts! A rule that I've done my best to ignore!
As Vladimir Tasi\'c says in his book \emph{Mathematics and the Roots of Postmodern Thought,}
a great deal of 20th century philosophy seems enamored with formalism and then senses
that G\"odel has pulled the rug out from under it, and therefore truth is relative.
--- Or, as he puts it, it \textbf{could} have happened this way. ---
Tasi\'c presents 20th century thought as, in effect, a dialogue between Hilbert and Poincar\'e\ldots\
But truth is relative is \textbf{not} the correct conclusion. The correct conclusion is that Hilbert was
wrong and that Poincar\'e was right: intuition cannot be eliminated from mathematics,
or from human thought in general. And it's not all the same, all intuitions are not equally valid.
Truth is \textbf{not} reinvented by each culture or each generation, 
it's not merely a question of fashion.
 
Let me repeat: formal axiomatic systems are a failure!  Theorem proving algorithms \textbf{do not
work}.  One can publish papers about them, but they only prove trivial 
theorems.  
And in the case-histories in this book,
we've seen that
the essence of math resides in its creativity, in imagining new concepts, in changing viewpoints,
not in mindlessly and mechanically grinding away deducing all the 
possible consequences of a fixed set of rules and ideas.
 
Similarly, proving correctness of software using formal methods is
hopeless.  Debugging is done experimentally, by trial and error.
And cautious managers insist on running a new system in parallel 
with the old one
until they believe that the new system works.  
 
Except---what we did in a marvelous project at IBM---for constantly eating
your own cooking. We were constantly running the latest version of our software,
constantly compiling our optimizing compiler through itself.
By constantly using it ourselves, we got instantaneous feedback on the performance
and on the design,
which was constantly evolving.
In my opinion that's the only way to develop a large piece of software:
a totalitarian top-down approach cannot work.
You have to design it as you go, not at the very beginning,
before you write a single line of code.
 
And I feel that my experience in the real world debugging software, yields 
a valuable philosophical lesson:
Experimentation is the only way to ``prove'' that software is correct.
Traditional mathematical proofs are only possible in toy worlds, not in the real world.
The real world is too complicated.
A physicist would say it like this: Pure math can only deal with the hydrogen atom.
One proton, one electron, that's it!
The quasi-empirical view of math may be controversial in the math community, but it is old hat
in the software business.
Programmers already have a quasi-empirical attitude to proof.
Even though software is pure mind-stuff, not physical,
programmers behave like physicists, not mathematicians, when it comes to debugging.
 
On the other hand, a way to minimize the debugging problem is to try at all costs to keep software
intellectually manageable, as illustrated by my discussion of LISP.
In our IBM project, I did this by re-writing all my code
from scratch each time that I advanced in my understanding of the problem.
I refused to use \emph{ad hoc} code and tried instead to base everything 
as much as possible on clean, systematic mathematical algorithms.
My goal was crystalline clarity.
In my opinion,
what counts most in a design is conceptual integrity, being faithful to an idea, not confusing
the issue!
 
The previous paragraph is very pro-math. However, my design evolved,
as I said, based on computer experiments,
and the experimental method is used by physicists, not by mathematicians.
It's too bad that mathematicians feel that way about experimentation!
 
Remember, I think that
math is not that different from physics, for we must be willing to add new axioms:
\begin{center}
\begin{tabular}{|c|}
\hline
\\
   \textbf{\large Physics}: laws $\rightarrow$ {\large Computer} $\rightarrow$ universe
\\ \\
   \textbf{\large Math}: axioms $\rightarrow$  {\large Computer} $\rightarrow$ theorems
\\
\\
\hline
\end{tabular}
\end{center}
To get more out, put more in!
 
Another lesson I'd like you to take from this book is that everything is connected,
all important ideas are---and that fundamental
questions go back millennia and are \textbf{never} resolved.
For example, the tension between the continuous and the discrete, or 
the tension between the world of
ideas (math!)\ and the real world (physics, biology).
You can find all this discussed in ancient Greece. And I suspect 
we could even trace it back to ancient Sumer,
if more remained of
Sumerian math than the scrap paper jottings on the clay tablets that are all we have,
jottings that give hints of surprisingly sophisticated methods and of a love for calculation that
seems to far outstrip any possible practical application.\footnote
{Did Sumer inherit 
its mathematics from an \textbf{even older}
civilization---one more advanced than the ancient Greeks---that was destroyed by
the glaciers, or when the glaciers suddenly melted, or by some other natural catastrophe? 
There is no way for such sophisticated computational techniques to appear
out of nowhere, without antecedents.}

\section*{On Creativity}

Having emphasized and re-emphasized the importance of creativity,  it would be nice if I had
a theory about it.  Nope, but \textbf{I do have} some experience being
creative. So let me try to share that with you.
 
The message of G\"odel incompleteness, 
as I've said again and again in this book, is that a static fixed FAS cannot
work.  You have to add new information, new axioms, new concepts.  Math is constantly evolving.
The problem with current metamath is that it deals only with---it refutes---static FAS's.
So where do new math ideas come from?  Can we have a theory about that?
A dynamic rather than static view of math, a dynamic rather than a static metamath, a kind of
dynamic FAS perhaps?
 
Since I don't have that theory, I think an anecdotal approach might be best.  This book
is full of amazing case studies of new, unexpected math ideas that reduced the complicated
to the obvious. And I've come up with a few of these ideas myself.  How does it feel to do that?
 
Well, you can't find them if you don't look for them, if you don't \textbf{really believe} in them.
 
Is there some way to train for it, like a sport?!
No, I don't think so!  You have to be seized by a demon, and our society doesn't want
too many people to be like that!
 
Let me describe what it feels like right now while I'm writing this book.
 
First of all, the ideas that I'm discussing seem very concrete, real and tangible to me.
Sometimes they even feel more real than the people around me.  They certainly feel
more real than newspapers, shopping malls and TV programs---those 
always give me a tremendous feeling of \textbf{unreality}!  
In fact, I only really feel alive when I'm working on a new idea, when I'm making love
to a woman (which is also working on a new idea, 
the child we might conceive), or when I'm going up a mountain!
It's intense, very intense.
 
When I'm working on a new idea I push everything else away. I stop swimming in the morning,
I don't pay the bills, I cancel my doctor appointments.
As I said, everything else becomes unreal!
And I don't have to force myself to do this.
 
On the contrary,
it's pure sensuality, pure pleasure.
I put beautiful new ideas in the same category with beautiful women and beautiful art.
To me it's like an amazing ethnic cuisine I've never tasted before.
 
I'm not depriving myself of anything, I'm not an ascetic.
I don't look like an ascetic, do I?
 
And you can't force yourself to do it, anymore than a man can force himself to make love
to a woman he doesn't want.
 
The good moments are very, very good!  Sometimes when I'm writing this I don't know
where the ideas come from. I think that it can't be me, that I'm just a channel
for ideas that want to be expressed.
--- But  
I \textbf{have} been concentrating on these questions for a long time. ---
I feel inspired, energized by the ideas.
People may think that something's wrong with me, but I'm okay, I'm more than okay.
It's pure enthusiasm! That's ``God within'' in Greek.
Intellectual elation---like summiting on a high peak!
 
And I'm a great believer in the subconscious, in sleeping on it, in going to bed at 3am
or 5am
after working all night, and then getting up the next morning full of new ideas, ideas that come
to you in waves while you're taking a bath, or having coffee.
Or swimming laps.
So mornings are very important to me, and I prefer to spend them at home.  
Routine typing and e-mail, I do in my office, not at home.
And when I get too tired to stay in the office, then I print out the
final version of the chapter I'm working on, bring it home---where there is \textbf{no}
computer---and lie in bed for hours reading it, thinking about it, making corrections,
adding stuff.
 
Sometimes the best time is lying in bed in the dark with my eyes closed,
in a half dreamy, half awake state that seems
to make it easier for new ideas, or new combinations of ideas, to emerge.
I think of the subconscious as a chemical soup that's constantly making new combinations,
and interesting combinations of ideas stick together, and eventually percolate up into 
full consciousness. --- That's not too different from a biological population in which
individuals fall in love and combine to produce new
individuals. --- My guess is that all this activity 
takes place at a molecular level---like DNA
and information storage in the immune system---not at the cellular level. That's why the
brain is so powerful, because that's where the real information processing is, at a
molecular level. The cellular level, that's just the front end\ldots
 
Yes, I believe in ideas, in the power of imagination and new ideas.  And I don't believe in
money or in majority views or the consensus.  Even if all you are interested in is money, I 
think that new ideas are vital in the long run, which is why a commercial enterprise like IBM
has a Research Division and has supported my work for so long.
Thank you, IBM!
 
And I think that this also applies to human society.  
 
I think that the current zeitgeist is  very dangerous, because people are really desperate
for their lives to be meaningful. They need to be creative, they need to help other people,
they need to be part of a community, they need to be adventurous explorers, all things that
tribal life provided.
That's why much of the art in my home is from so-called ``primitive'' cultures like Africa.
 
So you're not going to be surprised to hear that I think 
that we desperately need new ideas about how human society should be organized,
about what it's all for and how to live.
 
When I was a child, I read a science fiction story by Chad Oliver.
He was an anthropologist (from Greek ``anthropos'' = human being), 
not someone interested in hardware, or in spaceships,
but in what was going on inside, in the human soul.
 
The story was called \emph{Election Day}, and it was about a typical American family consisting of
a working father, a stay-at-home mother, and two children, a boy and a girl.---That was the 1950's,
remember!---Everything seemed normal, 
until you learned that the father, who was running for office,
actually wasn't himself a candidate, no, he had created a social system.  And the election wasn't
electing people, it was voting on the choice of the next social system.  And this wasn't
really 1950's America, they had voted to live that way for a while, and then they were going
to change and try a completely different system, right away, soon after election day!
 
What a remarkable idea!
 
And that's what I think we desperately need: Remarkable new ideas!
Paradigm shifts!
Meaning, intuition, creativity!
We need to reinvent ourselves! 
 
Thank you for reading this book and for taking this journey with me!

\chapter*{Reading a Note in the Journal \emph{Nature} I Learn}
\addcontentsline{toc}{chapter}{Poem on Omega by Robert Chute}
\markright{Poem on Omega by Robert Chute}

\emph{
$\;\;\;\;$
by Robert M. Chute 
}

\begin{verse}
Omega, Omega-like, and Computably Enumerable\\
$\;\;\;\;\;\;\;\;\;$
Random Real Numbers all may be\\
$\;\;\;\;\;\;\;\;\;\;\;\;\;\;\;\;\;\;\;$ 
a single class.

Should I be concerned? Is this the sign\\
$\;\;\;\;\;\;\;\;\;$
of a fatal fault line\\
$\;\;\;\;\;\;\;\;\;\;\;\;\;\;\;\;\;\;\;$
in the logic of our world?

You shouldn't worry, I'm told,\\
$\;\;\;\;\;\;\;\;\;$
about such things, but how to be indifferent\\
$\;\;\;\;\;\;\;\;\;\;\;\;\;\;\;\;\;\;\;$
if you don't understand?
 
Standing here, we do not sense\\
$\;\;\;\;\;\;\;\;\;$
any tectonic rearrangement\\
$\;\;\;\;\;\;\;\;\;\;\;\;\;\;\;\;\;\;\;$
beneath our feet.
 
Despite such reassurance I feel\\
$\;\;\;\;\;\;\;\;\;$
increasing dis-ease. Randomness is rising\\
$\;\;\;\;\;\;\;\;\;\;\;\;\;\;\;\;\;\;\;$ 
around my knees.
  
Remember how we felt when we learned\\
$\;\;\;\;\;\;\;\;\;$
that one infinity\\
$\;\;\;\;\;\;\;\;\;\;\;\;\;\;\;\;\;\;\;$
could contain another?
 
Somewhere between the unreachable\\
$\;\;\;\;\;\;\;\;\;$
and the invisible\\
$\;\;\;\;\;\;\;\;\;\;\;\;\;\;\;\;\;\;\;$
I had hoped for an answer.\footnote
{From \emph{Beloit Poetry Journal,} Spring 2000 
(Vol.\ 50, No.\ 3, p.\ 8).
The note in \emph{Nature} in question is C. S. Calude, G. J. Chaitin,
``Randomness everywhere,'' \emph{Nature,} 22 July 1999 (Vol.\ 400, pp.\ 319--320).}
\end{verse}

\chapter*{Math Poem}
\addcontentsline{toc}{chapter}{Math Poem by Marion Cohen}
\markright{Math Poem by Marion Cohen}

\begin{verse}
  Someone wrote a book called \emph{The Joy of Math.}
  
  Maybe I'll write a book called \emph{The Pathos of Math.}

  For through the night I wander

  between intuition and calculation
  
  between examples and counter-examples
  
  between the problem itself and what it has led to.

  I find special cases with no determining vertices.

  I find special cases with only determining vertices.

  I weave in and out.

  I rock to and fro.

  I am the wanderer

  with a lemma in every port.
\end{verse}
\vspace{\baselineskip}
\emph{---Marion D. Cohen}\footnote
{From her collection of math poetry
\emph{Crossing the Equal Sign} 
at http://mathwoman.com.
Originally published in the April 1999 
\emph{American Mathematical Monthly.}}

\chapter*{Further Reading --- Books / Plays / Musicals on Related Topics}
\addcontentsline{toc}{chapter}{Further Reading}
\markright{Further Reading}

Instead of providing an infinitely long bibliography, I decided to concentrate
mostly on \textbf{recent} books that caught my eye.
\begin{itemize}
\item 
Stephen Wolfram, \emph{A New Kind of Science,} Wolfram Media, 2002.
[A book concerned with many of the topics discussed here,
but with quite a different point of view.]
\item 
Joshua Rosenbloom and Joanne Sydney Lesser,
\emph{Fermat's Last Tango. A New Musical,}
York Theatre Company, New York City, 2000.
CD: Original Cast Records OC-6010, 2001.
DVD: Clay Mathematics Institute, 2001.
[A humorous and playful presentation of one mathematician's obsession with math.]
\item
David Foster Wallace, \emph{Everything and More: A Compact History of Infinity},
Norton, 2003. [An American writer's take on mathematics.]
\item
Dietmar Dath, \emph{H\"ohenrausch.\ Die Mathematik des XX.\ Jahrhunderts in zwanzig Gehirnen}, 
Eichborn, 2003.
[A German writer's reaction to modern math and modern mathematicians.]
\item
John L. Casti, \emph{The One True Platonic Heaven.\ A Scientific Fiction on the Limits of
Knowledge,} Joseph Henry Press, 2003.
[G\"odel, Einstein and von Neumann at the Princeton Institute for Advanced Study.] 
\item
Apostolos Doxiadis, \emph{Incompleteness, A Play and a Theorem,}
http://www.apostolosdoxiadis.com.
[A play about G\"odel from the author of \emph{Uncle Petros and Goldbach's Conjecture}.]
\item
Mary Terall, \emph{The Man Who Flattened the Earth.\ Maupertuis and the Sciences in the
Enlightenment}, University of Chicago Press, 2002.
[Newton vs.\ Leibniz, the generation after; a portrait of an era.] 
\item
Isabelle Stengers, \emph{La Guerre des sciences aura-t-elle lieu?\ Scientifiction,}
Les Emp\^echeurs de penser en rond/Le Seuil, 2001.
[A play about Newton vs.\ Leibniz.]
\item
Carl Djerassi, 
David Pinner, \emph{Newton's Darkness.\ Two Dramatic Views,}
Imperial College Press, 2003.
[Two plays, one about Newton vs.\ Leibniz.]
\item
Neal Stephenson, \emph{Quicksilver,} Morrow, 2003.
[A science fiction novel about Newton vs.\ Leibniz. Volume 1 of 3.]
\item
David Ruelle, \emph{Chance and Chaos,} 
Princeton Science Library,
Princeton University Press, 1993.
[A physicist's take on randomness.]
\item
Cristian S. Calude, \emph{Information and Randomness,} 
Springer-Verlag, 2002.
[A mathematician's take on randomness.]
\item
James D. Watson, Andrew Berry, \emph{DNA: The Secret of Life,} Knopf, 2003.
[The role of information in biology.]
\item
Tom Siegfried, \emph{The Bit and the Pendulum.\ From Quantum Computing to M Theory---The
New Physics of Information,} Wiley, 2000.
[The role of information in physics.]
\item
Tor N{\o}rretranders, \emph{The User Illusion.\ Cutting Consciousness Down to Size}, Viking, 1998.
[Information theory and the mind.]
\item
Hans Christian von Baeyer, \emph{Information: The New Language of Science}, 
Weidenfeld \& Nicholson, 2003. 
\item
Marcus du Sautoy, \emph{The Music of the Primes: Searching to Solve the Greatest Mystery in
Mathematics}, HarperCollins, 2003. [Of the recent crop of books on the Reimann
hypothesis, the only one that considers the possibility that G\"odel incompleteness
might apply.]
\item
Douglas S. Robertson, \emph{The New Renaissance: Computers and the Next Level of Civilization},
Oxford University Press, 1998.
\item
Douglas S. Robertson, \emph{Phase Change: The Computer Revolution in Science and Mathematics},
Oxford University Press, 2003. 
 
[These two books by Robertson discuss the 
revolutionary social transformations provoked by improvements in information transmission technology.]
\item
John Maynard Smith and E\"ors Szathm\'ary,
\emph{The Major Transitions in Evolution},
Oxford University Press, 1998.
\item
John Maynard Smith and E\"ors Szathm\'ary,
\emph{The Origins of Life: From the Birth of Life to the Origin of Language},
Oxford University Press, 1999.
 
[These two books by Maynard Smith and 
Szathm\'ary discuss evolutionary progress in terms of radical improvements 
in the representation of biological information.]
\item
John D. Barrow, Paul C. W. Davies, Charles L. Harper, Jr.,
\emph{Science and Ultimate Reality: Quantum Theory, Cosmology, and Complexity,}
Cambridge University Press, 2004.
[Essays in honor of John Wheeler.] 
\item
Gregory J. Chaitin, \emph{Conversations with a Mathematician,}
Springer-Verlag, 2002.
[Previous book by the author of this book.]
\item
Robert Wright,
\emph{Three Scientists and Their Gods: Looking for Meaning in an Age of Information,}
HarperCollins, 1989.
[On Fredkin's thesis that the universe is a computer.]
\item
Jonathan M. Borwein, David H. Bailey,
\emph{Mathematics by Experiment: Plausible Reasoning in the 21st Century,}
A. K. Peters, 2004.
[How to discover new mathematics. Volume 1 of 2.]
\item
Vladimir Tasi\'c, \emph{Mathematics and the Roots of Postmodern Thought,}
Oxford University Press, 2001.
[Where I learnt about Borel's know-it-all number.]
\item
Thomas Tymoczko, \emph{New Directions in the Philosophy of Mathematics,}
Princeton University Press, 1998.
\item
Newton C. A. da Costa, Steven French,
\emph{Science and Partial Truth,}
Oxford University Press,  2003.
\item
Eric B. Baum, \emph{What is Thought?,} MIT Press, 2004.
\end{itemize}

\end{document}